\documentclass{memo-l}

\usepackage[utf8]{inputenc}
\usepackage{amsmath, amssymb, bbm}
\usepackage{amsthm}
\usepackage{mathtools}
\usepackage{bbm}
\usepackage{blkarray}
\usepackage[matrix,arrow,curve]{xy}
\usepackage{color}
\usepackage{xcolor}
\usepackage{graphicx}
\definecolor{darkblue}{rgb}{0,0,0.6}
\usepackage[ocgcolorlinks,colorlinks=true, citecolor=darkblue, filecolor=darkblue, linkcolor=darkblue, urlcolor=darkblue]{hyperref}
\usepackage[capitalize,noabbrev]{cleveref}
\usepackage{comment}
\usepackage{adjustbox}
\usepackage{enumerate}

\usepackage{tikz}
\usetikzlibrary{
	cd,
	calc,
	positioning,
	arrows,
	decorations.pathreplacing,
	decorations.markings,
}

\newcounter{commentcounter}

\newcommand{\ignore}[1]{}
\renewcommand{\epsilon}{\varepsilon}
\renewcommand{\phi}{\varphi}
\newcommand{\bbZ}{\mathbbm{Z}}
\newcommand{\Z}{\mathbbm{Z}}

\newcommand{\R}{\mathbbm{R}}
\newcommand{\CP}{\mathbb{CP}}
\newcommand{\RP}{\mathbb{RP}}

\newcommand{\bbN}{\mathbbm{N}}

\newcommand{\raiso}{\xrightarrow{\,\cong\, }}
\newcommand{\wt}{\widetilde}
\newcommand{\ol}{\overline}
\newcommand{\Sq}{\mathrm{Sq}}
\newcommand{\red}{\mathrm{red}}
\newcommand{\Ext}{\mathrm{Ext}}
\newcommand{\Tor}{\mathrm{Tor}}
\newcommand{\sm}{\setminus}
\newcommand{\pt}{\mathrm{pt}}
\newcommand{\ter}{\mathfrak{ter}}
\newcommand{\msec}{\mathfrak{sec}}
\newcommand{\pri}{\mathfrak{pri}}

\newcommand{\RC}{\mathcal{RC}}
\newcommand{\scong}{{\cong_s}}
\newcommand{\wh}{\widehat}
\newcommand{\ba}{\begin{array}}
	\newcommand{\ea}{\end{array}}
\newcommand{\ra}{\longrightarrow}
\newcommand{\imra}{\looparrowright}
\newcommand{\RPT}{\mathbb{RP}^2}

\newcommand{\e}{\epsilon}
\newcommand{\cH}{\mathcal{H}}

\DeclareMathOperator{\Aut}{Aut}

\DeclareMathOperator{\Hom}{Hom}
\DeclareMathOperator{\id}{Id}
\DeclareMathOperator{\im}{im}
\DeclareMathOperator{\coker}{coker}
\DeclareMathOperator{\Id}{Id}

\DeclareMathOperator{\Sesq}{Sesq}

\DeclareMathOperator{\Tate}{Tate}

\DeclareMathOperator{\Arf}{Arf}

\DeclareMathOperator{\pr}{pr}
\DeclareMathOperator{\ex}{ex}
\DeclareMathOperator{\tr}{tr}
\DeclareMathOperator{\fr}{fr}
\DeclareMathOperator{\ev}{ev}

\DeclareMathOperator{\Wh}{Wh}

\newcommand{\spin}{\mathrm{Spin}}
\newcommand{\pin}{\mathrm{Pin}}

\newtheorem{thm}[equation]{Theorem}
\newtheorem{pre-thm}[equation]{Pre-Theorem}
\newtheorem{theorem}[equation]{Theorem}
\newtheorem{prop}[equation]{Proposition}
\newtheorem{conj}[equation]{Conjecture}
\newtheorem{question}[equation]{Question}

\newtheorem{cor}[equation]{Corollary}
\newtheorem{lemma}[equation]{Lemma}
\newtheorem{proposition}[equation]{Proposition}
\newtheorem{sproperty}[equation]{Secondary Property}
\newtheorem{tproperty}[equation]{Tertiary Property}
\theoremstyle{definition}

\newtheorem{defi}[equation]{Definition}
\newtheorem{definition}[equation]{Definition}
\newtheorem{rem}[equation]{Remark}
\newtheorem{remark}[equation]{Remark}

\newtheorem{cond}[equation]{Condition}
\newtheorem*{claim}{Claim}

\crefname{lemma}{Lemma}{Lemmas}
\crefname{section}{Section}{Sections}
\crefname{definition}{Definition}{Definitions}
\crefname{defi}{Definition}{Definitions}
\crefname{prop}{Proposition}{Propositions}
\crefname{theorem}{Theorem}{Theorems}
\crefname{thm}{Theorem}{Theorems}
\crefname{cor}{Corollary}{Corollaries}

\newcommand{\tmfrac}[2]{\mbox{\large$\frac{#1}{#2}$}} 

\numberwithin{section}{chapter}
\numberwithin{equation}{chapter}
\makeindex

\begin{document}

\frontmatter

\title{Algebraic criteria for stable diffeomorphism of spin $4$-manifolds}

\author{Daniel Kasprowski}
\address{School of Mathematical Sciences, University of Southampton, \newline \indent Southampton SO17 1BJ, United Kingdom}
\email{d.kasprowski@soton.ac.uk}

\author{Mark Powell}
\address{School of Mathematica and Statistics, University of Glasgow, \newline \indent University Place, Glasgow, G12 8QQ, United Kingdom}
\email{mark.powell@glasgow.ac.uk}

\author{Peter Teichner}
\address{Max Planck Institut f\"{u}r Mathematik, Vivatsgasse 7, 53111 Bonn, \newline \indent Germany}
\email{teichner@mac.com}

\date{\today}

\def\subjclassname{\textup{2020} Mathematics Subject Classification}
\expandafter\let\csname subjclassname@1991\endcsname=\subjclassname
\subjclass{Primary 57K40}

\keywords{Stable diffeomorphism, 4-manifold}


\begin{abstract}
We study closed, connected, spin $4$-manifolds up to stabilisation by connected sums with copies of $S^2 \times S^2$.
For a fixed fundamental group, there are primary, secondary and tertiary obstructions, which together with the signature lead to a complete stable classification.  The primary obstruction exactly detects $\mathbb{CP}^2$-stable diffeomorphism and was previously related to algebraic invariants by Kreck and the authors.
	
In this article we formulate conjectural relationships of the secondary and tertiary obstructions with algebraic invariants: the secondary obstruction should be determined by the (stable) equivariant intersection form and the tertiary obstruction via a $\tau$-invariant recording intersection data between 2-spheres, with trivial algebraic self-intersection, and their Whitney discs.
	
We prove our conjectures for the following classes of fundamental groups: groups of cohomological dimension at most $3$, right-angled Artin groups, abelian groups, and finite groups with quaternion or abelian $2$-Sylow subgroups.

We apply our theory to give a complete algebraic stable classification of spin $4$-manifolds with fundamental group $\Z \times \Z/2$.
\end{abstract}

\maketitle

\tableofcontents


\mainmatter

\chapter{Introduction}

Two closed, connected, smooth $4$-manifolds are \emph{stably diffeomorphic} if they become diffeomorphic after taking connected sums with finitely many copies of $S^2 \times S^2$, where we allow different numbers of copies for the two manifolds. An analogous notion of \emph{stable homeomorphism} applies in the topological category.  Kreck~\cite[Theorem~C]{surgeryandduality} reduced these classification problems to computations of bordism groups, which we review in  \cref{section:review} in the spin case.  He showed that two spin 4-manifolds $N_1$ and $N_2$ with fundamental group $\pi$ are stably diffeomorphic if and only there are maps $N_i \to B\pi$ and choices of spin structures such that the $N_i$ represent equivalent elements of the bordism group $\Omega_4^{\spin}(B\pi)$.  While this is extremely useful, explicitly calculating invariants of spin bordism groups is highly geometric, and usually intractable.

In this article we introduce new \emph{algebraic} invariants for spin $4$-manifolds.
The invariants are independent of category (smooth or topological) as one might expect in the stable setting. We shall focus on the smooth category, and explain why our classification also applies in the topological category in \cref{sec:topological-case}.

\begin{thm}\label{thm:1.1}
Fix a group $\pi$ with finite 3-dimensional classifying space, or an abelian group, a right angled Artin group, or a finite group $\pi$ with abelian or quaternion $2$-Sylow subgroup.

One can decide whether two closed, connected, spin, smooth 4-manifolds $N_1$ and $N_2$ are stably diffeomorphic over their fundamental group $\pi$ by computing, for all possible spin structures $\alpha_i$ on $N_i$, the following algebraic invariants on a 1-skeleton sum $M := (N_1,\alpha_1) \#_1 (N_2,-\alpha_2)$:
\begin{enumerate}
\setcounter{enumi}{-1}
\item\label{item:prethm0} the signature $\sigma(M) = \sigma(N_1) -\sigma(N_2)$ of $M$, an integer;
\item\label{item:prethm1} the extension class of the stable $\Z\pi$-module $\pi_2(M)$, a  `linear' invariant;
\item\label{item:prethm2} the stable intersection form $\lambda_M$ on $\pi_2(M)$, a $\Z\pi$-valued `quadratic' form;
\item\label{item:prethm3} the Kervaire-Milnor invariant $\tau_M$, a $\Z/2$-valued `cubical' refinement of $\lambda_M$.
\end{enumerate}
\end{thm}

There are handle decompositions of $N_i$ with diffeomorphic 1-skeleta (the union of the $0$- and $1$-handles). To obtain a 1-skeleton sum $(N_1,\alpha_1) \#_1 (N_2,-\alpha_2)$, remove the 1-skeleta, and identify the boundaries by a spin structure preserving diffeomorphism.  The stable diffeomorphism class turns out not to depend on the choices.

The invariants in \cref{thm:1.1} precisely detect the steps in an obstruction theory to decide whether $M$ is trivial in the spin bordism group $\Omega_4^{\spin}(B\pi)$.  While the signature is an invariant of each $N_i$, it remains an open problem to formulate the other invariants solely in terms of the $N_i$, in such a way that they detect stable diffeomorphism.  Our solution is to use the 1-skeleton sum described above.
We shall give precise definitions of the invariants in points \eqref{item:prethm1}, \eqref{item:prethm2}, and \eqref{item:prethm3} in \cref{section:alg-pri-intro,section:alg-sec-intro,section:alg-ter-intro}.

The linear invariant in \eqref{item:prethm1} exactly detects $\mathbb{CP}^2$-stable diffeomorphism of closed, connected $4$-manifolds. This was the focus of \cite{KPT18}, so in this article we shall focus on the other invariants.

The algebraic obstructions that we will discuss are generalisations of the ones studied in \cite{KLPT15} by the present authors together with Land. There the fundamental group was that of a closed, oriented, aspherical $3$-manifold, and we obtained our results by computing the obstructions on explicit models for the stable diffeomorphism classes.  The current investigation proceeds more abstractly, as suggested for the secondary obstruction in the third author's thesis~\cite{teichnerthesis}.

\section{The stable classification for fundamental group \texorpdfstring{$\Z\times \Z/2$}{ZxZ/2}}

As a concrete application, we provide an intrinsic stable classification for spin 4-manifolds with fundamental group $\Z \times \Z/2$.
We will recall in \cref{lem:Omega-4-ZxZ2} that $\Omega_4^{\spin}(B(\Z\times\Z/2))$ is isomorphic to $16\Z\oplus\Omega_2^{\pin^-}\cong 16\Z\oplus\Z/8$, where the isomorphism is given by the signature and the Brown invariant~\cite{Brown-invariant} of a $\Z/4$-valued quadratic enhancement of the $\Z/2$-valued intersection form of a surface, determined by the $\pin^-$ structure. By \cref{lem:action}, changing the spin structure acts on the $\Z/8$ summand by $-1$. Thus by Kreck's theorem, smooth, closed, connected, spin $4$-manifolds with fundamental group $\Z\times\Z/2$ are stably classified by their signature and the Brown invariant in $(\Z/8)/\pm1$. In particular, there are 5 stable diffeomorphism classes with zero signature that have Brown invariant $0,\pm 1,\pm 2,\pm 3$ and $4$. We apply our classification results to obtain an intrinsic description of the Brown invariant in this case.
Let $\pi:=\Z\times\Z/2=\langle T,t\mid [T,t],T^2\rangle$.

\begin{thm}
	\label{thm:ZxZ2-intro}
	For a smooth, closed, connected, spin $4$-manifold with fundamental group $\pi$, the Brown invariant $\beta(M)\in \Z/8$ is detected as follows:
	\begin{enumerate}
		\item $\beta(M)=\pm 1$ if and only if $\pi_2(M)$ is stably free and $\lambda_M$ is stably isometric to a form induced up to $\Z\pi$ from a form over $\Z$;
		\item $\beta(M)=\pm 3$ if and only if $\pi_2(M)$ is stably free but $\lambda_M$ is not stably induced from a form over $\Z$;
		\item $\beta(M)=\pm 2$ if and only if $\pi_2(M)$ is not stably free and $\lambda_M$ is odd;
		\item $\beta(M)\in\{0,4\}$ if and only if $\pi_2(M)$ is not stably free and $\lambda_M$ is even.
	\end{enumerate}
	In the last case, $\pi_2(M)$ is stably isomorphic to $I\oplus I^*$, where $I:=I\pi$ is the augmentation ideal, and there exists such a stable isomorphism for which the restriction of $\lambda_M$ to $I^*$ is trivial. Let $x\in\pi_2(M)$ be the image of $\phi\in I^*$ under such an isomorphism, where $\phi\in I^*$ is determined by $\phi(1-t)=1+T$ and $\phi(1-T)=0$. Then
	\begin{enumerate}
		\item [(5)] $\beta(M)=0$ if and only if $\tau_M(x)=0$, where $\tau_M$ is the Kervaire--Milnor invariant.
	\end{enumerate}
\end{thm}

\begin{rem}
  If $\pi_2(M)$ is stably free, we can determine whether the first case holds by computing the signature $\sigma(M)$, then asking whether $\lambda_M$ becomes stably hyperbolic after subtracting $\sigma(M)$ copies of $\Z\pi \otimes_{\Z} E_8$. In particular, if $\sigma(M)=0$, then the first case is equivalent to $\lambda_M$ being stably hyperbolic.
\end{rem}

The Kervaire--Milnor invariant is also detected by the stable equivariant intersection form for 2-dimensional \cite{HKT,KPT21} and 3-dimensional groups \cite{KLPT15,Hamb-Hild}. This leads us to the following question.
\begin{question}
    Is the difference between $\beta(M)=0$ and $\beta(M)=4$ also detected by the stable intersection form?
\end{question}

\section{Review of stable classification}\label{section:review}

As mentioned above, Kreck~\cite[Theorem~C]{surgeryandduality} showed  that two closed, connected, spin $4$-manifolds with fundamental group $\pi$ are stably diffeomorphic if and only if there are choices of spin structures and identifications of the fundamental groups with $\pi$, giving rise to equal elements in the bordism group $\Omega_4^{\spin}(\pi):=\Omega_4^{\spin}(B\pi)$, as we shall explain in \cref{section:AHSS}.  To understand this group of bordism classes of pairs $(M,c)$, where $M$ is a closed 4-manifold with spin structure and $c\colon M\to B\pi$ classifies the universal cover, we consider the Atiyah-Hirzebruch spectral sequence (AHSS) computing $\Omega_4^{\spin}(B\pi)$ in terms of $E_{p,q}^2 = H_p(B\pi;\Omega_q^{\spin})$. The AHSS gives rise to a filtration whose iterated graded quotients are
\[
\Z \cong \Omega_4^{\spin} \underbrace{\subseteq}_{E_{2,2}} F_2  \underbrace{\subseteq}_{E_{3,1}} F_3   \underbrace{\subseteq}_{E_{4,0}} \Omega_4^{\spin}(\pi).
\]
The first isomorphism is determined by the signature. More precisely, it is given by the signature divided by 16 in the smooth case and the signature  divided by 8 in the topological case; the closed topological $E_8$-manifold cannot be smoothed. This divisibility is the only difference between the stable classification of smooth and topological spin 4-manifolds and therefore we can ignore it in the sequel.
The signature extends to the entire group $\Omega_4^{\spin}(\pi)$ and so we reduce our study to $\widetilde\Omega_4^{\spin}(\pi)$, the kernel of the signature map, which is independent of category.
The AHSS then reduces to a shorter filtration
\[
E_{2,2} \underbrace{\subseteq}_{E_{3,1}} F   \underbrace{\subseteq}_{E_{4,0}} \widetilde\Omega_4^{\spin}(\pi),
\]
where the subgroup $F$ consists of bordism classes represented by signature zero 4-manifolds $M$ with spin structure such that $c\colon M\to (B\pi)^{(3)}$ lands in the 3-skeleton of the classifying space $B\pi$. Similarly, the smallest filtration term $E_{2,2}$ is represented by elements $(M,c)$ with $c\colon M\to (B\pi)^{(2)}$.
Since the $E^2_{p,q}$ term of the spectral sequence is $H_p(\pi;\Omega_q^{\spin})$, the $E^\infty_{p,q}$-terms are as follows:
\begin{itemize}
\item
$E_{2,2}:=E^\infty_{2,2}=H_2(\pi;\Z/2)/\im (d_2, d_3)$;
\item
$E_{3,1}:=E^\infty_{3,1}=H_3(\pi;\Z/2)/\im (d_2)$;
\item
$E_{4,0}:=E^\infty_{4,0}= \ker (d_2\colon  H_4(\pi;\Z)\to H_2(\pi;\Z/2))$.
\end{itemize}
Following \cite{teichnerthesis}, we obtain the primary invariant $\pri(M)=c_*[M] \in E_{4,0}$, the secondary invariant $\msec(M) \in E_{3,1}$ and the tertiary invariant $\ter(M) \in E_{2,2}$. The challenge is to recast these obstructions in terms of algebraic topological data of the 4-manifold. As a preliminary observation, the equivariant intersection form $\lambda_M$ on $\pi_2(M)$  changes by orthogonal sum with a hyperbolic form (on a free $\Z\pi$-module of rank 2) if we add $S^2 \times S^2$ to $M$.  So the signature $\sigma(M) \in \Z$ is a stable invariant, as is the isomorphism class of $\pi_2(M)$ up to stabilisation by free $\Z\pi$-modules, and the isometry class of $\lambda_M$ up to stabilisation with hyperbolic forms on free modules.

\section{Translating the bordism invariants into algebra}
Let us discuss the three obstructions in order of appearance.

\begin{enumerate}
   \item\label{item:prim-obstruction-description} The \emph{Primary Obstruction \cref{thm:prim-obstruction-thm}} \cite{KPT18} reinterprets the invariant $\pri(M)=c_*[M]$ as the extension class of a short exact sequence of $\Z\pi$-modules whose central module is stably isomorphic to $\pi_2(M)$. In loc.\ cit.\ we gave examples of various fundamental groups for which $\pri(M)$ is and is not determined by the stable isomorphism class of the $\Z\pi$-module $\pi_2(M)$.

  \item We say that a group $\pi$ has the \emph{Secondary Property} if $\msec(M)$ is detected by $\lambda_{M,s}$ for all spin $4$-manifolds $M$ with fundamental group $\pi$.  Here $\lambda_{M,s}$ is the restriction of $\lambda_M$ to a certain summand of $\pi_2(M)$. This summand arises as the image of a splitting $s$ of the extension in \eqref{item:prim-obstruction-description}, which exists if and only if $\pri(M)=0$. We shall show in \cref{section:alg-sec-intro} that if $\msec(M)$ vanishes then there are choices of a splitting for which $\lambda_M$ vanishes on the image.

\item We say that a group $\pi$ has the \emph{Tertiary Property} if $\ter(M)$ is detected by the Kervaire-Milnor invariant $\tau_{M,s}$ for all spin $4$-manifolds~$M$ with fundamental group~$\pi$.
 Here $\tau_{M,s}$ records intersection data between 2-spheres (with trivial algebraic self-intersection) and their Whitney discs. It is the Kervaire-Milnor invariant $\tau_M$ from \cite{Freedman-Quinn}, as corrected by \cite{Stong} and \cite{schneiderman-teichner-tau}, restricted to the image of a splitting in \eqref{item:prim-obstruction-description} on which $\lambda_M$ vanishes. This cubical intersection invariant is defined if and only if $\pri(M)=0=\msec(M)$, and it does not depend on the choice of splitting, as we will discuss in \cref{section:alg-ter-intro}.
\end{enumerate}

We will state these properties of a group in detail as \cref{secondary-property} and \cref{tertiary-property}. Then the following statement will be made precise.

\begin{thm}\label{pre-theorem=2-list-of-groups}
Consider a group $\pi$ with finite 3-dimensional classifying space $B\pi$, or an abelian group, a right angled Artin group, or a finite group $\pi$ with abelian or quaternion $2$-Sylow subgroup.
Then $\pi$ has the Secondary and Tertiary properties.  As a consequence, Theorem~\ref{thm:1.1} holds.
\end{thm}

Next we give more detailed descriptions of the obstructions and the statements of our theorems.

\subsection{The algebraic primary obstruction}\label{section:alg-pri-intro}

Let $M$ be a closed, connected, oriented $4$-manifold with fundamental group $\pi$.
There is a map $c \colon M \to B\pi$ inducing the identification on fundamental groups, that is well defined up to based homotopy.
The primary obstruction $\pri(M)$ corresponds to the edge homomorphism in the AHSS and is given by the image $c_*[M]$ of the fundamental class $[M]\in H_4(M;\Z)$ in $H_4(B\pi;\Z)=H_4(\pi;\Z)$.
We will build on the following result.

\begin{thm}[Primary obstruction theorem {\cite[Theorem~1.8]{KPT18}}]
	\label{thm:prim-obstruction-thm}
Let $K$ be a finite connected $2$-complex with fundamental group $\pi$. There is an isomorphism
	\[
	\Ext^1_{\Z\pi}(H^2(K;\Z\pi),\pi_2(K)) \cong H_4(\pi;\Z)
	\]
	mapping $\pri(M)=c_*[M]$ to an extension
	\[
0 \ra \pi_2(K) \ra (\Z\pi)^r \oplus \pi_2(M)  \ra H^2(K;\Z\pi) \ra 0.
	\]
\end{thm}

Recall that a finite presentation of $\pi$ gives a 2-complex $K$ as above by using a single 0-cell, one 1-cell for each generator and one 2-cell for each relation. Using Tietze transformations, we showed that the choice of presentation disappears when considering the extension group above. To make this precise, we need the following notion and the next lemma.

 \begin{definition}
 	\label{def:stableiso}
 	We say that two $R$-modules $P$ and $Q$ are \emph{stably isomorphic}, and write $P \scong Q$, if there exist non-negative integers $p$ and $q$ such that $P \oplus R^p \cong Q \oplus R^q$.
 \end{definition}

\begin{lemma}\cite[(40)]{lms197}
\label{lem:stable}
Let $K_i$ be finite connected  $2$-complexes with fundamental group $\pi$. Then there exist $k_i \in \bbN_0$ such that $K_1\vee \bigvee_{r=1}^{k_1} S^2\simeq K_2\vee\bigvee_{r=1}^{k_2} S^2$. In particular, $\pi_2(K_1)\scong \pi_2(K_2)$ and $H^2(K_1;\bbZ\pi)\scong H^2(K_2;\bbZ\pi)$.
\end{lemma}

In \cite[Lemma~5.11]{KPT18} we checked that wedge sum with $S^2$ does not change the extension group in \ref{thm:prim-obstruction-thm} and  that a homotopy self-equivalence of $K$ inducing the identity on $\pi_1(K)$ determines the identity map on this extension group~\cite[Lemma~5.12]{KPT18}.
We concluded that the image of $\pri(M)$ in this extension group is a well-defined algebraic invariant of $M$ in a group that depends only on the group~$\pi$.

\subsection{The algebraic secondary obstruction}\label{section:alg-sec-intro}

Consider the equivariant intersection form
\[
\lambda_M \colon \pi_2(M) \times \pi_2(M) \to \Z\pi
\]
of a closed, connected, oriented $4$-manifold $M$ with fundamental group $\pi$. It is given by either counting geometric intersections (with signs and group elements) between transverse 2-spheres in $M$ or by identifying $\pi_2(M) \cong H_2(M;\Z\pi)$ via the Hurewicz homomorphism and setting
\[
\lambda_M(x,y) := \langle PD^{-1}(y), x \rangle \quad \text{for } x,y \in H_2(M;\Z\pi),
\]
where $PD \colon H^2(M;\Z\pi) \to H_2(M;\Z\pi)$ is the Poincar\'{e} duality isomorphism.
The intersection form $\lambda_M$ is \emph{sesquilinear}, meaning it is additive in each variable and $\lambda_M(ax, by)= a\lambda_M(x,y) \overline{b}$ holds for all $a,b\in\Z\pi$, where the involution $a\mapsto \overline{a}$ on the group ring $\Z\pi$ is determined by $\overline{g} := g^{-1}$ for $g\in \pi$.
Note that $\lambda_M$ is frequently singular. In fact, its kernel and cokernel are determined by the fundamental group $\pi$ because by the universal coefficient spectral sequence, its adjoint fits into the exact sequence
\begin{equation}\tag{UCSS} \label{eq:UCSS}
0 \to H^2(\pi;\Z\pi) \to \pi_2(M) \overset{\lambda_M^{\operatorname{ad}}}{\ra} \Hom_{\Z\pi}(\pi_2(M), \Z\pi) \to H^3(\pi;\Z\pi)\to 0.
\end{equation}
The intersection form $\lambda_M$ is also hermitian in the following sense.
The dual $\lambda^*$ of a sesquilinear form $\lambda$ is given by
\[
\lambda^*(x,y) := \overline{\lambda(y,x)}.
\]
and $\lambda$ is called \emph{hermitian} if $\lambda^*=\lambda$. We say that $\lambda$ is \emph{even} if $\lambda=q + q^*$ for some sesquilinear form $q$ on $\pi_2(M)$.

\begin{remark}\label{rem:q}
If $\lambda_M$ is even, there is a {\it quadratic refinement} $\widehat{q} \colon \pi_2(M) \to \Z\pi/\{a - \overline{a} \mid a\in \Z\pi\}$ of $\lambda_M$  defined as  $\widehat q(x):=[q(x,x)]$ for some choice of sesquilinear form $q$ with $\lambda_M=q+q^*$. This means that, for $x,y\in \pi_2(M)$, we have
\begin{enumerate}
\item $\lambda_M(x,x) = \widehat q(x) + \overline{\widehat{q}(x)}\in\Z\pi$ and
\item $ \widehat{q}(x+y)= \widehat{q}(x) + \widehat{q}(y) + \lambda_M(x,y) \in \Z\pi/\{a - \overline{a} \mid a\in \Z\pi\}.$
\end{enumerate}
\end{remark}

On the other hand, it is not necessary for $\lambda_M$ to be even in order for it to have a quadratic refinement.
Every class $x\in \pi_2(M)$ can be represented by a generic immersion $f\colon S^2\looparrowright M$. If $\widetilde{M}$ is spin, one can arrange for $f$ to have trivial normal bundle as follows. Since $w_2(x)=0$, the Euler number $e(f)$ of the normal bundle to $f$ is even and hence can be changed to zero by a non-regular homotopy of $f$, performing $e(f)/2$ cusp homotopies.
As a consequence, the intersection form $\lambda_M$  has a quadratic refinement
\[
\mu_M \colon \pi_2(M) \to \Z\pi/\{a - \overline{a} \mid a\in \Z\pi\}
\]
given by counting self-intersections (with group elements and signs) of a generic immersion $f \colon S^2 \looparrowright M$ with $e(f)=0$.
The expression $\mu_M(x)+\overline{\mu_M(x)}$ can be lifted canonically to $\Z\pi$ and then equals $\lambda_M(x,x)$.  To see this, using that $e(f)=0$ choose a nowhere vanishing section of $\nu f$, and let $f'$ denote a push-off of $f$ along this section.  To compute $\lambda_M(x,x)$ we may count intersections between $f$ and $f'$.
Each double point of $f$ with double point loop $\pm g$ contributes $\pm(g+ g^{-1})$ to $\lambda_M(x,x)$, and there are no other contributions because we used a nowhere vanishing section to define $f'$. See \cite[Chapter~5]{Wall} for more details.

It is not hard to show algebraically that $\mu_M$, if it exists, is determined by $\lambda_M$ via the equation $\lambda_M(x,x) = \mu_M(x)+\overline{\mu_M(x)}$.
Thus we shall not keep $\mu_M$ in the notation, but will remember that $\lambda_M$ is a \emph{weakly even hermitian} form, meaning by definition that a quadratic refinement exists.

\begin{remark}
The use of cusp homotopies to fix the Euler number of the normal bundle amounts to a normalisation of $\mu_M(x)$ at $1\in \pi$. It is important to note that the above normalisation does not in general not lead to an obstruction to representing the homotopy class $x\in \pi_2(M)$ by an embedding $f\colon S^2 \hookrightarrow M$. However, $\mu_M(x)=0$ does hold if $x$ is represented by a framed embedding $f$ where $e(f)=0$. For example, this is important if one plans to do surgery on $f$ to simplify $\pi_2(M)$.

To get an embedding obstruction also in the non-framed case, one can use the following alternative normalisation at $1 \in \pi$. Using cusp homotopies, change $f\colon S^2\looparrowright M$ to arrange that
the self-intersection number of $f$ at $1$ is trivial.  This can be done even if the universal covering of $M$ is not spin. Taking self-intersections gives an element $\mu'_M(x)\in \Z\pi/\{a-\ol{a} \mid a \in \Z\pi\}$, which
again only depends on the homotopy class $x\in\pi_2(M)$ and is an obstruction for representing $x$ by an embedding. In the setting of this paper, our 4-manifolds are spin, thus so are their universal covers.
For $\wt{M}$ spin we have the relation $\mu'_M(x) + \tmfrac{1}{2} \lambda_M(x,x)_1 = \mu_M(x)$, because by construction $\mu'(x)$ vanishes at $1\in\pi$ and $\mu_M$ is a quadratic refinement. Here, $\lambda_M(x,x)_1$ denotes the coefficient of $\lambda_M(x,x)$ at $1\in\pi$.
In particular when $\lambda_{M}(x,x)_1=0$, the two conventions agree.
\end{remark}

In general, weakly even forms are not even. If $\pi_2(M)$ happens to be a free $\Z\pi$-module on an ordered basis $\{e_i\}$ then weakly even forms~$\lambda$ are even, as can be seen by lifting $\mu(e_i)$ to $q(e_i,e_i)\in\Z\pi$ and setting $q(e_i,e_j):=\lambda(e_i,e_j), q(e_j,e_i)=0$ for $i<j$. But as we shall see, for many spin 4-manifolds $M$ the intersection form~$\lambda_M$ is not even.
To make this precise, consider a closed, connected spin 4-manifold $M$ with fundamental group $\pi$ and $\pri (M)=0$. In the notation of \cref{thm:prim-obstruction-thm}, where we chose a 2-complex $K$, also choose  a splitting
\[
s=(s_1,s_2)\colon H_K := H^2(K;\Z\pi)\to (\Z\pi)^r \oplus \pi_2(M)
\]
of the short exact sequence $0 \ra \pi_2(K) \ra (\Z\pi)^r \oplus \pi_2(M)  \ra H^2(K;\Z\pi) \ra 0$. Consider the sesquilinear form $\lambda_{M,s} \colon H_K \times H_K \to \Z\pi$ on $H_K$ induced  by $s_2$,  $\lambda_M(s_2(-),s_2(-))$.
We write $\Sesq(H)$ for the group of $\Z\pi$-sesquilinear forms on the $\Z\pi$-module $H$. Sending a form to its hermitian conjugate induces a $\Z/2 = \langle T \rangle$ action on $\Sesq(H)$, and we consider the corresponding Tate group
\begin{align*}
  \widehat{H}^0(\Sesq(H)) := \widehat{H}^0(\Z/2;\Sesq(H)) &= \ker(1-T)/\im(1+T) \\  &=  \{\text{Hermitian forms}\}/\{\text{even forms}\}.
\end{align*}
In Definition~\ref{defn:map-A} we will construct, purely algebraically, a homomorphism
\[
A_K \colon H_3(\pi;\Z/2) \ra \widehat{H}^0(\Sesq(H_K)).
\]
For every map $\phi\colon K\to K'$ inducing the given identification with $\pi$ on fundamental groups, the diagram
\[\xymatrix{H_3(\pi;\Z/2)\ar[r]^-{A_K}\ar[dr]_{A_{K'}}&\widehat{H}^0(\Sesq(H_K))\ar[d]^{\phi_*}\\
&\widehat{H}^0(\Sesq(H_{K'}))}\]
commutes by \cref{lem:ker-A-well-defined}. We obtain a group \[\Tate(\pi):= H_3(\pi;\Z/2)/\ker(A_K) \cong \im A_{K},\] which is independent of the choice of~$K$. We write the projection as
\[
A_\pi \colon H_3(\pi;\Z/2) \twoheadrightarrow \Tate(\pi).
\]
Although it will not play a role for this article, note that $\im A_{K}$ lies in the subgroup of weakly even forms on $H_K$ by \cref{lem:weaklyeven}. One might hope to identify $\Tate(\pi)$ with a particular class of such weakly even forms on $H_K$.

 We will show in \cref{subsection:proof-thm-mainsec} that for every spin 4-manifold~$M$ with fundamental group~$\pi$ and $\pri (M)=0$,
a splitting $s\colon H_K\to (\Z\pi)^r\oplus\pi_2(M)$ gives rise (after making some choices) to a map $f \colon M\to (B\pi)^{(3)}$. This determines an element $\msec(M,f)\in H_3(\pi;\Z/2)$, lifting the element $\msec(M) \in H_3(\pi;\Z/2)/\im(d_2)$.
We will show that the two elements $\lambda_{M,s}$ and $A_\pi(\msec(M,f))$ agree in $\Tate(\pi)\subseteq \widehat{H}^0(\Sesq(H_K))$.

\begin{sproperty}
\label{secondary-property}
We say that a group $\pi$ has the \emph{Secondary Property} if for all closed, connected spin 4-manifolds $M$ with fundamental group $\pi$ and $\pri (M)=0$,  $\lambda_{M,s} \in \Tate(\pi)$ is independent of the splitting $s$ and $\msec(M)$ is determined by $\lambda_{M,s} \in \Tate(\pi)$.
\end{sproperty}

The Secondary Property of a group $\pi$ follows from the purely algebraic assertion that the differential $d_2$ in the AHSS makes the following sequence exact:
\begin{equation} \tag{$\msec_\pi$} \label{eq:sec}
H_5(\pi;\Z) \overset{d_2}{\ra} H_3(\pi;\Z/2) \overset{A_\pi}{\ra} \Tate(\pi).
\end{equation}
In \cref{sec:ex,sec:abelian-groups} we will prove this exactness (and hence the Secondary Property) for the following classes of groups, proving the secondary part of Theorem~\ref{pre-theorem=2-list-of-groups}.

\begin{thm}
\label{thm:sec}
The \cref{secondary-property} holds for the following classes of groups.
\begin{itemize}
\item Groups $\pi$ with a finite $3$-dimensional model for $B\pi$.
\item Finitely generated abelian groups and right angled Artin groups.
\item Finite groups whose $2$-Sylow subgroups are abelian or quaternion.
\end{itemize}
\end{thm}

When $\pi$ has geometric dimension at most two, the Secondary Property holds since  $H_3(\pi;\Z/2)=0$.
If $\pi$ is the fundamental group of a closed, oriented, aspherical $3$-manifold then $H_5(\pi;\Z)=0$ and $H_K \cong_s I\pi$, the augmentation ideal in~$\Z\pi$. We showed in \cite{KLPT15} that $A_\pi$ is an isomorphism between two groups of order 2.

The case of finite groups with quaternion $2$-Sylow subgroups was proven in \cite[Theorem~6.4.1]{teichnerthesis}; see also our \cref{sec:gen-quaternion-groups}. Emboldened by our success, we make the following conjecture cf.~\cite[Conjecture~B]{teichnerthesis}.

\begin{conj}\label{secondary-conjecture}
The sequence \ref{eq:sec} is exact for every finitely presented group $\pi$.
\end{conj}

\subsection{The algebraic tertiary obstruction}\label{section:alg-ter-intro}
We start by giving a brief description of the \emph{Kervaire--Milnor} invariant $\tau$ that appears in the Tertiary Property.
Recall that in the smooth category, a generic map of a surface to a $4$-manifold looks locally like a linear subspace $\R^2 \subset \R^4$ or like a transverse double point $\R^2 \times \{0\} \cup \{0\} \times \R^2 \subset \R^4$. We know that, stably, every spin manifold with signature zero has a smooth structure, so we may assume that all our surfaces are smoothly embedded with respect to such a structure.

Let $f \colon S^2 \imra M$ be a generically immersed 2-sphere with vanishing algebraic self-intersection, that is~$\mu_M(f)=0$. Then all the double points of~$f$ can be paired by a family of generically immersed Whitney discs $\{W_i\}$. They can be chosen to be disjointly embedded, framed, and to intersect $f$ transversely, by boundary twisting to fix the framing and then pushing any intersections down to $f$.
We define
\[\tau(f,\{W_i\}):=\sum_i |\mathring{W}_i\cap f| \mod 2.\]
This count depends a priori on the choice of Whitney discs, so one has to be careful.

Recall that for a closed, connected 4-manifold $M$, since the $\Z/2$ intersection form is nonsingular the homomorphism
\[
H_2(M;\Z/2) \ra\Z/2, \quad x \mapsto x \cdot x := \langle PD^{-1}(x), x \rangle
\]
can be written as $x\cdot x = c\cdot x$ for a unique \emph{characteristic} element $c\in H_2(M;\Z/2)$. This motivates the terminology in the next definition.

\begin{definition}\label{defn:char}
 For $\alpha \in H_2(M;\Z/2)$, if $x\cdot x = \alpha \cdot x$ for all \emph{spherical} $x$, meaning those $x$ represented by a map $S^2 \to M$, we call $\alpha$ \emph{$S^2$-characteristic}.

 Similarly we call $\alpha$ \emph{$\RPT$-characteristic} if the same equation holds for all classes $x$ represented by a map $\RP^2\to M$. Note that $\RPT$-characteristic implies $S^2$-characteristic.

 We say that an element $\beta \in \pi_2(M)$ is $S^2$- or $\RPT$-characteristic if its image under the modulo 2 Hurewicz map $\pi_2(M) \to H_2(M;\Z/2)$ is $S^2$- or $\RPT$-characteristic respectively.
\end{definition}

It was shown in \cite{Stong} and \cite{schneiderman-teichner-tau} that
$\tau(f,\{W_i\})$ is independent of the choice of $\{W_i\}$
exactly when $f$ is $\RPT$-characteristic. Moreover, the resulting invariant \[\tau([f]):=\tau(f,\{W_i\})\in \Z/2\] only depends on the homotopy class $[f]\in \pi_2(M)$. It was called the Kervaire-Milnor invariant by Freedman-Quinn, who gave a slightly incomplete account of it in \cite[Section~10]{Freedman-Quinn}. This was later resolved by Stong~\cite{Stong}.

Let us return to our obstruction theory for closed, connected spin 4-manifolds $M$ with $\pri(M)=0=\msec(M)$. For the algebraic secondary invariant, we chose a 2-complex $K$ and a stable splitting $s \colon H_K\to (\Z\pi)^r \oplus \pi_2(M)$
of the short exact sequence from \cref{thm:prim-obstruction-thm} and defined the hermitian form $\lambda_{M,s}$ on $H_K$ as the composition of $\lambda_M$ with $s_2 \colon H_K \to \pi_2(M)$. Up to even forms, this did not depend $K$ nor on the splitting $s$ and there are splittings $s$ for which $\lambda_{M,s}$ actually vanishes; see \cref{remark-section-as-required}.

We want to define $\tau$ on elements in the image of such splittings $s$.
Start with an element $x\in H^2(\pi;\Z/2)$, restrict it using $K \subseteq B\pi$ to $H^2(K;\Z/2)$ and lift that to $\wt{x}\in H_K=H^2(K;\Z\pi)$, using that $H^3(K;-)=0$, so the Bockstein potentially obstructing this lifting vanishes. Given a splitting $s$ with $\lambda_{M,s}\equiv 0$, let $y := s_2(\wt{x}) \in \pi_2(M)$, and  then compute $\tau_M(y)\in\Z/2$. This is well defined if
$y$ is $\RPT$-characteristic.

We will show in \cref{def:tauM} that the above construction determines a well-defined map
\[
\tau_{M,s} \colon h^2(\pi) \to \Z/2
\]
on the set
\[h^2(\pi):= \ker (Sq^2 \colon H^2(\pi;\Z/2) \to H^4(\pi;\Z/2)).\]
Here $h^2(\pi)$ is dual to the $E^3_{2,2}$-term in our AHSS, since
\[
Sq^2=(d^2_{4,1})^* \colon H^2(\pi;\Z/2) \to H^4(\pi;\Z/2)
\]
is the dual of the $d^2_{4,1}$ differential.  The duals of the other relevant differentials in the AHSS for $\Omega_4^{\spin}(B\pi)$ are
\[
(d^2_{5,0})^* \colon H^3(\pi;\Z/2) \to \Hom_{\Z}(H_5(\pi;\Z),\Z)  \quad \text{ and }
\]
\[
(d^3_{5,0})^* \colon h^2(\pi) = \ker(d^2_{4,1})^* \to \ker(d^2_{5,0})^*.
\]
There is a canonical isomorphism
\begin{align*}
\omega \colon E_{2,2} = H_2(\pi;\Z/2)/\im (d^2_{4,1}, d^3_{5,0}) &\xrightarrow{\cong} \Hom_{\Z/2}(H^2(\pi;\Z/2),\Z/2)/\im (d^2_{4,1}, d^3_{5,0})\\ &\xrightarrow{\cong} \Hom_{\Z/2}(\ker(d_{5,0}^3)^*,\Z/2).
\end{align*}
Since $\ker(d_{5,0}^3)^* \subseteq h^2(\pi)$, we can restrict $\tau_{M,s}$ to a map $\tau_{M,s}| \colon  \ker(d_{5,0}^3)^* \to \Z/2$.  We want to show that this coincides with the image of $\ter(M)$ under~$\omega$.  Again, we give the statement in the form of a property of a group.

\begin{tproperty}
\label{tertiary-property}
We say that $\pi$ has the \emph{Tertiary Property} if $\ter(M)\in E_{2,2}$ is sent to $\tau_{M,s}|$ via~$\omega$ for all closed, connected, spin 4-manifolds $M$ with fundamental group $\pi$ and $\pri(M) =0 = \msec(M)$.
\end{tproperty}

We remark that $\tau_{M,s}(0)=0$ because $0 \in h^2(\pi)$ leads to computing $\tau_M$ on $0 \in \pi_2(M)$, which is represented by a trivial (embedded) sphere. So in the case that $\ker(d_{5,0}^3)^* = 0$, the Tertiary Property automatically holds.
We are unable to show in general that $\tau_{M,s}|$ is a homomorphism.  This follows from the Tertiary Property, when it holds, since $\ter(M)$ maps to a homomorphism.

In \cref{lem:ter} we will show that the Tertiary Property holds whenever the splitting $s$ is induced by a map $f \colon M \to K$, in the sense that
\[s_2 = PD \circ f^* \colon H^2(K;\Z\pi) \to H^2(M;\Z\pi) \to H_2(M;\Z\pi) \cong \pi_2(M)\]
for $K$ a 2-complex with $f_* \colon \pi_1(M) \to \pi_1(K)$ an isomorphism.
In \cref{thm:ter-alg-condition} we give a condition under which every splitting is realised in this way, and we use this to deduce the following theorem.

\begin{thm}
\label{thm:ter}
The Tertiary Property holds for the following classes of groups.
\begin{itemize}
\item Groups $\pi$ with a finite $3$-dimensional model for $B\pi$.
\item Finitely generated abelian groups and right angled Artin groups.
\item Finite groups whose 2-Sylow subgroups are abelian or quaternion.
\end{itemize}
\end{thm}

Our verification of \cref{thm:ter} in the case of abelian groups relies on computations of Whitehead's $\Gamma$-groups~\cite{whitehead} made by the first and second authors with Ben Ruppik~\cite{KPR}.
In the special case that $\pi$ is the fundamental group of a closed, oriented, aspherical $3$-manifold, we proved in \cite{KLPT15} that $\pi$ has the \cref{tertiary-property}. To do this, we computed $\ter(M)$ and $\tau_{M,s}$ on sufficiently many concrete examples for which the primary and secondary obstructions vanish. We close this section with the following conjecture, analogous to \cref{secondary-property}.

\begin{conj}\label{tertiary-conjecture}
 Every group has the Tertiary Property.
\end{conj}

\section{Context for our work}

The stable classification of 4-manifolds was reduced to a bordism computation by Kreck in \cite{surgeryandduality}. $4$-manifolds with finite fundamental group were studied by Hambleton and Kreck in \cite{Ham-kreck-finite}, as well as in the PhD thesis of the third author~\cite{teichnerthesis}.  The case of geometrically $2$-dimensional groups was solved by Hambleton, Kreck and the last author in~\cite{HKT}.  Groups of cohomological dimension at most 3, and in particular right angled Artin groups with this property, were studied by Hambleton and Hildum~\cite{Hamb-Hild} in the case that the equivariant intersection form is even, which by our results in the current paper is equivalent to the secondary obstruction vanishing. Previous work on the stable classification question also includes work of Cavicchioli, Hegenbarth and Repov\v{s} \cite{CAR-95}, Spaggiari \cite{Spaggiari-03} and Davis \cite{Davis-05}.

The papers \cite{Ham-kreck-finite}, \cite{HKT} and \cite{Hamb-Hild} focussed on the normal 2-type -- roughly this means they looked at the entire intersection form of the $4$-manifold -- also obtaining unstable classification results. On the other hand, the algebraic invariants that we study are tailored to the stable question, so we are able to handle more fundamental groups with more easily computable obstructions.

\subsection{When stable homeomorphism implies homeomorphism}

By combining our results with results of Hambleton-Kreck, Khan and Crowley-Sixt, we obtain results on the unstable homeomorphism classification.  If we consider manifolds that are already sufficiently stabilised, in a sense to be made precise presently, then it turns out that stable homeomorphism implies homeomorphism.

Recall that a finitely presented group $\pi$ is \emph{polycyclic-by-finite} if it has a subnormal series where the quotients are either cyclic or finite. The
number of infinite cyclic quotients in such a subnormal series turns out to be an invariant of $\pi$, called the \emph{Hirsch number} $\mathfrak{h}(\pi)$. Define $\mathfrak{h}'(\pi)=1$ if $\pi$ is finite, and define $\mathfrak{h}'(\pi)=\mathfrak{h}(\pi)+3$ if $\pi$ is infinite. The following theorem is due to Hambleton and Kreck \cite[Theorem B]{hambleton-kreck93} for $\pi$ finite, and due to Crowley and Sixt \cite[Theorem 1.1]{crowley-sixt} for $\pi$ infinite.

\begin{theorem}[{\cite{hambleton-kreck93,crowley-sixt}}]
Let $M$ and $N$ be closed, connected $4$-manifolds with polycyclic-by-finite fundamental group $\pi$ such that  $\chi(N) +2k = \chi(M)$ for some $k\geq \mathfrak{h}'(\pi)$. If $M$ and $N$ are stably homeomorphic then $M$ is homeomorphic to~$N \# k(S^2 \times S^2)$.
\end{theorem}

Note that \cite[Theorem 1.1]{crowley-sixt} is only stated in the smooth category, but as remarked at the beginning of \cite[Section 2.1]{crowley-sixt}, the theorem also holds in the topological category since polycyclic-by-finite groups are good in the sense of Freedman~\cite[Section~2.9~and~Theorem~5.1A]{Freedman-Quinn}.

For $\pi$ virtually abelian and infinite, there is a similar statement to \cite[Theorem~1.1]{crowley-sixt} by Khan \cite[Corollary 2.4]{khan}, with a slightly better bound on the number of stabilisations needed.

\subsection*{Acknowledgements}

We are delighted to have the opportunity to thank Diarmuid Crowley, Jim Davis, Fabian Hebestreit, Matthias Kreck, Markus Land, Ian Hambleton, Henrik R\"uping, Mark Ullmann and Christoph Winges for many useful and interesting discussions on this work. We also thank the anonymous referee for carefully reading the paper and providing many helpful suggestions. 

The authors thank the Max Planck Institute for Mathematics and the Hausdorff Institute for Mathematics in Bonn for financial support and their excellent research environments.

The second author was partially supported by EPSRC New Investigator grant EP/T028335/2 and EPSRC New Horizons grant EP/V04821X/2.


\chapter{Background}

\section{Stable diffeomorphism, bordism groups, and spectral sequences}\label{section:AHSS}

Throughout the paper $M$ denotes a smooth, closed, compact, connected, oriented $4$-manifold with fundamental group $\pi_1(M)\cong \pi$.   The universal cover of $M$ will be denoted $\wt{M}$.  If we fix an identification of $\pi_1(M)$ with $\pi$, the identification determines a homotopy class of maps $c \colon M \to B\pi$ classifying $\wt{M}$.

Two smooth 4-manifolds $M$ and $N$ are called \emph{stably diffeomorphic} if there exist integers $m,n \in \mathbb{N}_0$ such that stabilising $M$ and $N$ with copies of $S^2 \times S^2$ yields diffeomorphic manifolds
\[M \# m (S^2 \times S^2) \cong N \# n (S^2 \times S^2).
\]
We require that the diffeomorphism respects orientations.
Note that, unlike elsewhere in the literature, we do not require that $m=n$.

\subsection{Bordism groups}
The starting point for our investigation, coming from Kreck's modified surgery, is that stable diffeomorphism can be understood in terms of bordism over the normal $1$-type.

\begin{definition}
Let $M$ be a closed oriented manifold of dimension $n$. A \emph{normal $1$-type} is a fibration over $BSO$, denoted by $\xi \colon B \to BSO$, through which the stable normal bundle $\nu_M \colon M \to BSO$ factors as:
\[\xymatrix@R=.5cm{ & B \ar[d]^\xi \\ M \ar@/^.7pc/[ur]^-{\wt{\nu}_M} \ar[r]_-{\nu_M} & BSO }\]
such that $\wt{\nu}_M$ is $2$-connected and $\xi$ is $2$-coconnected. A choice of a map $\wt{\nu}_M$ is called a \emph{normal $1$-smoothing} of $M$.
\end{definition}


The different normal $1$-types of a fixed $M$ are fibre homotopy equivalent to one another.

\begin{thm}[{\cite[Theorem C]{surgeryandduality}}]\label{thm:stablediffeoclasses}
Two closed $4$-dimensional manifolds with the same Euler
characteristic and the same normal $1$-type $\xi \colon B \to BSO$, admitting bordant normal $1$-smoothings, are diffeomorphic after connected sum of both with $r$ copies of $S^2\times S^2$ for some $r$.
\end{thm}

Here is a summary of the proof.
Given a bordism over the normal $1$-type, one can perform surgery below the middle dimension until the bordism becomes $1$-connected.  Excise $2$-spheres by tubing a neighbourhood to the boundary.  This makes the bordism into an $s$-cobordism, which is stably a product, at the expense of connect summing the boundary $4$-manifolds with copies of $S^2 \times S^2$.  A key point to check is that the normal smoothing data allows us to excise $2$-spheres with a framed normal bundle, resulting in connect sums with $S^2 \times S^2$ and not the twisted bundle $S^2 \wt{\times} S^2$.  The converse, that stably diffeomorphism manifolds admit bordant $1$-smoothings, can be found in Crowley-Sixt~\cite[Lemma~2.3(ii)]{crowley-sixt}.

By taking account of the different choices of normal $1$-smoothing, and allowing different numbers of $S^2\times S^2$ to be added to either side to allow the Euler characteristics of the initial manifolds to differ, we obtain the following.

\begin{theorem}\label{cor:stablediffeoclasses}
The stable diffeomorphism classes of $4$-manifolds with normal $1$-type $\xi$ are in one-to-one correspondence with $\Omega_4(\xi)/\Aut(\xi)$.
\end{theorem}

Next we need to understand the normal $1$-type.
The following lemma, determining the normal $1$-type of spin 4-manifolds, is well-known to the experts. We refer to  \cite[Lemma~3.5]{KLPT15} for a proof.

\begin{lemma}\label{lem:normal-spin-case}
Let $\pi$ be a finitely presented group. A normal $1$-type of a spin $4$-manifold with fundamental group $\pi$ is given by
\[\xymatrix{B\pi\times BSpin\ar[r]^-{\gamma \circ \pr_2}&BSO,}\]
where $\pr_2$ is the projection onto $BSpin$ and $\gamma$ is the canonical map $BSpin\to BSO$.
\end{lemma}

This article gives algebraic invariants that obstruct null-bordism in $\Omega_4(\xi)$ in the above two cases.

\subsection{Spectral sequences}
The main tool for understanding the bordism groups is the James spectral sequence.  However, as explained in \cite[Section~2]{KLPT15}, this is isomorphic to the Atiyah-Hirzebruch spectral sequence for $\Omega_4^{\spin}(B\pi)$.

Let $h_*$ be a generalised homology theory (for us spin bordism $\Omega^{\spin}_*$), and let $X$ be a CW complex.  The Atiyah-Hirzebruch spectral sequence arising from the homology Leray-Serre spectral sequence for the fibration $\pt \to X \to X$ has $E_2$ page given by
\[E^2_{p,q} := H_p(X;h_q(\pt)).\]
The $d^r$ differential has $(p,q)$-bidegree $(-r,r-1)$, and the sequence converges to $h_*(X)$.  We denote the differentials by
\[d^r_{p,q} \colon H_p(X;h_q(\pt)) \to H_{p-r}(X;h_{q+r-1}(\pt)).\]
Denote the filtration on the abutment of an Atiyah-Hirzebruch spectral sequence by
\[ 0 \subset F_{0,n} \subset F_{1,n-1} \subset \dots \subset F_{n-q,q} \subset \dots \subset F_{n,0} = h_n(X) = \Omega_n(\xi).\]
Recall that $F_{n-q,q}/F_{n-q-1,q+1} \cong E^{\infty}_{n-q,q}$.

Let $B:= B\pi \times BSpin$. Denote the restriction of $\pr_1 \colon B \to B\pi$ to the inverse image of the $p$-skeleton of $B\pi$ by $B|_p$, and let $\xi|_p \colon B|_p \to BSO$ be the restriction of $\xi$ to $B|_p$. An element of $\Omega_n(\xi)$ lies in $F_{p,n-p}$ if and only if it is in the image of the map $\Omega_n(\xi|_p) \to \Omega_n(\xi)$.  That is, if the element lies in the image of $\Omega_n^{\spin}(B\pi^{(p)})$.
This follows from the naturality of the spectral sequence applied to the map of fibrations induced by the inclusion of $B\pi^{(p)} \to B\pi$.

In our situation, being in $F_{p,4-p}$ means that a $4$-manifold $M$ together with its classifying map $c \colon M \to B\pi$ is spin
bordant to a $4$-manifold $M'$ over $B\pi$ where the map to $B\pi$ factors through the $p$-skeleton $B\pi^{(p)}$ of $B\pi$.
For $2\leq p\leq 4$ we can perform surgeries on $M'$ to convert the map $M' \to B\pi^{(p)}$ to a map $M'' \to B\pi^{(p)}$ that induces an isomorphism on fundamental groups.  By Kreck's modified surgery (\cref{thm:stablediffeoclasses}),  $M$ and $M''$ are stably diffeomorphic, and thus after connected sums with copies of $S^2 \times S^2$, $c$ is homotopic to a map factoring through $B\pi^{(p)}$.

The converse also holds, that is if one can find a map $f \colon M \to B\pi^{(p)}$ for $2\leq p\leq 4$ inducing the chosen isomorphism on fundamental groups, then $(M,c)$ lies in $F_{p,4-p}$.

Now, to state the next lemma, let $X$ be any CW complex and let $X^{(p)}$ be its $p$-skeleton. Denote the barycentres of the $p$-cells $\{e_i^p\}$ of $X$ by $\{b_i^p\}_{i \in I}$. Given an element $[M \to X^{(p)}] \in \Omega_n^{\spin}(X^{(p)})$, denote the regular preimage of the barycentre $\{b_i^p\} \in X^{(p)}$ by $N_i \subset M$. Note that $[N_i] \in \Omega_{n-p}^{\spin}$, since the normal bundle of $N_i$ in $M$ is trivial, and hence $N_i$ inherits a spin structure from $M$.

The following lemma is a consequence of \cite[Lemma~2.5]{KLPT15}.

\begin{lemma}\label{lem:arf}
The canonical map $\Omega_n^{\spin}(X^{(p)}) \to H_p(X^{(p)};\Omega_{n-p}^{\spin})$ that comes from the spectral sequence (see below) coincides with the map
\[\begin{array}{rcl} \Omega_n^{\spin}(X^{(p)}) &\to&  H_p(X^{(p)};\Omega_{n-p}^{\spin}) \\
\lbrack M \to X^{(p)}\rbrack  &\mapsto&  \Big[\sum\limits_{i \in I} [N_i]\cdot e_i^p\Big]. \end{array}\]
The same holds with oriented bordism replacing spin bordism.
\end{lemma}

The map in the statement of the lemma $\Omega_n^{\spin}(X^{(p)}) \to H_p(X^{(p)};\Omega_{n-p}^{\spin})$ is sometimes called an edge homomorphism. It arises as follows.  The abutment of the Atiyah-Hirzebruch spectral sequence $\Omega_n(\xi|_{X^{(p)}}) = F_{n,0}$ maps to its quotient by the first filtration step $F_{p,n-p}$ that differs from $F_{n,0}$.  This term is indeed $F_{p,n-p}$, since the homology of $X^{(p)}$ vanishes in degrees greater than $p$, thus $E^2_{s,t} = E^{\infty}_{s,t} =0$ for all $s>p$.   We have $F_{n,0}/F_{p,n-p} \cong E^{\infty}_{p,n-p}$.  The target, $H_p(X^{(p)};\Omega_{n-p}^{\spin})$, is the $E^2_{p,n-p}$ term of the spectral sequence.  Since no differentials have image in $E^2_{p,n-p}$, we have that $E^{\infty}_{p,n-p} \subseteq E^{2}_{p,n-p} = H_p(X^{(p)};\Omega_{n-p}^{\spin})$, and so the composition
\[\Omega_n^{\spin}(X^{(p)}) = F_{n,0} \to F_{n,0}/F_{p,n-p} \xrightarrow{\simeq} E^{\infty}_{p,n-p} \to E^{2}_{p,n-p} = H_p(X^{(p)};\Omega_{n-p}^{\spin})\]
gives the desired map.

We will need the following standard fact, so we give the argument here separately.
Recall that $\pi$ can be any finitely presented group.

\begin{lemma}\label{lem:coefficients-split}
The subgroup $F_{0,4} \subset \Omega_4^{\spin}(B\pi)$ is isomorphic to $\Omega_4^{\spin}$, and is a direct summand of $\Omega_4^{\spin}(B\pi)$.
\end{lemma}

\begin{proof}
The term $E^2_{0,4} \cong H_0(B\pi;\Omega_4^{\spin}) \cong \Omega_4^{\spin}$.  The maps $\pt \to B\pi$ and $B\pi \to \pt$ induce maps between the Atiyah-Hirzebruch spectral sequences. But the spectral sequence for a point is supported in a single column, so there are maps $\Omega_4^{\spin} \to \Omega_4^{\spin}(B\pi) \to \Omega_4^{\spin}$.  The composition $\pt \to B\pi \to \pt$ is the identity, and maps between the spectral sequences are natural, so the map $\Omega_4^{\spin} \to \Omega_4^{\spin}(B\pi)$ splits. Thus $\Omega_4^{\spin}$ is a direct summand as claimed.
\end{proof}

We consider the Atiyah-Hirzebruch spectral sequence for $\Omega_4^{\spin}(B\pi)$.

\begin{prop}\label{prop:spin-AHSS}
There is a filtration
 \[0 \subseteq F_3 \subseteq F_2 \subseteq F_1 \subseteq \Omega_4^{\spin}(B\pi)\]
by subgroups such that the following holds.
\begin{enumerate}
\item The quotient $\Omega_4^{\spin}(B\pi)/F_1$ is isomorphic to \[\ker\big( d^3_{4,0} \colon \ker \big(\Sq_2\circ \red_2 \colon H_4(\pi;\Z) \to H_2(\pi;\Z/2)\big) \to H_1(\pi;\Z/2)\big).\]  The image of an element $(M,c)$ in this subgroup of $H_4(B\pi;\Z)$ is $c_*([M])$, and is denoted $\pri(M)$.
Here $\Sq_2 \colon H_4(\pi;\Z/2) \to H_2(\pi;\Z/2)$ denotes the map on homology that is dual to the usual cohomology Steenrod operation~$\Sq^2 \colon H^2(\pi;\Z/2) \to H^4(\pi;\Z/2)$ and $\red_2$ is the reduction modulo two of the coefficients.
\item The quotient $F_1/F_2$ is isomorphic to $H_3(\pi;\Z/2)/\im(d^2_{5,0})$, where
\[d^2_{5,0} = \Sq_2\circ \red_2 \colon H_5(\pi;\Z) \to H_3(\pi;\Z/2).\]
The image of an element $(M,c)$ is denoted $\msec(M)$. If $\pri(M)=0$, then $\msec(M)$ can be computed by \cref{lem:arf}.
\item The quotient $F_2/F_3$ is isomorphic to $H_2(\pi;\Z/2)/\im(d^2_{4,1},d^3_{5,0})$, where
\[d^2_{4,1} = \Sq_2 \colon H_4(\pi;\Z/2) \to H_2(\pi;\Z/2).\] If $\pri(M)=\msec(M)=0$, then $\ter(M)$ can be computed by \cref{lem:arf}.
\item The subgroup $F_3$ is a direct summand isomorphic to $\Z$, and the map to $\Z$ is given by $\sigma/16$, that is take the signature and divide by $16$.
\end{enumerate}
\end{prop}

\begin{proof}
The proposition follows from computing using the Atiyah-Hirzebruch spectral sequence, with $E^2_{p,q} = H_p(B\pi;\Omega_q^{\spin})$.  We know from \cref{lem:coefficients-split} that the coefficients split, so recalling Rochlin's theorem that $\Omega_4^{\spin} \cong 16\Z$, this proves the last item.

Spin bordism in the lower dimensions is given by $\Omega_3^{\spin}=0$, $\Omega_2^{\spin} =\Z/2$, $\Omega_1^{\spin} =\Z/2$ and $\Omega_0^{\spin} \cong \Z$~ \cite{Milnor-spin}, \cite{ABP-Spin}.

For the spin bordism Atiyah-Hirzebruch spectral sequence,
the $E^2$-page differential $d_2$ is dual to the Steenrod square $\Sq^2$, together with reduction of the coefficients mod 2 if necessary, according to~\cite[Theorem 3.1.3]{teichnerthesis}. This implies that the differential vanishes when the codomain is $H_p(B\pi;\Omega_q^{\spin})$ and $p<2$, since $\Sq^r \colon H^n \to H^{n+r}$ is zero whenever $n <r$.

There are a couple of potentially nontrivial $d^2$ and $d^3$ differentials, that are recorded in the proposition.
\end{proof}

Now, computing the secondary and tertiary invariants as inverse images is not a very feasible task.
A particular difficulty is that one typically has to find a bordant manifold such that the map to $B\pi$ can be homotoped onto the $p$-skeleton, where $p=3$ for $\msec$ and $p=2$ for $\ter$.  Then one takes inverse images with framings.  The goal of this paper, as explained in the introduction, is to give algebraic criteria that decide whether these invariants vanish.  A particularly nice feature is that we show that the algebraic secondary obstruction does not depend on the way in which the primary invariant vanishes. That is, while it needs the primary obstruction to vanish in order to be well-defined, it does not depend on the choice of vanishing i.e.\ of splitting of the short exact sequence.  Similarly, the algebraic tertiary obstruction needs the first two obstructions to vanish in order to be well-defined, but it does not depend on the way in which they vanish.
These features increase the computability of our obstructions.
\cref{section:stable-classn-ZxZ2} provides an explicit instance of this phenomenon.

\section{Topological 4-manifolds and stable homeomorphism}\label{sec:topological-case}

Our results apply to the stable homeomorphism classification of topological spin $4$-manifolds, as we explain in this section.
First we point out that the invariants from \cref{thm:1.1} apply in the topological category, as asserted in the introduction.  Note that the signature $\sigma(M)$, $\pi_2(M)$, and $\lambda_M$ are all homotopy invariant, so give rise to stable homeomorphism invariants.  The Kervaire-Milnor invariant $\tau_M$ is defined using intersections: by topological transversality~\cite[Theorem~9.5A]{Freedman-Quinn}, such intersections can be counted in topological manifolds.  \cref{thm:1.1} requires that we compute these invariants on all possible 1-skeleton connected sums $(N_1,\alpha_1) \#_1  (N_2,-\alpha_2)$.  To avoid handle structures, one can instead take the connected sum $N_1 \# -N_2$ and then perform surgery on circles, framed using the $\alpha_i$, to make the fundamental group isomorphic to $\pi$, and in this way construct the test manifolds on which we claim one must compute our algebraic invariants in order to decide stable homeomorphism.

Now we explain why the classifications coincide.
To start, we have forgetful map
\[\Z \cong \Omega_4^{\spin} \xrightarrow{\cdot 2} \Omega_4^{TOPSpin} \cong \Z.\]
Recall that the isomorphism with $\Z$ is given by the signature, taking account of its divisibility. That is, we have isomorphisms $\sigma  \colon \Omega_4^{SO} \raiso \Z$, $\sigma/16 \colon \Omega_4^{\spin} \raiso \Z$ and $\sigma/8 \colon \Omega_4^{TOPSpin} \raiso \Z$.
The cokernel
\[\coker(\Omega_4^{\spin} \to \Omega_4^{TOPSpin}) \cong \Z/2\]
is detected by the Kirby-Siebenmann invariant.  Thus in particular the Kirby-Siebenmann invariants must coincide for stably homeomorphic $4$-manifolds, but for spin manifolds the Kirby-Siebenmann invariant is determined by the signature.

We want to apply our stable classification, but for this we need smooth manifolds.  Consider two $4$-manifolds $M$ and $M'$ that are spin and have nonvanishing Kirby-Siebenmann invariant.  Suppose that $\pi_1(M) \cong \pi_1(M')$ has the Secondary and Tertiary Properties.  We want to know if the two manifolds are stably homeomorphic.
Let $W$ be the $E_8$-manifold.
The manifolds $M \# W$ and $M' \# W$ have vanishing Kirby-Siebenmann invariants, and so are stably smoothable by \cite[Section~8.6]{Freedman-Quinn}.  We then have smooth manifolds, and so we can apply our classification programme to decide whether they are stably diffeomorphic.  This involves algebraic invariants that are independent of the smooth structure.   Suppose that we discover $M \# W$ and $M'\# W$ to be stably diffeomorphic.  Then $M \# W \# -W$ is stably homeomorphic to $M' \# W \# -W$.  But then $W \# -W$ is $B$-null bordant, so is stably homeomorphic to $S^4$.  It follows that $M$ and $M'$ are stably homeomorphic.   Thus the same programme as in the smooth case determines whether two closed, spin topological $4$-manifolds, whose fundamental groups satisfy \cref{secondary-property} and \cref{tertiary-property}, are stably homeomorphic.


\section{Whitehead's certain exact sequence and the \texorpdfstring{$\Gamma$}{Gamma}  functor}\label{sec:certain-exact-sequence}

In this section, we give some background on the ``certain exact sequence'' of Whitehead~\cite{whitehead}, material that we will make extensive use of in our investigation of the Secondary Property, and its proof for some families of groups.

Let $Y$ be a simply connected CW complex. Let $C_n:=\pi_n(Y^{(n)},Y^{(n-1)})$ be the $n$th cellular chain group. Let $d_n\colon C_n\to C_{n-1}$ denote the canonical boundary map arising from the long exact sequence of the pair. Let $i_n \colon \pi_n(Y^{(n)})\to C_n$ be the map induced by the inclusion of pairs. Let $Z_n:=\ker(d_n)$ and let $B_n:= \im (d_{n+1})$. Then of course $H_n(Y;\Z)\cong Z_n/B_n$.
Let $h_n\colon \pi_n(Y)\to Z_n/B_n$ be the Hurewicz map, and let $b_n\colon Z_n/B_n\to \ker(i_{n-1})$ be induced from the boundary map $\pi_n(Y^{(n)},Y^{(n-1)})\to \pi_{n-1}(Y^{(n-1)})$. Furthermore, let $\iota_n\colon \ker(i_n)\to \pi_n(Y)$ be induced by the inclusion $Y^{(n)}\to Y$.  These maps enable us to formulate the ``certain exact sequence'' of J.H.C.~Whitehead, which describes the kernel and cokernel of the Hurewicz map.

\begin{thm}[{\cite[Section 10]{whitehead}}]
	 There is a long exact sequence
	\[\cdots\to \pi_4(Y)\xrightarrow{h_4} H_4(C_n,d_n)\xrightarrow{b_4}\ker(i_3) \xrightarrow{\iota_3} \pi_3(Y) \xrightarrow{h_3}H_3(C_n,d_n)\to 0.\]
\end{thm}

The extra ingredient that makes this sequence extremely useful is a description of~$\ker(i_3)$ in terms of~$\pi_2(Y)$. This uses the $\Gamma$ group, which we will now define.

\begin{defi}
Let $A$ be an abelian group. Then $\Gamma(A)$ is an abelian group with generators the elements of $A$. We write $a$ as $v(a)$ when we consider it as an element of $\Gamma(a)$. The group $\Gamma(A)$ has the following relations:
	 \[\{v(-a)-v(a)\mid a\in A\}  \quad \text{ and }
	\]
	 \[\{v(a+b+c)-v(b+c)-v(c+a)-v(a+b)+v(a)+v(b)+v(c)\mid a,b,c\in A\}.\]
\end{defi}

We remark that the symbol $v$ has no meaning on its own, rather it is used to differentiate the generating set for $\Gamma(A)$ from the generating set for $A$.

The functor $\Gamma$ is the universal quadratic functor, in the following sense.  A map $f\colon A\to B$ of abelian groups is called quadratic if $f(a)=f(-a)$ for all $a\in A$ and for all $a,a'\in A$ the map
\[A\times A\to B,~(a,a')\mapsto f(a+a')-f(a)-f(a')\]
is bilinear. The map $j\colon A\to \Gamma(A)$ sending $a$ to $v(a)$ is quadratic.
The functor $\Gamma$ satisfies the universal property that for every quadratic map $f\colon A\to B$, there is a unique homomorphism $\Gamma(f)\colon \Gamma(A)\to B$ with $f=\Gamma(f)\circ j$.

\begin{lemma}[{\cite[p.~62]{whitehead}}]
	\label{lem:gammafree}
	If $A$ is free abelian with basis $\mathcal{B}$, then $\Gamma(A)$ is free abelian with basis
	\[\{v(b), v(b+b')-v(b)-v(b')\mid b,b'\in \mathcal{B}\}.\]
\end{lemma}

		In particular, if $A$ is free abelian, then sending $v(a)$ to $a\otimes a$ defines an isomorphism between $\Gamma(A)$ and the subgroup of symmetric tensors of $A\otimes_\Z A$, that is, the subgroup generated by $\{a\otimes a\mid a\in A\}$.

\begin{thm}
	[{\cite[Sections 10 and 13]{whitehead}}]
	Let $\eta\colon S^3\to S^2$ be the Hopf map.
	The map $\widehat{\eta}\colon \Gamma(\pi_2(Y))\to \ker(i_3\colon \pi_3(Y^{(3)})\to \pi_3(Y^{(3)},Y^{(2)}))$ given by $v(\alpha)\mapsto \alpha\circ\eta$ is an isomorphism.
	In particular, for every simply-connected CW complex $Y$ we have the exact sequence:
	\[\ldots\pi_4(Y)\to H_4(Y;\Z)\xrightarrow{\widehat{\eta}^{-1}\circ b_4}\Gamma(\pi_2(Y))\xrightarrow{\iota_3 \circ \widehat{\eta}} \pi_3(Y)\to H_3(Y;\Z)\to 0.\]
\end{thm}

Now let $L$ be any CW complex and as above let $\eta\colon S^3\to S^2$ be the Hopf map.

\begin{lemma}\label{lem:gamma-eta}{\cite[Section 13]{whitehead}}
The map
\[\Gamma(\eta) \colon \Gamma(\pi_2(L))\to \pi_3(L); v(\alpha)\mapsto \alpha \circ \eta\]
yields a well defined homomorphism.
\end{lemma}

For later use we restate the Whitehead exact sequence in the form that we will need it, namely for general CW complexes that need not be simply connected.

\begin{theorem}\label{thm:whitehead-seq}
For a CW complex $L$, the following sequence is exact:
\[H_4(\wt{L};\Z) \to \Gamma(\pi_2(L)) \xrightarrow{\Gamma(\eta)} \pi_3(L) \to H_3(\wt{L};\Z) \to 0.\]
Each of the terms in functorial in $L$ and the maps in the sequence above are natural.
\end{theorem}

To finish the section, we record a couple of preliminary facts on the $\Gamma$ groups that we will use throughout the rest of the paper.

\begin{cor}\label{lem:S-iso}
Let $K$ be a 2-complex. Then the map $\Gamma(\eta) \colon \Gamma(\pi_2(K)) \to \pi_3(K)$ is an isomorphism.
\end{cor}

\begin{proof}
  For $K$ a $2$-complex, both the third and fourth homology groups of $\wt{K}$ vanish, so $\Gamma(\eta)$ is an isomorphism by \cref{thm:whitehead-seq}.
\end{proof}

\begin{cor}\label{lem:rewriting-pi3-K}
Let $K$ be a $2$-complex.  Every element in $\pi_3(K)$ is a sum of elements of the form $\beta \circ \eta$
with $\beta\in\pi_2(K)$.
\end{cor}

\begin{proof}
Every element of $\Gamma(\pi_2(K))$ is a sum of elements $v(a)$ with $a\in \pi_2(K)$ by definition. Since $\Gamma(\eta)$ is a surjection and $\Gamma(\eta)(v(a))=a\circ \eta$, the corollary follows.
\end{proof}

Let \[\begin{array}{rcl}  T \colon N \otimes_{\Z} N' &\to & N' \otimes_{\Z} N \\ n \otimes n' & \mapsto & n' \otimes n \end{array}\]
be the transposition map.

\begin{lemma}
\label{lem:gammasum}
Let $N,N'$ be free $\Z$-modules. Then \[\Gamma(N\oplus N')\cong \Gamma(N)\oplus (N\otimes N') \oplus \Gamma(N').\] Moreover the inclusion into \[(N\oplus N')\otimes_\Z(N\oplus N')\cong (N\otimes_\Z N)\oplus (N\otimes_\Z N')\oplus(N'\otimes_\Z N)\oplus (N'\otimes_\Z N')\] is given by the direct sum of the inclusion $\Gamma(N)\to N\otimes_\Z N$, the diagonal map $(1,T)\colon N\otimes_\Z N'\to (N\otimes_\Z N')\oplus(N'\otimes_\Z N)$ and the inclusion $\Gamma(N')\to N'\otimes_\Z N'$.
\end{lemma}

\begin{proof}
Since $N \oplus N'$ is free, $\Gamma(N \oplus N')$ is the symmetric tensors. Embed $\Gamma(N \oplus N')$ into $(N \oplus N') \otimes (N \oplus N')$ and observe that the subgroup of symmetric tensors is isomorphic to $\Gamma(N) \oplus N \otimes N' \oplus \Gamma(N')$ with \[(n_1 \otimes n_2,m \otimes m',\ell_1' \otimes \ell_2') \mapsto n_1 \otimes n_2 + m \otimes m' + m' \otimes m + \ell_1' \otimes \ell_2'\] the inverse of the embedding.
\end{proof}


\chapter{The secondary obstruction}
\label{sec:sec}

Here is an outline of this chapter. First we discuss the equivariant intersection form and some of its different guises. In \cref{sec:sesq-forms-tate-cohomology} we associate an element of Tate cohomology  to the intersection form of a $4$-manifold with vanishing primary obstruction.  The Tate group measures whether or not the intersection form $\lambda_M$ is even, when restricted to the $H^2(K;\Z\pi)$ summand  of the stabilised second homotopy group $\pi_2(M) \oplus \Z\pi^n$.  The Secondary Property asserts that this element of the Tate group detects $\msec(M)$.

In \cref{cond}, we give a condition that we then prove in the remainder of the section is sufficient for the \cref{secondary-property} to hold for a group $\pi$.  In \cref{chapter:examples} we will use this condition to show that many families of groups have the \cref{secondary-property}.

\section{The equivariant intersection form}\label{subsection:equiv-int-form}

Let $M$ be a smooth, closed, spin, based $4$-manifold together with an identification $\pi_1(M) \raiso \pi$.
The equivariant intersection form
\[\lambda_M \colon \pi_2(M) \times \pi_2(M) \to \Z\pi\]
is defined as follows.  Identify $\pi_2(M) \cong H_2(M;\Z\pi)$ via the Hurewicz theorem.
Then for classes $x ,y \in H_2(M;\Z\pi)$, we have by definition
\[\lambda_M(x,y) = \langle PD^{-1}(y), x \rangle,\]
where $PD \colon H^2(M;\Z\pi) \to H_2(M;\Z\pi)$ is the Poincar\'{e} duality isomorphism given by cap product with the fundamental class $[M] \in H_4(M;\Z)$.
Here, and indeed throughout the article, we use the involution in the definition of cohomology $H^n(M;\Z;\pi) := H^n(\Hom_{\Z\pi}(C_*(M;\Z\pi)^t,\Z\pi))$ to consider chains as a right $\Z\pi$-module $C_*(M;\Z\pi)^t$, so that cohomology still carries a left $\Z\pi$-module structure.
Using naturality of the cap product and of evaluation, it is not hard to show that an orientation preserving homotopy equivalence $f \colon M \to M'$ of closed manifolds induces an isometry $\lambda_M \cong \lambda_{M'}$.

We identify the equivariant intersection form on $M$ with $\Z\pi$ coefficients with the intersection form on $\wt{M}$.
Pick a lift of each cell of $M$ to obtain an identification $\theta \colon C_*(M;\Z\pi) \raiso C_*(\wt{M};\Z)$.  We also have an isomorphism
\[\begin{array}{rcl}
\Psi \colon C_{cs}^*(\wt{M};\Z) & \xrightarrow{\cong} & \Hom_{\Z\pi}(C_*(\wt{M};\Z),\Z\pi)   \\
f &\mapsto & \big(a \mapsto \sum_{g \in \pi}f(g^{-1} a)\cdot g\big)
\end{array}\]
for $a \in C_*(\wt{M};\Z)$.  The inverse is given by sending $\phi \in \Hom_{\Z\pi}(C_*(\wt{M};\Z),\Z\pi)$ to the homomorphism that maps $x \in C_*(\wt{M};\Z)$ to the coefficient of the neutral group element of $\phi(x) \in \Z\pi$.
These two maps induce an isomorphism $\theta^* \circ \Psi_* \colon H^*_{cs}(\wt{M};\Z) \raiso H^*(M;\Z\pi)$, which fits into the following commuting diagram.
\[\xymatrix{H_2(M;\Z\pi)  \ar[r]^{PD^{-1}}_{\cong} \ar[d]_{\theta}^{\cong}  &  H^2(M;\Z\pi) \ar[r]^-{\ev} & \Hom_{\Z\pi}(H_2(M;\Z\pi),\Z\pi) \\
H_2(\wt{M};\Z) \ar[r]^{PD^{-1}}_{\cong} &  H^2_{cs}(\wt{M};\Z) \ar[r]^-{\ev \circ \Psi} \ar[u]^{\theta^* \circ \Psi}_{\cong} & \Hom_{\Z\pi}(H_2(\wt{M};\Z),\Z\pi) \ar[u]^{\theta^*}_{\cong}
}\]
In the diagram we also write $PD \colon H^2_{cs}(\wt{M};\Z) \to H_2(\wt{M};\Z)$ for Poincar\'{e} duality in~$\wt{M}$.

Thus, by taking the two routes from $H_2(M;\Z\pi)$ to $\Hom_{\Z\pi}(H_2(M;\Z\pi),\Z\pi)$, we see that for $x,y \in H_2(M;\Z\pi)$, we have
\[\lambda_M(x,y) = \langle PD^{-1}(y),x \rangle = \sum_{g \in \pi} \langle PD^{-1}(\theta(y)),g^{-1} \theta(x) \rangle \cdot g.\]
Now we will consider $w= \theta(x)$ and $u=\theta(y)$ in $H_2(\wt{M};\Z)$.  Note that $H_0(\wt{M};\Z) \cong H^4_{cs}(\wt{M};\Z) \cong \Z$.
In the final expression in the equation above, we write $PD \colon H^2_{cs}(\wt{M};\Z) \to H_2(\wt{M};\Z)$ for Poincar\'{e} duality.  Then since the cup product is signed commutative, we have, for $u,w \in H_2(\wt{M};\Z)$:
\begin{align*}
\langle PD^{-1}(u),g^{-1}w \rangle &= PD^{-1}(u) \cap g^{-1} w \\
&= PD^{-1}(u) \cup (g^{-1}  PD^{-1}(w) \cap [M]) \\
&= (g PD^{-1}(u) \cup PD^{-1}(w)) \cap [M] \\
&= (PD^{-1}(w) \cup g PD^{-1}(u)) \cap [M].
\end{align*}
We record the outcome in the following proposition.

\begin{prop}\label{prop:equivalence-intersection-forms}
  The map $\theta_* \colon H_2(M;\Z\pi) \to H_2(\wt{M};\Z)$ induces an isometry of the equivariant intersection form $\lambda_M(x,y) = \langle PD^{-1}(y),x \rangle$ with  the form
\[\lambda_{\wt{M}}(\theta_*(x),\theta_*(y)) = \lambda_{\wt{M}}(w,u)  = \sum_{g \in \pi} \big( (PD^{-1}(w) \cup g PD^{-1}(u)) \cap [M]\big) g  \in \Z\pi. \]
\end{prop}

It is well-known that $(PD^{-1}(x)\cup PD^{-1}(y))\cap [M] \in H_0(\wt{M};\Z) \cong \Z$ agrees with the geometric intersection of $S_x$ and $S_y$, where $S_x$ and $S_y$ are transverse, generically immersed surfaces in $\wt{M}$ representing $x$ and $y$ respectively.

Here each intersection point $p$ is counted with a sign depending on whether the orientation of $T_pS_x\oplus T_pS_y$ agrees with $T_pM$ or not.

Let $x_0$ be the chosen base point of $M$ and let $\wt x_0$ be a chosen lift of $x_0$ in $\wt M$. Given two transverse, generically immersed spheres $\alpha,\beta\in \pi_2(M,x_0)$, any intersection point $p$ determines an element of $\pi$ by choosing a path in $\alpha$ from $x_0$ to $p$ and concatenating it with a path from $p$ to $x_0$ in $\beta$. This element is $g$ precisely if the lifts $\wt \alpha\in \pi_2(\wt M, \wt x_0)$ and $g\wt \beta \in \pi_2(\wt M, gx_0)$ intersect at a lift $\wt p$ of $p$.  Hence we have the following statement, which we will use later in our discussion of the tertiary invariant.

\begin{prop}
  The above geometric count of intersections in $\Z\pi$ computes the equivariant intersection form.
\end{prop}

Using the description
\[\lambda_M(x,y)=\sum_{g\in\pi}((PD^{-1}(x)\cup gPD^{-1}(y))\cap [M])g,\]
it is easy to see that the form $\lambda_M$ is sesquilinear and hermitian.
But it is often not nonsingular or even nondegenerate, as can be computed by the universal coefficient spectral sequence~\cite[Theorem~2.3]{Levine77}.

\begin{proposition}\label{prop:int-form-singular}
  There is an exact sequence
\[0 \to H^2(\pi;\Z\pi) \to H_2(M;\Z\pi) \to \Hom_{\Z\pi}(H_2(M;\Z\pi),\Z\pi) \to H^3(\pi;\Z\pi) \to 0\]
where the middle map is the adjoint of the intersection form.
\end{proposition}

\begin{proof}
Recall that the universal coefficient spectral sequence has $E^2$ page $E_2^{p,q} \cong \Ext^q_{\Z\pi}(H_p(M;\Z\pi),\Z\pi)$, differential $d^r$ of degree $(1-r,r)$, and the sequence converges to $H^{p+q}(M;\Z\pi)$. Since $H_1(M;\Z\pi) =0$, the spectral sequence yields a filtration
$0 \subseteq F_{0,2} \subseteq H^2(M;\Z\pi)$ with $F_{0,2} \cong \Ext^2_{\Z\pi}(H_0(M;\Z\pi),\Z\pi)$ and
\[H^2(M;\Z\pi)/F_{0,2} \cong \ker\big(\Hom_{\Z\pi}(H_2(M;\Z\pi),\Z\pi) \to \Ext^3_{\Z\pi}(H_0(M;\Z\pi),\Z\pi)\big).\]
Note that $\Ext^i_{\Z\pi}(H_0(M;\Z\pi),\Z\pi) = H^i(\pi;\Z\pi)$ to obtain the exact sequence claimed.
\end{proof}

\section{Sesquilinear forms and Tate cohomology}\label{sec:sesq-forms-tate-cohomology}

In this section, and indeed throughout the rest of the article, when not specified every $\Z\pi$-module is assumed by default to be a left $\Z\pi$-module.

Let $N$ be a left $\Z\pi$-module. We write $N^*$ for $\Hom_{\Z\pi}(N,\Z\pi)$. This would a priori be a right $\Z\pi$-module, where $\Z\pi$ acts by right multiplication on the target, but we turn it into a left module using the involution on $\Z\pi$.

Below, and indeed also throughout the article, when we consider the tensor product $N \otimes_{\Z\pi} N'$ of two left $\Z\pi$-modules $N,N'$, we use the involution to turn $N$ into a right $\Z\pi$-module, so that the tensor product makes sense.

We denote the group of sesquilinear forms on $N$ by $\Sesq(N)$, that is the group of maps
\[\lambda\colon N \times N\to  \Z\pi\]
with $\lambda(am,bn)=a\lambda(m,n)\ol{b}$ for $a,b\in\Z\pi$ and $m,n\in N$. The group operation is the obvious addition, defined by $(\lambda + \lambda')(m,n) = \lambda(m,n) + \lambda'(m,n)$.  Equivalently $\Sesq(N)\cong \Hom_{\Z\pi}(N,N^*)$ by $\lambda(n)(m):=\lambda(m,n)$. This is a contravariant functor, where a map $f\colon N\to N'$ is sent to the map $f^*\colon \Sesq(N')\to \Sesq(N)$ with $(f^*\lambda)(m,m')=\lambda(f(m),f(m'))$.
 The group $\Z/2$ acts on $\Sesq(N)$ via $(T\lambda)(m,n)=\ol{\lambda(n,m)}$. Thus we can form the \emph{Tate cohomology group}
\[\widehat{H}^0(\Sesq(N)):=\widehat{H}^0(\Z/2;\Sesq(N))=\ker(1-T)/\im(1+T).\]
Note that $\ker(1-T)$ is precisely the hermitian forms, that is those for which  $\lambda(m,n)=\ol{\lambda(n,m)}$ for all $m,n \in N$.  In addition, every hermitian form has order two in the Tate cohomology, since $\lambda = T\lambda$ implies that $2\lambda = \lambda + T\lambda = (1+T)\lambda$.

\begin{definition}\label{defn:even}
  A hermitian form is called \emph{even} if it is in the image of $1+T \colon \Sesq(N) \to \Sesq(N)$, that is if it vanishes in $\widehat{H}^0(\Sesq(N))$.
\end{definition}

We therefore obtain a contravariant functor $\widehat{H}^0(\Sesq(-))$ from $\Z\pi$-modules to $\Z/2$-modules.
Let $(-)^*\otimes_{\Z\pi}(-)^*$ be the contravariant functor sending $N$ to $N^*\otimes_{\Z\pi}N^*$, where we consider the first $N^*$ as a right $\Z\pi$-module as described above; equivalently, we do not use the involution on this~$N^*$.

\begin{defi}
For left $\Z\pi$-modules $N$ and $N'$, define
\[\ba{rcl} \Phi_{N,N'}\colon N^*\otimes_{\Z\pi} N' &\to & \Hom_{\Z\pi}(N,N') \\
f\otimes n' &\mapsto & \big(n \mapsto f(n)n'\big).\ea\]
We also define
\[\ba{rcl}\Phi_N:=\Phi_{N,N^*}\colon N^*\otimes_{\Z\pi}N^* &\to & \Hom_{\Z\pi}(N,N^*)\cong \Sesq(N).\\
 f_1\otimes f_2 &\mapsto & \big((n_1,n_2) \mapsto f_1(n_1)\ol{f_2(n_2)}\big) .\ea\]
\end{defi}

\noindent The next lemma follows directly from the definition.

\begin{lemma}
\label{lem:phisum}
For two left $\Z\pi$-modules $N$ and $N'$, the map $\Phi_{N \oplus N'}$ is equivalent to the direct sum $\Phi_{N} \oplus \Phi_{N'} \oplus \Phi_{N,(N')^*} \oplus \Phi_{N',N^*}$, under the obvious isomorphisms of the domains and the codomains.
\end{lemma}

\begin{lemma}\label{lem:defn-Phi}
The map
\[\Phi_N\colon N^*\otimes_{\Z\pi}N^*\to \Sesq(N)\]
defines a natural transformation
\[\Phi\colon (-)^*\otimes_{\Z\pi}(-)^*\Rightarrow \Sesq(-).\]
\end{lemma}

\noindent The proof of \cref{lem:defn-Phi} is straightforward and we omit the details.

\begin{lemma}\label{lem:Theta-Iso}
The map
\[\begin{array}{rcl} \Theta_N\colon \Z/2\otimes_{\Z\pi}N^* &\to & \widehat{H}^0(N^*\otimes_{\Z\pi}N^*) \\
1\otimes f &\mapsto & [f\otimes f] \end{array}\]
determines a natural isomorphism \[\Theta\colon \Z/2\otimes_{\Z\pi}(-)^*\Rightarrow \widehat{H}^0((-)^*\otimes_{\Z\pi}(-)^*).\]
\end{lemma}

\begin{proof}
First we show that $\Theta_N$ is well-defined.
For all $f,f'\in N^*$ and for all $g\in\pi$, we have
\[\Theta_N(1\otimes gf)=[gf\otimes gf]=[fg^{-1}\otimes gf]=[f\otimes f]=\Theta_N(1\otimes f)\]
and
\begin{align*}
\Theta_N(1\otimes(f+f'))&=[(f+f')\otimes (f+f')]\\
&=[f\otimes f+f'\otimes f'+(1+T)(f\otimes f')]\\
&=[f\otimes f]+[f'\otimes f']=\Theta_N(1\otimes f)+\Theta_N(1 \otimes f')
\end{align*}
 This shows that $\Theta_N$ is well-defined as desired.

It is easy to see that $\Theta$ is natural, so it remains to prove that it is an isomorphism. For this we will define an inverse.
Let $N^*/2:=\Z/2\otimes_{\Z\pi}N^*$. The projection $N^*\to N^*/2$ defines a map $N^*\otimes_{\Z\pi}N^*\to N^*/2\otimes_{\Z/2}N^*/2$, and then applying $\widehat{H}^0$ we obtain a map
\[\widehat{H}^0(N^*\otimes_{\Z\pi}N^*)\to \widehat{H}^0(N^*/2\otimes_{\Z/2}N^*/2).\]
Since $N^*/2$ is a $\Z/2$-vector space, we have an isomorphism $\widehat{H}^0(N^*/2\otimes_{\Z/2}N^*/2)\cong N^*/2$ given by sending $[[f]\otimes [f]]\to [f]$. To see that this map makes sense and is an isomorphism, observe that symmetric modulo even forms are determined by the diagonal entries. So every element in $\widehat{H}^0(N^*/2\otimes_{\Z/2}N^*/2)$ is a sum of elements of the form $[[f]\otimes [f]]$.
The composition
\[\widehat{H}^0(N^*\otimes_{\Z\pi}N^*)\to \widehat{H}^0(N^*/2\otimes_{\Z/2}N^*/2) \raiso N^*/2\]
sends $[f\otimes f]\in \widehat{H}^0(N^*\otimes_{\Z\pi}N^*)$ to $[f]=1\otimes f\in N^*/2=\Z/2\otimes N^*$, and is therefore the desired inverse to $\Theta_N$.
\end{proof}

\begin{defi}\label{defn:Psi-H}
Define the natural transformation
\[\Psi:=\widehat{H}^0(\Phi)\circ\Theta\colon \Z/2\otimes_{\Z\pi}(-)^*\Rightarrow \widehat{H}^0(\Sesq(-)).\]
\end{defi}



\begin{lemma}
\label{lem:injphi}
Let $P,N$ be left $\Z\pi$-modules, and suppose in addition that $P$ is finitely generated projective. Then:
\begin{enumerate}[(i)]
\item \label{item:Phi-iso-one} $\Phi_{P,N}$ and $\Phi_{N,P}$ are isomorphisms.
\item \label{item:Phi-iso-two} $\Psi_P\colon  \Z/2 \otimes_{\Z\pi} P^* \to \widehat{H}^0(\Sesq(P))$ is an isomorphism.
\item \label{item:Phi-iso-three} $\Phi_P$ is an isomorphism.
\end{enumerate}
\end{lemma}


\begin{proof}~
\begin{enumerate}[(i)] \item
First suppose that $P \cong \Z\pi$. In this case $\Phi_{\Z\pi,N}$ and $\Phi_{N,\Z\pi}$ are obviously isomorphisms. Since $\Phi_{N,\Z\pi^n}$ is the direct sum of $n$ times $\Phi_{N,\Z\pi}$, it is an isomorphism. The same holds for $\Phi_{\Z\pi^n,N}$.

Now let $P$ be finitely generated projective and let $Q$ be such that $P\oplus Q$ is finitely generated free. Then the isomorphism $\Phi_{N,P\oplus Q}$ is the direct sum of $\Phi_{N,P}$ and $\Phi_{N,Q}$. Thus both these maps have to be isomorphisms as well. By the same argument $\Phi_{P,N}$ is an isomorphism.
\item Since $\Phi_P$ is an isomorphism by part (\ref{item:Phi-iso-one}), so is $\widehat{H}^0(\Phi_P)$. Furthermore, $\Theta$ is a natural isomorphism by \cref{lem:Theta-Iso}, and thus also the composition $\Psi_P=\widehat{H}^0(\Phi_P)\circ \Theta_P$ is an isomorphism.
\item Immediate from (\ref{item:Phi-iso-one}).  \qedhere
\end{enumerate}
\end{proof}

\begin{cor}\label{cor:same-kernel}
Let $N$ be a left $\Z\pi$-module.  Then the inclusion $N \to N \oplus \Z\pi$ induces isomorphisms
\[\ker \Phi_N^{\Z/2} \raiso \ker \Phi_{N \oplus \Z\pi}^{\Z/2} \text{ and }
\ker \widehat{H}^0(\Phi_N) \raiso \ker \widehat{H}^0(\Phi_{N \oplus \Z\pi}), \]
where $\Phi_{N}^{\Z/2}$ is the map induced on $\Z/2$-fixed points by $\Phi_N$.
\end{cor}

\noindent This corollary will be used in the proof of \cref{lem:ker-A-well-defined}.

\begin{proof}
We have a commutative diagram
{\small \[\centerline{\xymatrix @C-0.2cm{\big((N^* \oplus \Z\pi) \otimes_{\Z\pi} (N^*\oplus \Z\pi)\big)^{\Z/2} \ar[r]^-{\cong} \ar[d]_{\Phi_{N\oplus \Z\pi}} & (N^* \otimes_{\Z\pi} N^*)^{\Z/2} \oplus (\Z\pi \otimes_{\Z\pi} N^*) \oplus (\Z\pi \otimes_{\Z\pi} \Z\pi)^{\Z/2} \ar[d] \\
\Sesq(N \oplus \Z\pi)^{\Z/2} \ar[r]^-{\cong} & \Sesq(N)^{\Z/2} \oplus \Hom_{\Z\pi}(\Z\pi,N^*) \oplus \Sesq(\Z\pi)^{\Z/2}.
}}\]}
In the two terms in the right hand column, the middle summands are isomorphic, as are the third summands by \cref{lem:injphi}~\eqref{item:Phi-iso-one}. The map between the first summands is $\Phi_{N}^{\Z/2}$.  Since the horizontal maps are isomorphisms, the first part of the lemma follows from the commutativity of the diagram.

The diagram descends to a similar diagram on $\widehat{H}^0$ groups.  By \cref{lem:injphi}~\eqref{item:Phi-iso-two} and \cref{lem:Theta-Iso}, the analogous statements hold, namely the maps induced by the right hand vertical map in the diagram (with $\Z/2$ fixed points replaced by $\widehat{H}^0$) splits along the direct summands, and the maps on the second and third summands are isomorphisms.  The second part of the lemma then also follows from commutativity.
\end{proof}

\section{Relating \texorpdfstring{$H_3(\pi;\Z/2)$}{the third homology} with symmetric modulo even forms}

Let $K$ be a connected finite 2-complex with $\pi_1(K) \cong \pi$.  Let $(D_*,d_*)$ be the cellular $\Z\pi$-module chain complex of~$K$.
In this subsection we establish a map $H_3(\pi;\Z/2) \to \widehat{H}^0(\Sesq(H^2(K;\Z\pi)))$ and we use this map to formulate \cref{cond}.  This condition will then later be shown to imply the \cref{secondary-property}.

Given $x\in \ker d_2\cong \pi_2(K)$, we define $f_x\colon \Z\pi\to \ker d_2$ by $f_x(a)=ax$. By dualizing we obtain a map $f_x^*\colon H^2(K,\Z\pi)\to \Z\pi^*\cong\Z\pi$, where the map $\Z\pi^* \to \Z\pi$ is given by $f \mapsto \overline{f(1)}$, which is a left $\Z\pi$-module homomorphism. Hence $f_x^*$ is a left $\Z\pi$-module homomorphism. Note that for $\phi\in H^2(K;\Z\pi)$ we have $f_x^*(\phi)=\overline{\phi(x)}$.

\begin{lemma}\label{lem:pi2-H2-dual}
Let $K$ be a connected finite $2$-complex with $\pi_1(K)\cong\pi$. The canonical map
\[\begin{array}{rcl}
  I\colon \pi_2(K) &\to& \Hom_{\Z\pi}(H^2(K;\Z\pi),\Z\pi) \\
  x &\mapsto & f_x^* \end{array}\]
 is an isomorphism.
\end{lemma}

\begin{proof}
To see that $I$ is a left $\Z\pi$-module homomorphism, let $g \in \pi$.
We compute that
\[I(gx)(\varphi)=f_{gx}^*(\varphi) = \overline{\varphi(gx)} = \overline{g \cdot \varphi(x)} = \overline{\varphi(x)} \cdot \overline{g}=f_x^*(\varphi) \cdot \overline{g}=(g \cdot f_x^*)(\varphi)=gI(x)(\varphi).\]
Now let $(D_*,d_*)$ denote the cellular $\Z\pi$-module chain complex of $K$. Consider the exact sequence
\[D^1\xrightarrow{d^2}D^2\to H^2(K;\Z\pi)\to 0.\]
Since the functor $\Hom_{\Z\pi}(-,\Z\pi)$ is left exact, we obtain an exact sequence
\[0\to \Hom_{\Z\pi}(H^2(K;\Z\pi),\Z\pi)\to\Hom_{\Z\pi}(D^2,\Z\pi)\xrightarrow{(d^2)^*}\Hom_{\Z\pi}(D^1,\Z\pi)\]
and since $D_1$ and $D_2$ are finitely generated, free $\Z\pi$-modules, this is isomorphic to the sequence
\[0\to\Hom_{\Z\pi}(H^2(K;\Z\pi),\Z\pi)\to D_2\xrightarrow{d_2}D_1.\]
Hence $\Hom_{\Z\pi}(H^2(K;\Z\pi),\Z\pi)$ is isomorphic to $\ker d_2 = H_2(K;\Z\pi) \cong \pi_2(K)$.
The isomorphism $\Hom(D^2,\Z\pi)\to D_2$ sends $f_x^*$, considered as a function $D^2\to \Z\pi$, to $f(1)=x$, considered as an element in $D_2$.
\end{proof}

On the other hand, if we start with $\coker d^2 = H^2(K;\Z\pi)$, the evaluation map is not necessarily an isomorphism.

\begin{lemma}
	\label{lem:dual-identification}
	Let $(C_*,d_*)$ be a free resolution of $\Z$ as a $\Z\pi$-module. There is an exact sequence
	\[0\to H^2(\pi;\Z\pi)\to \coker d^2\xrightarrow{\ev} \Hom_{\Z\pi}(\ker d_2;\Z\pi)\to H^3(\pi;\Z\pi)\to 0\]
	where the map $\ev$ sends $[f]\in \coker d^2$ with $f\in C^2$ to $(x\mapsto f(x))$.
\end{lemma}

\begin{proof}
	Consider the dual chain complex $(C^*,d^*)$. We have an exact sequence
\begin{equation}\label{eqn:coker-d2-kerd4}
0\to H^2(\pi;\Z\pi)\to \coker d^2\xrightarrow{d^3} \ker d^4\to H^3(\pi;\Z\pi)\to 0.
\end{equation}
	From the exactness of $(C_*,d_*)$ it follows that $\ker d_2=\im d_3$ and that $C_4\xrightarrow{d_4}C_3\xrightarrow{d_3}\im d_3\to 0$ is exact. By left exactness of $\Hom_{\Z\pi}(-,\Z\pi)$, we deduce that $d^3\colon \Hom_{\Z\pi}(\im d_3,\Z\pi)\xrightarrow{\cong} \ker d^4$. Therefore $\ker d^4 \cong \Hom_{\Z\pi}(\ker d_2,\Z\pi)$.  Substituting into \eqref{eqn:coker-d2-kerd4} we obtain the exact sequence
	\[0\to H^2(\pi;\Z\pi)\to \coker d^2\to \Hom_{\Z\pi}(\ker d_2;\Z\pi)\to H^3(\pi;\Z\pi)\to 0.\]
	The middle map now used $d^3$ and $(d^3)^{-1}$, so is given simply by the evaluation as claimed.
\end{proof}

\begin{rem}
Let $M$ be a closed oriented $4$-manifold with $\pi_1(M)=\pi$ and let $(C_*^M,d_*^M)$ be a $\Z\pi$-chain complex arising from a handle decomposition. 
Let $C_*^M[0,2]$ be the the 2-dimensional truncation of $C_*^M$ and let $j\colon C_*^M[0,2]\to C_*^M$ be the inclusion. Then the exact sequences from \cref{lem:dual-identification} and from \cref{prop:int-form-singular} fit into the following commutative diagram.
{\small \[\xymatrix @R-0.5cm@C-0.57cm {
0\ar[r]&H^2(\pi;\Z\pi)\ar[r]\ar[dd]^-=&H_2(M;\Z\pi)\ar[dd]_-{j^*\circ PD^{-1}}\ar[dr]_-{PD^{-1}}\ar[rr]^-{\lambda_M}&&\Hom_{\Z\pi}(H_2(M;\Z\pi),\Z\pi)\ar[r]\ar[dd]^-{(j_*)^*}&H^3(\pi;\Z\pi)\ar[r]\ar[dd]^-=&0\\
&&&H^2(M;\Z\pi)\ar[dl]_-{j^*}\ar[ur]_-{\ev}&&&\\
0\ar[r]&H^2(\pi;\Z\pi)\ar[r]&\coker d^2_M\ar[rr]^-{\ev}&&\Hom_{\Z\pi}(\ker d^M_2,\Z\pi)\ar[r]&H^3(\pi;\Z\pi)\ar[r]&0
}\]}
This identifies the failure of the equivariant intersection form to be nonsingular with the failure of the map $\ev$ from \cref{lem:dual-identification} to be an isomorphism.
\end{rem}

As above, consider a CW-complex model for $B\pi$ with finite $2$-skeleton $K$.
Now we define a map $\iota$, which will be important for relating the secondary invariant to pairings.

\begin{defi}\label{defn:iota}
Consider the Leray-Serre spectral sequence computing $H_*(K;\Z/2)$ from the fibration $\wt K\to K\to B\pi$ with $E^2$-term $E^2_{p,q}:=H_p(\pi;H_q(\wt K;\Z/2))$. Identify $H_2(\wt K;\Z/2)$ with $\Z/2\otimes_\Z \pi_2(K)$. Then $H_0(\pi;H_2(\wt K;\Z/2))\cong \Z/2\otimes_{\Z\pi}\pi_2(K)$. Since $H_1(\wt K;\Z/2)=0$, the entire $q=1$ row of the $E^2$ page vanishes, so the $d_3$ differential gives a map:
\[\iota:=d_3\colon H_3(\pi;\Z/2)\to \Z/2\otimes_{\Z\pi} \pi_2(K).\]
\end{defi}

\begin{lemma}\label{lem:iota-inj}
The map $\iota\colon H_3(\pi;\Z/2)\to \Z/2\otimes_{\Z\pi} \pi_2(K)$ is injective.
\end{lemma}

\begin{proof}
In the spectral sequence, $\ker \iota \cong E^{\infty}_{3,0}$. Since $E^{\infty}_{3,0}$ is isomorphic to the quotient of $H_3(K;\Z/2)$ by largest subgroup in the filtration of it determined by the Leray-Serre spectral sequence, we obtain a surjection $H_3(K;\Z/2)\to \ker \iota$. Since $K$ is $2$-dimensional, we have $H_3(K;\Z/2)=0$ and thus $\iota$ is injective.
\end{proof}

\begin{defi}\label{defn:map-A}
Let $H:=H^2(K;\Z\pi)$.  Define
\[A \colon H_3(\pi;\Z/2)\to \widehat{H}^0(\Sesq(H))\]
to be the composition $\Psi_H\circ (\Id \times I) \circ \iota$.
\end{defi}

\begin{lemma}\label{lem:ker-A-well-defined}
For any other choice $K'$ as the 2-skeleton of $B\pi$ and every map $\phi\colon K\to K'$ that induces the identity on fundamental groups, the diagram
\[\xymatrix{
H_3(\pi;\Z/2)\ar[r]^A\ar[dr]^{A'}&\wh H^0(\Sesq(H))\ar[d]^{\phi_*}\\
&\wh H^0(\Sesq(H'))
}\]
commutes, where $A'$ and $H'$ are the analogous definitions of $A$ and $H$ for $K'$ instead of $K$. In particular, the kernel of $A$ does not depend on the choice of finite 2-complex $K$ with $\pi_1(K) = \pi$, but only on $\pi$.
\end{lemma}

\begin{proof}
Consider the following diagram where $\iota'$ denotes the map $\iota$ for $K'$ instead of $K$.
\[\xymatrix @C+0.2cm {H_3(\pi;\Z/2) \ar @{^{(}->}[r]^-{\iota} \ar @{_{(}->} [dr]^-{\iota'} & \Z/2 \otimes_{\Z\pi} \pi_2(K) \ar[r]^-{\Psi_H} \ar[d] & \widehat{H}^0(\Sesq(H^2(K;\Z\pi))) \ar[d] \\
& \Z/2 \otimes_{\Z\pi} \pi_2(K') \ar[r]^-{\Psi_{H'}} & \widehat{H}^0(\Sesq(H^2(K';\Z\pi)))
}\]
Here we use the isomorphism $I \colon \pi_2(K)\cong H^*$ from \cref{lem:pi2-H2-dual} in order to consider $\Psi_H$ as a map from $\Z/2\otimes \pi_2(K)$ to $\widehat{H}^0(\Sesq(H))$, and similarly for $H' := H^2(K';\Z\pi)$.
The triangle commutes by naturality of the spectral sequence defining $\iota$, and the square commutes by naturality of $\Psi$.

It follows that the kernel of $A$ is contained in the kernel of $A'$. By switching the r\^oles of $K$ and $K'$, the kernel of $A$ is independent of the choice of $K$.
\end{proof}

\noindent In light of the previous lemma, $\Tate(\pi):=H_3(\pi;\Z/2)/\ker(A)$ is independent of the choice of $K$ and we can formulate the following condition.

\begin{cond}
\label{cond}
The sequence
\[H_5(\pi;\Z)\xrightarrow{\Sq_2\circ \red_2}H_3(\pi;\Z/2)\xrightarrow{A_{\pi}}\Tate(\pi)\]
is exact at $H_3(\pi;\Z/2)$, where $A_{\pi}$ is the projection.
\end{cond}

\begin{rem}\label{rem:sec-cond}
	By definition, \cref{cond} is equivalent to the exactness of the sequence
	\[H_5(\pi;\Z)\xrightarrow{\Sq_2\circ \red_2}H_3(\pi;\Z/2)\xrightarrow{A}\widehat H^0(\Sesq(H))\]
	for any choice of 2-complex $K$ and $H:=H^2(K;\Z\pi)$ as in \cref{defn:map-A}.	
\end{rem}

\noindent By the end of this section we will have proven the following theorem.

\begin{thm}
\label{thm:mainsec}
If \cref{cond} holds for $\pi$, then $\pi$ has the \cref{secondary-property}.
\end{thm}

Thus we will have reduced the Secondary Property, for a given group $\pi$, to verifying \cref{cond} for that group.  The classes of group for which we have been able to show that \cref{cond} holds are described in \cref{sec:ex,sec:abelian-groups}.  In \cite{KLPT15}, we proved the Secondary Property by constructing model manifolds and computing their intersection forms and their secondary invariant in $H_3(\pi;\Z/2)$.  Here we cannot construct model manifolds, but as long as we can compute $H_3(\pi;\Z/2)$ and the map $\Sq_2 \circ \red_2$, then by \ref{thm:mainsec} we can test the Secondary Property on ``model'' intersection forms produced by the map~$A$.

The following lemma will not be used in this section, but it may become useful in the future because, unlike the Tate group $\widehat{H}^0(\Sesq(H))$, the subgroup of \emph{weakly even} hermitian forms modulo even forms on a $\Z\pi$-module $H$ is unchanged under adding projective modules $H\mapsto H\oplus P$.

Recall the exact sequence from the Leray-Serre spectral sequence for $\widetilde K\to K\to B\pi$
\[
0\ra H_3(\pi;\Z/2)  \overset{\iota}{\ra}   \Z/2 \otimes_{\Z\pi} H_2(K;\Z\pi) \overset{H_2(\e_2)}{\ra}  H_2(K;\Z/2),
\]
where $\e_2\colon \Z\pi\to\Z/2$ is the nontrivial ring homomorphism, inducing the map on coefficients in the final map. Also recall that $H:=H^2(K;\Z\pi)$ has dual module $H^* \cong H_2(K;\Z\pi) \cong \pi_2(K)$ and that we defined a map $\Psi_H\colon  \Z/2 \otimes_{\Z\pi} H^* \to \widehat{H}^0(\Sesq(H))$ via
\[
\Psi_H(1\otimes f)(h_1, h_2) = f(h_1)\cdot \overline{f(h_2)} \text{ for } f\in H^*,\, h_i\in H.
\]

\begin{lemma}\label{lem:weaklyeven}
For every $f\in H^*$, the hermitian form $\Psi_H(1\otimes f)$ is weakly even if and only if $1\otimes f$ is in the image of $\iota$, or equivalently, if and only if $H_2(\e_2)(1\otimes f)=0$.
\end{lemma}

\begin{proof}
For $h\in H$, write $f(h)= \sum_{i=1}^n m_i g_i$ and compute in $\Z\pi$
\[
f(h)\cdot\overline{f(h}) = \sum_{i,j=1}^n m_i m_j (g_i\bar g_j) = \sum_{i=1}^n (m_i)^2  + \sum_{i<j} m_i m_j (g_i \bar g_j + g_j \bar g_i).
\]
The second sum is of the form $a + \bar a$ and the first sum is just a multiple of the unit.
That multiple $\sum_{i=1}^n (m_i)^2$ is modulo 2 the same as $\sum_{i=1}^n m_i \equiv \e_2(f(h))$, so we conclude that the above element is of the form $a+\overline{a}$ if and only if $\e_2(f(h))=0$.

The property of $\Psi_H(1\otimes f)$ being weakly even means that $\Psi_H(1\otimes f)(h, h)$ lies in the subgroup $\{a+\overline{a} \mid a\in \Z\pi\}$ for all $h\in H$. By the above, this is equivalent to the vanishing of the composition $\e_2\circ f\colon H \to \Z\pi\to\Z/2$. Finally, we use that $H^2(\e_2)\colon H^2(K;\Z\pi) \to H^2(K;\Z/2)$ is surjective because the next term in the Bockstein sequence is $H^3(K;\ker(\e_2))=0$. Hence the vanishing of $\e_2 \circ f$ is equivalent to the vanishing of $H_2(\e_2)(f)=H_2(\e_2)(1\otimes f)$, which is our original claim.
\end{proof}

\section{Proof of the Secondary Property Theorem}\label{subsection:proof-thm-mainsec}

Our aim for this section is to prove \cref{thm:mainsec}.
However, before we start the proof, we switch to a dual picture that we find convenient for the statements and the proofs of some imminent lemmas.
By  \cite[Proposition~1.7]{KPT18}, the short exact sequence
\[0\to \ker d_2\to C_2\oplus H_2(M;\Z\pi)\to \coker d_3 \to 0.\]
splits if and only if the primary obstruction vanishes. Recall that we write $C_* := C_*(M;\Z\pi)$ for the handle chain complex coming from a choice of handle decomposition of $M$, with boundary maps $d_i$.

In the discussion of the secondary obstruction we work with a slightly different ``dual'' formulation.
For this, let $(C'_*,d'_*)$ denote the $\Z\pi$-module chain complex coming from a dual handle decomposition of~$M$. We denote $M$ endowed with the dual handle decomposition by $M^d$.  The above sequence for $M^d$ has the form
\[0\to \ker d'_2\to C'_2\oplus H_2(M^d;\Z\pi)\to \coker d'_3 \to 0.\]
We can identify the chain complex of the dual handle decomposition with the dual of the chain complex of the first handle decomposition, i.e.\ $C'_*\cong C^{4-*}$. Apply this to the above sequence to yield
\[0\to \ker d^3\to C^2\oplus H^2(M;\Z\pi)\to \coker d^2 \to 0.\]
Since the primary invariant is independent of the chosen handle decomposition of $M$, the last sequence splits if and only if $\pri(M)=0$.
Note that $\coker d^2 \cong H^2(M^{(2)};\Z\pi)$.

Next, we show that every splitting of this dual short exact sequence can be realised by a geometric map $M^{(3)} \to M^{(2)}$.

\begin{prop}
\label{prop:geosplit}
For every splitting
\[s=(s_1,s_2)\colon H^2(M^{(2)};\Z\pi) \to C^2\oplus H^2(M;\Z\pi),\]
there is a map $f\colon M^{(3)}\to M^{(2)}$ such that
\[f^*=s_2\colon H^2(M^{(2)};\Z\pi)\to H^2(M^{(3)};\Z\pi)\cong H^2(M;\Z\pi).\]
\end{prop}

\begin{proof}
Let $j=(j_1,j_2) \colon C^2 \oplus H^2(M;\Z\pi) \to H^2(M^{(2)};\Z\pi)$ be the map from the short exact sequence above, so $j$ is split by $s= (s_1,s_2)$.  That is, $j_1$ is the projection $C^2 \to C^2/\im d^2$, and $j_2 \colon \ker d^3/\im d^2 \hookrightarrow C^2/\im d^2$ is the inclusion.

Since $s$ is a splitting, we have
\[j_1(c) = j_1  s_1 j_1(c) + j_2 s_2 j_1(c), \text{ so that }
j_1(c-s_1j_1(c)) = j_2 s_2 j_1 (c).\]
Thus $j_1(c-s_1j_1(c))\in  \im j_2$,  and so in particular $c-s_1j_1(c) \in \ker d^3$.

Define a map $\rho\colon C^2 \to \ker d^3 \subset C^2$ by $c\mapsto c-s_1j_1(c)$. Furthermore, for all $e \in C^1$,
\[(\rho\circ d^2)(e)=d^2(e)-s_1j_1 d^2(e)=d^2(e)-s_1(0)=d^2(e).\]

Therefore, we can consider the following commutative diagram, in which we use the dual $\rho^*$ of $\rho \colon C^2 \to C^2$, composed with the identification of $C_i$ with its double dual.
\[\xymatrix{ & C_2 \ar[r]^{d_2} & C_1 \ar[r]^{d_1} & C_0 \\
C_3 \ar[r]_{d_3} \ar[ru]^{0} & C_2 \ar[r]_{d_2} \ar[u]_{\rho^*} & C_1  \ar[u]_{\id} \ar[r]_{d_1} & C_0 \ar[u]_{\id}
}\]
We will show that this chain map can be realised as a map $f \colon M^{(3)}\to M^{(2)}$.
To see this, start with the identity map $M^{(2)}\to M^{(2)}$. Since $d_2\circ \rho^* - d_2=(\rho\circ d^2)^* -d_2=(d^2)^*-d_2=0$, we can change this map on the $2$-cells by elements of $\pi_2(M^{(2)})$ according to $\rho^*-\id$, i.e.\ pinch off a 2-sphere from each 2-cell of $M^{(2)}$, and map this 2-sphere to the image of the 2-cell under $\rho^*-\id\colon C_2 \to  \ker d_2 = \pi_2(M^{(2)})$. Since $\rho^*\circ d_3=(d^3\circ \rho)^*=0$, the attaching maps of the $3$-cells are null homotopic under the new map $M^{(2)} \to M^{(2)}$.  Therefore it can be extended to a map $f\colon M^{(3)} \to M^{(2)}$, as desired.

It remains to show that the map $f^*$ induced by $f$ on second cohomology is $s_2$.
Consider the following diagram.
\[\xymatrix @C+0.5cm {C^2/\im d^2 \ar[r]^{s_2} & \ker d^3/\im d^2  \ar @{_{(}->}[rd]^{j_2} & \\
C^2 \ar@{->>}[u]^{j_1} \ar@{->>}[d]_{j_1} \ar[r]^-{\rho = \id -s_1j_1} & \ker d^3 \subset C^2 \ar[d] \ar[r]^{j_1} & C^2/\im d^2 \\
C^2/\im d^2 \ar[r]^{f^*} & \ker d^3/\im d^2. \ar @{^{(}->}[ru]^{j_2} &
}\]
The bottom left square commutes since $f$ was constructed to realise the chain level map $\rho^*$.  The top trapezium commutes by the formula $j_1(c-s_1j_1(c)) = j_2 s_2 j_1 (c)$ from above.  It is straightforward to see that the bottom right triangle commutes.
Since $j_1$ is surjective and $j_2$ injective, it follows from the equality of the top and bottom routes that $f^*=s_2\colon H^2(M^{(2)};\Z\pi)\to H^2(M^{(3)};\Z\pi)\cong H^2(M;\Z\pi)$.
\end{proof}

 Let $X$ be a CW complex model for $B\pi$ with finite $2$-skeleton $K$.  Let $i\colon K\to X^{(3)}$ denote the inclusion of the $2$-skeleton into the $3$-skeleton. Let $\alpha\in\pi_3(M^{(3)})$ denote the class of the attaching map of the (unique) $4$-cell of $M$.

\begin{lemma}
\label{lem:extension3skeleton}
The map $i\colon K\to X^{(3)}$ induces the trivial map $\pi_3(K) \xrightarrow{0} \pi_3(X^{(3)})$. In particular, for every map $f\colon M^{(3)}\to K$ we have $i_*f_*(\alpha)=0 \in \pi_3(K)$, and therefore the map $i\circ f$ can be extended to a map $\widehat{f}\colon M\to X^{(3)}$.
\end{lemma}

\begin{proof}
Since  the maps in Whitehead's sequence (\cref{thm:whitehead-seq}) are natural, we have a commutative diagram:
\[\xymatrix{0=H_4(\wt{K};\Z)\ar[r]\ar[d]&\Gamma(\pi_2(K))\ar[r]^-{\cong}\ar[d]^{i_*} & \pi_3(K)\ar[r]\ar[d]^{i^*} & H_3(\wt{K};\Z)=0\ar[d]\\
0=H_4(\wt{X^{(3)}};\Z)\ar[r]&\Gamma(\pi_2(X^{(3)}))\ar[r] & \pi_3(X^{(3)})\ar[r] & H_3(\wt{X^{(3)}};\Z).}\]
But $\pi_2(X^{(3)}) = 0$, and hence $i_*\colon \pi_3(K)\to\pi_3(X^{(3)})$ is trivial.
\end{proof}

\noindent We obtain the following corollary to \cref{prop:geosplit}.

\begin{cor}\label{cor:map-exists-M3-K-iso-pi1}
There exists a map $f\colon M^{(3)}\to K$ that induces an isomorphism  $\pi_1(M^{(3)}) \to \pi_1(K)$ if and only if $\pri(M)=0$.
\end{cor}

\begin{proof}
Recall that all $2$-complexes with the same fundamental group become homotopy equivalent after wedging with sufficiently many $2$-spheres, by \cref{lem:stable}.  Let $K$ and $K'$ be two such 2-complexes. We have a sequence of maps
\[K' \to K' \vee \bigvee^m S^2 \raiso K \vee \bigvee^{m'} S^2 \to K,\]
where the first map is inclusion, the second map is a homotopy equivalence, and the final map is a collapse map.  The composition induces an isomorphism on fundamental groups, so it is enough to show that there is a map $M^{(3)} \to K'$ for one 2-complex $K'$ with $\pi_1(K')\cong \pi$.  We will use $K'=M^{(2)}$. If $\pri (M) =0$, then the sequence
\[0\to \ker d^3 \to C^2\oplus H^2(M;\Z\pi)\to \coker d^2 = H^2(M^{(2)};\Z\pi) \to 0\]
splits and there is a map $f \colon M^{(3)} \to M^{(2)}$ by \cref{prop:geosplit}, as required.

On the other hand, if there is a map $f \colon M^{(3)} \to K$, then we can extend $f$ to a map $\widehat f\colon M\to X^{(3)}$ by \cref{lem:extension3skeleton}.  If $f$ is an isomorphism on fundamental groups, then so is $\widehat{f}$. By \cref{prop:spin-AHSS} and the discussion preceding it, this implies $\pri(M)=0$.
\end{proof}

From now on in this section we will now only consider manifolds $M$ with $\pri(M)=0$. Let $f\colon M^{(3)}\to K$ be a fixed map that is an isomorphism on fundamental groups, which exists by \cref{cor:map-exists-M3-K-iso-pi1}. Let $\alpha\in\pi_3(M^{(3)})$ denote the class of the attaching map of the $4$-cell of $M$ (we may and will assume that~$M$ has a unique $4$-cell).

By \cref{lem:extension3skeleton}, $f$ can be extended to a map $\widehat{f} \colon M\to X^{(3)}$.
The invariant \[\msec(M) \in H_3(\pi;\Z/2)/\im(d^2_{5,0})\] from \cref{prop:spin-AHSS} is independent of the choice of homotopy of the classifying map $M \to B\pi$ to a map $M\to B\pi^{(3)} = X^{(3)}$.  If we fix a map $\widehat f\colon M\to X^{(3)}$, then we can consider an element $\msec(M,\widehat f)\in H_3(\pi;\Z/2)$, defined as in \cref{lem:arf}. Note that $\msec(M,\widehat f)\in H_3(\pi;\Z/2)$ is independent of the homotopy class of $\widehat f$ viewed as a map $M\to X^{(4)}$.

We briefly recall the construction here, since we will need the details in the near future.
Take a regular preimage $F_i$ of a barycentre for every $3$-cell $e_i$ of $X$. A framing of the normal bundle of this point pulls back to a framing of the normal bundle of the regular preimage $F_i\subseteq D^4\subseteq M$ and thus gives the preimage a framing of the normal bundle. Thus we can consider $[F_i]\in\Omega_1^{\fr} \cong \Omega_1^{\spin}\cong \Z/2$.
Then $\sum_i[F_i]e_i\in \Z/2 \otimes C_3(X)$ lies in $\ker(\Z/2 \otimes C_3(X)\to \Z/2 \otimes C_2(X))$, and so defines a homology class $\msec(M,\widehat f)\in H_3(\pi;\Z/2)$.

\begin{lemma}\label{lem:f-determines-sec}
The element $\msec(M,\widehat{f})\in H_3(\pi;\Z/2)$, where $\widehat{f}\colon M\to X^{(3)}$ is an extension of $i\circ f\colon M^{(3)} \to K \to X^{(3)}$, only depends on the map~$f$.
\end{lemma}

\begin{proof}
Note that any two extensions of $f$ to $\widehat{f}$ only differ by an element of $\pi_4(X^{(3)})$, that is obtained from the two choices of extension on the unique 4-cell of $M$. Consider the diagram
\[\xymatrix{
\pi_4(X^{(3)})\ar[r]\ar[d]^h&\pi_4(X^{(4)})\ar[d]^h_\cong\\
H_4(X^{(3)};\Z\pi)=0\ar[r]&H_4(X^{(4)};\Z\pi),}\]
where $h$ denotes Hurewicz maps, which is commutative by naturality of the Hurewicz homomorphism.  It follows that $\pi_4(X^{(3)})\to \pi_4(X^{(4)})$ is the trivial map, so any two choices of extension $\widehat f$ are homotopic over $X^{(4)}$. Hence $\msec(M,\widehat f)$ only depends on $f$, and not on the choice of~$\widehat f$.
\end{proof}

\begin{defi}
We define \[\msec(M,f):=\msec(M,\widehat{f})\in H_3(\pi;\Z/2),\]
for an extension $\widehat{f} \colon M \to X^{(3)}$ of $i\circ f \colon M^{(3)} \to X^{(3)}$. Such an extension exists by \cref{lem:extension3skeleton} and $\msec(M,f)$ is independent of the choice of extension by \cref{lem:f-determines-sec}.
\end{defi}

Consider the following diagram, which commutes by functoriality of Tate cohomology $\widehat{H}^0(-)$, and by the definition of~$A \colon H_3(\pi;\Z/2) \to \widehat{H}^0(\Sesq(H^2(K;\Z\pi)))$.  We will explain more about the maps in the diagram below.  The proof of \cref{thm:mainsec}, that \cref{cond} implies the Secondary Property, will be based on this diagram.
\begin{equation}\label{eq:master-sec-diagram}
\xymatrix{
\Z\otimes_{\Z\pi} \pi_3(M^{(3)})\ar[d]^{f_*}&H_3(\pi;\Z/2)\ar@{^(->}[r]^{\iota} \ar@{-->}[rdd]^<<<<<<<<<<{A} &\Z/2 \otimes_{\Z\pi} \pi_2(K)\ar[d]^{\Theta_{H^2(K;\Z\pi)}}\\
\Z \otimes_{\Z\pi}\pi_3(K) \ar[r]^-{\cong}_-{S} & \Z \otimes_{\Z\pi}\Gamma(\pi_2(K)) \ar[r] \ar[d]^{\Phi_{H^2(K;\Z\pi)}} & \widehat{H}^0(\pi_2(K) \otimes_{\Z\pi}\pi_2(K)) \ar[d]^{\widehat{H}^{0}(\Phi_{H^2(K;\Z\pi)})}\\
& (\Sesq(H^2(K;\Z\pi)))^{\Z/2} \ar[r]&\widehat{H}^0(\Sesq(H^2(K;\Z\pi)))\\
& (\Sesq(H^2(M;\Z\pi)))^{\Z/2} \ar[r]\ar[u]^{\Sesq(f^*)}&\widehat{H}^0(\Sesq(H^2(M;\Z\pi)))\ar[u]^{\widehat{H}^0(\Sesq(f^*))}}\end{equation}
As promised, some remarks on the diagram are in order.

\begin{enumerate}
\item In the bottom row, we could have written $M^{(3)}$ instead of $M$, since $M^{(3)}$ is the domain of the map $f$.  But attaching a 4-cell makes no difference to second cohomology, so we allow ourselves this slight abuse.
\item Recall that for a set $X$ and a group $G$ that acts on $X$, $X^G$ denotes the fixed points. Thus $(\Sesq(H^2(K;\Z\pi)))^{\Z/2}$ is the $\Z/2$-fixed points, that is the kernel of the $\Z[\Z/2]$-module homomorphism $1-T \colon \Sesq(H^2(K;\Z\pi)) \to \Sesq(H^2(K;\Z\pi))$. In other words, the hermitian forms.
\item Also recall that, by \cref{lem:pi2-H2-dual}, there is a canonical isomorphism \[I \colon \pi_2(K) \raiso (H^2(K;\Z\pi))^*,\] so that we can consider the map $$\Phi_{H^2(K;\Z\pi)} \colon \pi_2(K) \otimes_{\Z\pi}\pi_2(K) \to \Sesq(H^2(K;\Z\pi)).$$
In the diagram, we see the map induced on Tate cohomology by $\Phi_{H^2(K;\Z\pi)}$ in the right hand column, and we see the map restricted to the $\Gamma$ group in the middle column.
\item The map $S :=\Id_{\Z} \otimes \Gamma(\eta)$ is an isomorphism by \cref{lem:S-iso}, where we use that $H_4(\wt{K})=H_3(\wt{K})=0$.
\item The map $\iota \colon H_3(\pi;\Z/2) \to \Z/2 \otimes_{\Z\pi} \pi_2(K)$ is an injection by \cref{lem:iota-inj}.
\end{enumerate}

We will show that the image of $\msec(M,f)\in H_3(\pi;\Z/2)$ in $\widehat{H}^0(\Sesq(H^2(K;\Z\pi)))$ agrees with the image of $\lambda_M\in (\Sesq(H^2(M;\Z\pi)))^{\Z/2}$ in $\widehat{H}^0(\Sesq(H^2(K;\Z\pi)))$. Here the image of $\msec(M,f)$ is $A(\msec(M,f))$ (the map $A$ was introduced in \cref{defn:map-A}). The image of $\lambda_M$ is the intersection form restricted to $H^2(K;\Z\pi)$ along $f^*$, modulo even forms. To show that the images coincide we will show that both of them agree with a third element of $\widehat{H}^0(\Sesq(H^2(K;\Z\pi)))$, namely the image of $[\alpha]\in \Z \otimes_{\Z\pi} \pi_3(M^{(3)})$ in $\widehat{H}^0(\Sesq(H^2(K;\Z\pi)))$; recall that $\alpha \in \pi_3(M^{(3)})$ is the attaching map of the $4$-cell of $M$.  We abuse notation and write $\alpha$ for $1 \otimes \alpha \in \Z \otimes \pi_3(M^{(3)})$.

\begin{lemma}\label{lem:image-sec}
We have an equality
\[\Theta_{H^2(K;\Z\pi)}\circ \iota(\msec(M,f))=[S \circ f_*(\alpha)]\in \widehat{H}^0(\pi_2(K)\otimes_{\Z\pi}\pi_2(K)).\]
\end{lemma}

\begin{proof}
We have that $(M\#(S^2\times S^2))^{(3)}\simeq M^{(3)}\vee S^2\vee S^2$, and under this stabilisation the attaching map of the $4$-cell changes by the Whitehead product of two identity maps $S^2 \to S^2$. Let $\alpha' \in \pi_3(M^{(3)} \vee S^2 \vee S^2)$ denote the new attaching map of the $4$-cell. For any two elements $a,b\in \pi_2(K)$, we can consider the map $f\vee a\vee b \colon M^{(3)} \vee S^2 \vee S^2 \to K$ and obtain
\[S\circ (f\vee a\vee b)_*\alpha'=S\circ f_*\alpha+a\otimes b+b\otimes a\in \Z\otimes_{\Z\pi}\Gamma(\pi_2(K)).\]
Hence
\[[S\circ (f\vee a\vee b)_*\alpha']=[S\circ f_*\alpha]\in \widehat{H}^0(\pi_2(K)\otimes_{\Z\pi}\pi_2(K)).\]
On the other hand, if we extend $f$ and $f\vee a\vee b$ to $X^{(3)}$ via the inclusion $K \to X^{(3)}$, then they become bordant maps over $X^{(3)}$, since $X^{(3)}$ is built from $K$ by adding $3$-cells to kill $\pi_2(K)$. It follows that $\msec(M,f)=\msec(M\#(S^2\times S^2),f\vee a\vee b)$. Thus the lemma holds for $(M,f)$ if and only if it holds true for $(M\#(S^2\times S^2),f\vee a\vee b)$.

Up to elements of the form $a\otimes b+b\otimes a$, every element in $\Gamma(\pi_2(K))$ can be written as $\beta\otimes\beta$ for some $\beta\in\pi_2(K)$. Thus by the previous paragraph we can and will assume that the attaching map of the $4$-cell $\alpha \in \pi_3(M^{(3)})$ is mapped to $\beta\otimes\beta$ under the composition
\[\Gamma(\eta)^{-1} \circ f_*\colon \pi_3(M^{(3)})\to \pi_3(K) \to \Gamma(\pi_2 (K)).\]
It follows from \cref{lem:gamma-eta} that $f_*(\alpha)=\beta\circ\eta$, where $\eta \colon S^3\to S^2$ is the Hopf map.

Then $f$ extends to a map $\widehat f\colon M \to K\cup_\beta D^3$, where the map from the top cell of $M$ to $D^3$ is the cone on $\eta \colon S^3 \to S^2$.

Attach further $3$-cells to turn $K\cup_{\beta} D^3$
into the 3-skeleton $B\pi^{(3)}$ of a model for $B\pi$.  Therefore, we obtain $\msec(M,f)=\msec(M,\widehat f)\in H_3(\pi;\Z/2)$ by the construction given in the proof of \cref{lem:f-determines-sec}, where technically speaking we extend $\widehat{f}$ further by the inclusion $\widehat{f} \colon M \to K\cup_\beta D^3 \to B\pi^{(3)}$.  This is a particular choice of extension $\widehat{f}$, as in the proof of \cref{lem:f-determines-sec}.
Recall that for the definition of $\msec(M,f)$, we take a regular preimage $F_i$ of a barycentre for every $3$-cell $e_i$ of $X^{(3)}$, where $X \simeq B\pi$, and then $\sum_i[F_i]e_i$ determines $\msec(M,\widehat f)\in H_3(\pi;\Z/2)$, where $[F_i]$ is the framed or spin bordism class of $F_i$.

We assert that in the above construction $[F_i]=1$ if and only if $e_i=e_\beta$, the first cell attached along $\beta$. If $e_i\neq e_\beta$, then $F_i$ is empty since $\widehat{f}$ factors through $K \cup_{\beta} D^3 \to B\pi^{(3)}$. To see that $[F_\beta]=1$, note that it was a regular preimage of $C\eta\colon D^4\to D^3$, and thus of $\Sigma\eta$, since $\Sigma\eta \colon S^4 \to S^3$ arises from gluing together two copies of $C\eta \colon D^4 \to D^3$, and we can assume the point whose preimage we take to be in one of the $C\eta$ halves.
Hence $[F_\beta]=1$ follows from $$\Sigma\eta=1\in \pi_1^{st}\cong \Omega_1^{\fr} \cong \Omega_1^{\spin}\cong \Z/2.$$

Hence $\msec(M,f)=[e_\beta]\in H_3(\pi;\Z/2)$. Consider the diagram
\[\xymatrix{
H_3(K;\Z/2)\ar[r]\ar[d]&H_3(\pi;\Z/2)\ar[r]^-{\iota}\ar[d]&\Z/2\otimes_{\Z\pi}\pi_2(K)\ar[d]\\
H_3(K\cup_\beta D^3;\Z/2)\ar[r]&H_3(\pi;\Z/2)\ar[r]^-{\iota}&\Z/2\otimes_{\Z\pi}\pi_2(K\cup_\beta D^3)\\
}\]
Since $[e_\beta]$ comes from $H_3(K\cup_\beta D^3;\Z/2)$ it maps to zero in $\Z/2\otimes_{\Z\pi}\pi_2(K\cup_\beta D^3)$ and thus $\iota([e_\beta])\in \Z/2\otimes_{\Z\pi}\pi_2(K)$ lies in the kernel of the map to $\Z/2\otimes_{\Z\pi}\pi_2(K\cup_\beta D^3)$. This kernel precisely consists of $1\otimes \beta$. Hence $\iota(\msec(M,f))=1\otimes \beta$ and thus
\[\Theta_{H^2(K;\Z\pi)}\circ \iota(\msec(M,f))=[\beta\otimes\beta]=[S \circ f_*(\alpha)].\qedhere\]
\end{proof}

\begin{lemma}\label{lem:image-lambda}
Let $\alpha \in \pi_3(M^{(3)})$ be the attaching map of the $4$-cell of $M$, and let $f \colon M^{(3)} \to K$ be a map that induces an isomorphism on fundamental groups.  We have an equality:
\[[\Phi_{H^2(K;\Z\pi)}(S \circ f_*(\alpha))]=[\Sesq(f^*)(\lambda_M)]\in \widehat H^0(\Sesq(H^2(K;\Z\pi))).\]
\end{lemma}

\begin{proof}
Arguing as in the first two paragraphs of the proof of \cref{lem:image-sec}, it suffices to consider the case $f_*(\alpha)=\beta\circ \eta$ for some  $\beta\in\pi_2(K)$.
Denote the extension of $f \colon M^{(3)} \to K$ to $M\to K\cup_{\beta\circ\eta}D^4$ again by~$f$.
Then for any $x,y \in H^2(K;\Z\pi)$ we have
\begin{align*}
&\lambda_M(PD(f^*(x)),PD(f^*(y)))=\langle f^*(y),PD(f^*(x))\rangle\\
=&\langle y, f_*(f^*(x)\cap [M])\rangle=\langle y, x\cap f_*[M]\rangle.\end{align*}
Denote  the extension of $\beta\colon S^2\to K$ by $\wh\beta\colon \CP^2\to K\cup_{\beta\circ \eta} D^4$. Then $\wh\beta_*[\CP^2]=f_*[M]$ and we have
\[\langle y,x\cap \wh\beta_*[\CP^2]\rangle =\langle \wh\beta^*(y),\wh\beta^*(x)\cap [\CP^2]\rangle
=\wh\beta^*(x)\langle 1,1\cap [\CP^2]\rangle \overline{\wh\beta^*(y)},\]
where $1\in H^2(\CP^2;\Z\pi) \cong \Z\pi$ denotes a chosen generator. For the second equation we used the sesquilinearity of the evaluation.  By the cup product form of $\CP^2$ we have $\langle 1,1\cap [\CP^2]\rangle=1$ and since $\wh\beta$ is given by $\beta$ when restricted to the $2$-skeleton, we obtain
\[\wh\beta^*(x)\langle 1,1\cap [\CP^2]\rangle \overline{\wh\beta^*(y)}=\beta^*(x)\overline{\beta^*(y)}=\Phi_{H^2(K;\Z\pi)}(\beta^*\otimes \beta^*)(x\otimes y).\]
Combining the above equations we obtain
\begin{align*}
    \lambda_M(PD(f^*(x)),PD(f^*(y)))\rangle&=\Phi_{H^2(K;\Z\pi)}(\beta^*\otimes \beta^*)(x\otimes y)\\
    &=\Phi_{H^2(K;\Z\pi)}(S\circ(\beta\circ\eta))(x\otimes y).
\end{align*}
as claimed.
\end{proof}

\begin{thm}
\label{thm:lambdasec}
We have an equality
\[[\Sesq(f^*)(\lambda_M)]=A(\msec(M,f))\in\widehat{H}^0(\Sesq(H^2(K;\Z\pi))).\]
\end{thm}

\begin{proof}
  The theorem follows from combining  \cref{lem:image-sec} and \cref{lem:image-lambda}, together with a straightforward diagram chase in diagram \eqref{eq:master-sec-diagram}.
\end{proof}

\begin{cor}\label{cor:sec}~
\begin{enumerate}
  \item\label{cor:sec1} If $\ker A\subseteq \im (Sq^2\circ\red_2 \colon H_5(\pi;\Z)\to H_3(\pi;\Z/2))$, then $\msec(M)=0$ if and only if there exists $f \colon M^{(3)} \to K$ as above with
\[0=[\Sesq(f^*)(\lambda_M)]\in\widehat{H}^0(\Sesq(H^2(K;\Z\pi))).\]
\item\label{cor:sec2} If $\im (Sq^2\circ\red_2\colon H_5(\pi;\Z)\to H_3(\pi;\Z/2))\subseteq \ker A$, then the class
\[[\Sesq(f^*)(\lambda_M)]\in\widehat{H}^0(\Sesq(H^2(K;\Z\pi)))\]
is independent of the choice of $f$.
\end{enumerate}
\end{cor}

\begin{proof}~
\begin{enumerate}
  \item If $\msec(M)=0$, then there  exists an $f$ with $\Sesq(f^*)(\lambda_M)=0$.
This is because $\msec(M)=0$ implies that there is a map $f \colon M \to B\pi^{(2)}$, and we can take $K$ for $B\pi^{(2)}$.
The intersection form vanishes on $f^*(H^2(K;\Z\pi))$ by naturality of the evaluation, because $H_*(K;\Z)=0$.

Conversely, if there exists a map $f$ with \[0=[\Sesq(f^*)(\lambda_M)]\in\widehat{H}^0(\Sesq(H^2(K;\Z\pi))),\]
then by \cref{thm:lambdasec}, $\msec(M,f)\in \ker A\subseteq \im (Sq^2\circ\red_2\colon H_5(\pi;\Z)\to H_3(\pi;\Z/2))$. Since $Sq^2\circ \red_2$ is the second differential in the Atiyah-Hirzebruch spectral sequence computing $\Omega_4^{\spin}(B\pi)$ (see \cref{section:AHSS}), we have $\msec(M)=0\in H_3(\pi;\Z/2)/\im(Sq^2\circ\red_2)$.
\item Given two choices $f,f' \colon M^{(3)} \to K$, then there exists an $x\in \im (Sq^2\circ \red_2)$ with $\msec(M,f)=\msec(M,f')+x$ since $Sq^2\circ\red_2$ is the second differential of the Atiyah-Hirzebruch spectral sequence and $\msec(M)$ is well-defined in $H_3(\pi;\Z/2)/\im(S^2\circ \red_2)$.  Thus by \cref{thm:lambdasec} we have
\begin{align*}
[\Sesq(f^*)(\lambda_M)] &=A(\msec(M,f))=A(\msec(M,f))+A(x)=A(\msec(M,f'))\\ &=[\Sesq((f')^*)(\lambda_M)].
\end{align*}
Here we used that $A$ is a homomorphism, and $A(x)=0$ by assumption.\qedhere
\end{enumerate}
\end{proof}

We can now apply \cref{cor:sec} to prove \cref{thm:mainsec}, which says that under \cref{cond} the \cref{secondary-property} holds, which we recall here in an expanded form for the convenience of the reader:

\noindent \emph{If $M$ is a spin manifold with $\pi_1(M)\raiso \pi$ and $\pri (M)=0$, then $\msec (M)=0$ if and only if there exists a splitting
\[s=(s_1,s_2)\colon \coker d_3 \to C_2 \oplus H_2(M;\Z\pi)\]
of the extension from \cref{thm:prim-obstruction-thm} such that $\lambda_M$ restricted to the image of $s_2$ is even.
The restriction of $\lambda_M$ for one splitting $s$ is even if and only if it is even for every splitting.}

\emph{Moreover, for two such manifolds $M,M'$ with the same fundamental group we have $\msec (M)=\msec (M')$ if and only if the restrictions of $\lambda_M$ and $\lambda_{M'}$ are isomorphic modulo even forms for some splittings $s$ and $s'$.
}

\begin{proof}[Proof of \cref{thm:mainsec}]
\cref{cond} says precisely that the hypotheses of \cref{cor:sec}~(\ref{cor:sec1}) and (\ref{cor:sec2}) hold.
If $\pri(M)=0$, then there is a splitting $s \colon \coker d_3 \to C_2 \oplus H_2(M;\Z\pi)$.

  We want to apply \cref{prop:geosplit}  to obtain a map $f\colon M^{(3)} \to K$ for a 2-complex $K$ with $\pi_1(K) \cong \pi$.
  However note that \cref{prop:geosplit} was formulated for the dual exact sequence
\[0\to \ker d^3\to C^2\oplus H^2(M;\Z\pi)\to \coker d^2\to 0.\]
So to remedy this, for a splitting $s=(s_1,s_2)\colon \coker d_3\to C_2\oplus H_2(M;\Z\pi)$ of the sequence
\[0\to \ker d_2\to C_2\oplus H_2(M;\Z\pi)\to \coker d_3\to 0\]
from \cref{thm:prim-obstruction-thm}, we apply \cref{prop:geosplit} with the dual handle decomposition of~$M$.  Recall that $M^d$ denotes $M$ with the dual handle decomposition. We obtain a map $f\colon (M^d)^{(3)}\to (M^d)^{(2)}$  with $s_2=f^*\colon H^2((M^d)^{(2)};\Z\pi)\to H^2(M^d;\Z\pi)$ by \cref{prop:geosplit}, and use the identifications $H^2((M^d)^{(2)};\Z\pi) \cong \coker d_3$ and $H^2(M^{d};\Z\pi) \cong H_2(M;\Z\pi)$.

Now by \cref{cor:sec}~(\ref{cor:sec1}), $\Sesq(f^*)(\lambda_M)=\lambda_M|_{s_2(\coker d_3)}$ is even (i.e.\ zero in the Tate cohomology $\widehat{H}^0(\Sesq(H^2(K;\Z\pi)))$) if and only if $\msec(M)=\msec(M^d)=0$.   By \cref{cor:sec}~(\ref{cor:sec2}), $\Sesq(f^*)(\lambda_M)$ is independent of the choice of $f$, and is therefore independent of the choice of splitting~$s$.

The last sentence of the Secondary Property, that $\msec (M)=\msec (M')$ if and only if the restrictions of $\lambda_M$ and $\lambda_{M'}$ are equal for some choices of splitting $s$ and $s'$, also follows from \cref{thm:lambdasec}.  If $\msec(M)=\msec(M')$, then the difference lies in $\im (Sq^2\circ \red_2) =\ker A$, so for some choice of $f,f'$ we have $A(\msec(M,f))=A(\msec(M',f'))$.  Thus  $[\Sesq(f^*)(\lambda_M)] = [\Sesq((f')^*)(\lambda_{M'})]$ by \cref{thm:lambdasec}.  On the other hand, $[\Sesq(f^*)(\lambda_M)] = [\Sesq((f')^*)(\lambda_{M'})]$ implies that $A(\msec(M,f))=A(\msec(M',f'))$, so $\msec(M,f)-\msec(M',f') \in \ker A = \im (Sq^2\circ \red_2)$.  Thus $\msec(M) = \msec(M') \in H_3(\pi;\Z/2)/\im d^2_{5,0}$.
\end{proof}

\chapter{The tertiary obstruction}\label{sec:tertiary-obstruction}

Let $M$ be a spin $4$-manifold with an identification $\pi_1(M) \raiso \pi$, such that $\pri(M)=0$ and $\msec (M)=0$.
The identification determines a classifying map $c \colon M\to B\pi$ uniquely up to homotopy.

In this chapter, more precisely in \cref{thm:ter-alg-condition}, we will give an algebraic criterion on $\pi$, under which we can identify the tertiary obstruction $\ter(M)$ with a~$\tau$ invariant defined on a  subset of $\pi_2(M)$.
The reader should recall the definition of the Kervaire--Milnor invariant $\tau$ from \cref{section:alg-ter-intro}.

\section[The tertiary obstruction and the Kervaire--Milnor invariant]{Relating the tertiary obstruction with the Kervaire--Milnor invariant}

For the well-definedness of the upcoming \cref{def:tauM}, we have to recall the following statements from \cite{KPT21} and \cite{KLPT15}.
Let $K$ be a finite connected 2-complex homotopy equivalent to $M^{(2)}$.
We view $K$ as a 2-skeleton for $B\pi$ and denote the inclusion by $i\colon K\to B\pi$.


\begin{lemma}[{\cite[Lemma~4.2]{KPT21}}]
\label{lem:char}
Let $f\colon M^{(3)}\to K$ be a map that induces the given isomorphism on fundamental groups. Let $j\colon M^{(3)}\to M$ be the inclusion of the $3$-skeleton. For every $\phi\in \ker \Sq^2\subseteq H^2(B\pi;\Z/2) = H^2(B\pi;\Z/2)$ and every lift $\phi'\in H^2(K;\Z\pi)$ of $i^*\phi \in H^2(K;\Z/2)$, the element $PD\circ (j^*)^{-1}\circ f^*(\phi') \in H_2(M;\Z\pi) \cong \pi_2(M)$ is $\RPT$-characteristic.
\end{lemma}

In the statement of \cref{lem:char} we used the following sequence of maps:
\[\phi' \in H^2(K;\Z\pi) \xrightarrow{f^*} H^2(M^{(3)};\Z\pi) \xrightarrow{(j^*)^{-1}} H^2(M;\Z\pi) \xrightarrow{PD} H_2(M;\Z\pi) \cong \pi_2(M).\]

\begin{lemma}[{\cite[Lemma~8.3~and~Lemma~8.4]{KLPT15}}]
\label{lem:augtau}
\begin{enumerate}
\item  Let $x,y\in\pi_2(M)$ be such that $\lambda(x,y)=0$, $\mu(x)$ and $\mu(y)$ are trivial, and $x$ is $\RPT$-characteristic. Then for every element $\kappa\in \ker(\Z\pi\to \Z/2)$, we have $\tau(x)=\tau(x+\kappa y)\in\Z/2$.
\item Let $Y$ be a finite $2$-dimensional CW complex with fundamental group $\pi$. Every element in the kernel of $H^2(Y;\Z\pi)\to H^2(Y;\Z/2)$ can be written as $\sum_{i=1}^n\kappa_ix_i$ with $\kappa_i\in \ker(\Z\pi\to \Z/2)$ and $x_i\in H^2(Y;\Z\pi)$.
\end{enumerate}
\end{lemma}

\begin{remark}
	In \cite{KLPT15} we asserted that the invariant $\tau_M$ is well-defined, and the results of the next lemma, for $S^2$-characteristic spheres. But we ought to have required the stronger $\RPT$-characteristic. The results of that paper are for torsion-free groups, where the two notions coincide. But the proof of the previous lemma goes through unchanged.
\end{remark}

As in \cref{sec:sec}, in this section we will consider the dual sequence
\begin{equation}
	\label{eq:dualsequence}
	0 \to \ker d^3 \to C^2 \oplus H^2(M;\Z\pi) \to \coker d^2 \to 0
\end{equation}
to the sequence in \cref{thm:prim-obstruction-thm}.
As noted previously, this can be obtained from the exact sequence in \cref{thm:prim-obstruction-thm} by working with the dual handle decomposition~$M^d$.
Recall that we say a map $f \colon M^{(3)} \to M^{(2)}\simeq K$ realises a splitting $s \colon \coker d^2 \to C^2 \oplus H^2(M;\Z\pi)$ if $s_2 \colon \coker d^2 \to H^2(M;\Z\pi)$ coincides with  the map \[\coker d^2 \cong H^2(K;\Z\pi) \xrightarrow{f^*} H^2(M^{(3)};\Z\pi) \cong H^2(M;\Z\pi).\]
Recall that $\RC \subseteq \pi_2(M) \cong H_2(M;\Z\pi)$ denotes the subset of $\RPT$-characteristic elements $\alpha$ with $\mu(\alpha)=0$.  The following diagram should help to read the upcoming definition, whose aim is to define a map $\tau_{M,s} \colon \ker \Sq^2 \to \Z/2$.

\[\xymatrix{H^2(K;\Z\pi) \ar[dd]_-{p} \ar[r]^{\cong} & \coker d^2 \ar[r]^{s_2}  & H^2(M;\Z\pi) \ar[r]^{PD}& H_2(M;\Z\pi)   \\
 & (PD \circ s_2)^{-1}(\RC) \ar[rr]^-{(PD \circ s_2)|} \ar[d]_{p|} \ar @{^{(}->}[u] & & \RC \ar @{^{(}->}[u] \ar[d]^{\tau} \\
 H^2(K;\Z/2) & \ar@{_{(}->}[l] p\big((PD\circ s_2)^{-1}(\RC)\big) \ar[rr]^-{\tau'_{M,s}} && \Z/2 \\
  H^2(\pi;\Z/2) \ar[u]^{i^*}  & \ker \Sq^2  \ar @{_{(}->} [l] \ar[u]^{i^*|} \ar[urr]_{\tau_{M,s}} &&
  }\]

\begin{defi}
\label{def:tauM}
Let $s=(s_1,s_2)\colon \coker d^2\to C^2\oplus H^2(M;\Z\pi)$ be a section of the sequence from \cref{thm:prim-obstruction-thm} such that $\lambda_M$ vanishes on $\im s_2$ (see \cref{remark-section-as-required} below for why such a splitting exists stably whenever $\msec(M)=0$).
It follows from \cref{lem:augtau} that the map \[\tau \circ PD \circ s_2\colon (PD \circ s_2)^{-1}\RC \to \Z/2\] factors through the restriction of $p\colon H^2(K;\Z\pi)\to H^2(K;\Z/2)$ to $(PD \circ s_2)^{-1}(\RC)$.
Denote the map $p((PD \circ s_2)^{-1}\RC)\to \Z/2$ arising in this factorisation by $\tau_{M,s}'$.
Let $i^*\colon H^2(\pi;\Z/2)\to H^2(K;\Z/2)$ be the map induced by inclusion of the $2$-skeleton. Note that \[i^*(\ker \Sq^2)\subseteq p((PD \circ s_2)^{-1}(\RC))\] by \cref{lem:char}, since by \cref{prop:geosplit} every splitting $s$ can be realised by a map $f \colon M^{(3)} \to M^{(2)}$. Denote \[\tau_{M,s}:=\tau'_{M,s}\circ i^*| \colon \ker \Sq^2 \to \Z/2.\]
\end{defi}

\noindent We quickly make the following observation.

\begin{lemma}\label{lem:tau-0=0}
We have that $\tau_{M,s}(0)=0$ for any section $s$.
\end{lemma}

\begin{proof}
Computing $\tau_{M,s}(0)$ involves computing $\tau(S)$, where $S$ is an embedded sphere contained in a single chart.  But there are no intersections, therefore no Whitney discs, and so~$\tau$ evidently vanishes.
\end{proof}

For the convenience of the reader, we recall that the \emph{\cref{tertiary-property}} requires that for every spin $4$-manifold with $\pi_1(M) \raiso \pi$, and $\pri(M) = \msec(M) =0$, we have the following, where we have switched to the dual version of the short exact sequence from \cref{thm:prim-obstruction-thm} (see the preamble to \cref{subsection:proof-thm-mainsec}, just before \cref{prop:geosplit}).

\noindent \emph{
For every section \[s=(s_1,s_2)\colon \coker d^2\to C^2\oplus H^2(M;\Z\pi)\]
of the $($dual version of the$)$ short exact sequence from \cref{thm:prim-obstruction-thm}
with $\lambda_M|_{\im PD \circ s_2}\equiv 0$, the map $\tau_{M,s}|_{\ker(d_{5,0}^3)^*}$ is a homomorphism, and $\omega(\ter(M))=\tau_{M,s}|_{\ker(d_{5,0}^3)^*}$, where $\omega$ is the composition:
\begin{align*}
\omega \colon H_2(\pi;\Z/2)/\im (d^2_{4,1}, d^3_{5,0}) &\xrightarrow{\cong} \Hom_{\Z/2}(H^2(\pi;\Z/2),\Z/2)/\im (d^2_{4,1}, d^3_{5,0})\\ &\xrightarrow{\cong} \Hom_{\Z/2}(\ker(d_{5,0}^3)^*,\Z/2).
\end{align*}
}

\begin{lemma}\label{remark-section-as-required}
	Assume that $\pri(M)=\sec(M)=0$. Then there exists a section $s=(s_1,s_2)\colon \coker d^2\to C^2\oplus H^2(M;\Z\pi)$ of \eqref{eq:dualsequence} such that $\lambda_M$ restricted to the image of $PD\circ s_2$ is trivial.
\end{lemma}

This lemma tells us that a section of the form required for the formulation of the Tertiary Property exists.

\begin{proof}
	We view $M^{(2)}$ as a 2-skeleton for $B\pi$. Since $\pri(M)=\sec(M)=0$, by \cref{section:AHSS} there exists, after stabilisation, a map $f\colon M\to M^{(2)}$ inducing an isomorphism on fundamental groups. Since $f$ is an isomorphism on $\pi_1$, we can homotopy it on the 1-skeleton to be the identity. Extending this homotopy to all of $M$, we can assume that the map $f$ is the identity on the 1-skeleton. On cochains the map $f$ then induces a map $f^2\colon C^2\to \ker d^3\subseteq C^2$ such that $f^2\circ d^2=d^2$. Let $s_1':=(\id-f^2)\colon C^2\to C^2$ and $s_2'\colon C^2\xrightarrow{f^2} \ker d^3\twoheadrightarrow H^2(M;\Z\pi)$. Since $f^2\circ d^2=d^2$, the maps $s_1'$ and $s_2'$ factor through $\coker d^2$ inducing a map $s=(s_1,s_2)\colon \coker d^2\to C^2\oplus H^2(M;\Z\pi)$. Then $s$ is a section since $(\id-f^2)+f^2=\id$. It remains to show that $\lambda_M$ restricted to the image of $PD\circ s_2$ is trivial. By definition, $s_2=f^2\colon H^2(M^{(2)};\Z\pi)\to H^2(M;\Z\pi)$. For $x,y\in H^2(M^{(2)};\Z\pi)$ we have
	\begin{align*}
  \lambda_{M}(PD \circ f^*(x),PD \circ f^*(y)) & =\langle f^*(y),f^*(x)\cap [M]\rangle=\langle y,f_*(f^*(x)\cap[M])\rangle \\ &=\langle y,x\cap f_*[M]\rangle.
  \end{align*}
	Since $M^{(2)}$ is 2-dimensional, $f_*[M]=0$ and thus $\lambda_M$ restricted to the image of $PD\circ s_2$ is trivial as needed.
\end{proof}

\begin{remark}
\cref{lem:tau-0=0} implies that the \cref{tertiary-property} holds when $\Hom_{\Z/2}(\ker(d_{5,0}^3)^*,\Z/2) =0$.
\end{remark}

\begin{remark}
  In order to prove the \cref{tertiary-property} when the group  in which $\ter(M)$ lives, $\Hom_{\Z/2}(\ker(d_{5,0}^3)^*,\Z/2)$, is nontrivial, we will show the following.  Recall that $d^2_{4,1}$ is dual to  $\Sq^2 \colon H^2(\pi;\Z/2) \to H^4(\pi;\Z/2)$.
\begin{enumerate}
\item The map $\tau_{M,s}$ defines a homomorphism $\ker \Sq^2 \to \Z/2$.
\item The image of $\tau_{M,s}$ under the map \[\Hom(\ker \Sq^2;\Z/2)\cong H_2(\pi;\Z/2)/ \im(d^2_{4,1}) \to H_2(\pi;\Z/2)/\im(d^2_{4,1},d^3_{5,0})\] is independent of $s$.
\item The image of $\tau_{M,s}$ under the above map agrees with $\ter(M)$.
\end{enumerate}
We will need \cite[Theorem~1.5]{KPT21}, which we recall next, for convenience. This shows that the above  statements hold, in the case that the splitting $s$ comes from a map $M\to M^{(2)}$. \cref{cor:conjter} gives an algebraic condition that implies $s$ is geometrically realised in this manner.
\end{remark}

\begin{proposition}[{\cite[Theorem~1.5]{KPT21}}]
\label{lem:ter}
Suppose that there is a map $f\colon M\to K$ to a 2-complex $K$ that is an isomorphism on fundamental groups.  Let $i \colon K \to B\pi$ be a classifying map.
 \begin{enumerate}[(i)]
\item\label{it:main0} For each $\phi\in \ker (\Sq^2 \colon H^2(\pi;\Z/2) \to H^4(\pi;\Z/2))$, there exists a lift $\widehat{\phi} \in H^2(K;\Z\pi)$ of $i^*(\phi) \in H^2(K;\Z/2)$.
   \item\label{it:main1}  The element $PD(f^*(\widehat{\phi})) \in H_2(M;\Z\pi) \cong \pi_2(M)$ is $\RPT$-characteristic and has trivial self-intersection number, so that the Kervaire-Milnor invariant $\tau(PD\circ f^*(\widehat{\phi})) \in \Z/2$ is well-defined.
 \item\label{it:main2} The map
 \begin{align*}
   \tau_{M,f} \colon \ker \Sq^2 &\to \Z/2 \\
   \phi & \mapsto \tau(PD\circ f^*(\widehat{\phi}))
 \end{align*}
is well-defined $($i.e.\ is independent of the choice of $\widehat{\phi})$ and is a homomorphism.
\item\label{it:main3} Under the map
\[\Hom(\ker \Sq^2,\Z/2) \xrightarrow{\cong} H_2(\pi;\Z/2)/\im \Sq_2 \to H_2(\pi;\Z/2)/\im (d^2_{4,1},d^3_{5,0}),  \]
$\tau_{M,f}$ is sent to $\ter(M)$.  In particular the image of $\tau_{M,f}$ under this map is independent of the choices of $f$ and~$K$.
 \end{enumerate}
\end{proposition}

This shows that, as we are aiming for, $\tau_{M,s}$ agrees with $\ter(M)$, when there is a map $M \to K$ inducing the given splitting.  Now we give an algebraic criterion guaranteeing that such a map $M \to K$ exists. Note that if we did not want to show independence on the choice of splitting, we would be done.  But for our obstructions to be computationally useful, it ought to be possible to choose any splitting.

\begin{defi}
Let $N$ be a $\Z\pi$-module that is free as a $\Z$-module.  Then
\[\varpi_N\colon \Z\otimes_{\Z\pi}\Gamma(N)\to N \otimes_{\Z\pi}N\]
denotes the map induced from the inclusion $\omega \colon \Gamma(N)\to N \otimes_\Z N$ by tensoring with $\Z$ over $\Z\pi$.
As per our convention, we view the first $N$ as a right $\Z\pi$-module using the involution.
Here $\Z\otimes_{\Z\pi} (N \otimes_{\Z} N) \cong N \otimes_{\Z\pi} N$, since $\Z\pi$ acts on $N \otimes_{\Z} N$ on the left by the diagonal action.
\end{defi}

\begin{theorem}
	\label{cor:conjter}
Let $M$ be a spin $4$-manifold with $\pri(M)=0 = \msec(M)$, and let $s \colon \coker d^2 \to C^2 \oplus H^2(M;\Z\pi)$ be a splitting such that $\lambda_M|_{\im(PD \circ s_2)} \equiv 0$ vanishes.
If \[\Phi_{H^2(M^{(2)};\Z\pi)}\circ \varpi_{\pi_2(M^{(2)})}\colon \Z\otimes_{\Z\pi} \Gamma(\pi_2(M^{(2)})) \to \Sesq(H^2(M^{(2)};\Z\pi))\] is injective, then the map $\tau_{M,s}$ is a homomorphism, and maps to $\ter(M)$ under \[\Hom(\ker(\Sq^2),\Z/2) \cong \coker(\Sq_2) \cong E^3_{2,2} \rightarrow E^{\infty}_{2,2} = H_2(\pi;\Z/2)/\im (d^2_{4,1}, d^3_{5,0}).\]
Moreover, the image of $\tau_{M,s}$ in $E^\infty_{2,2}$ does not depend on the choice of $s$.
\end{theorem}

\begin{proof}
By \cref{prop:geosplit}, the splitting $s \colon \coker d^2 \to C^2 \oplus H^2(M;\Z\pi)$ can be realised as the splitting induced from a map $f \colon M^{(3)} \to M^{(2)}$, where $f_*$ is the identity on $\pi_1$.  Here a splitting is said to be realised by $f$ if $s_2= f^* \colon H^2(M^{(2)};\Z\pi) \to H^2(M;\Z\pi)$ after identifying $H^2(M^{(3)};\Z\pi) = H^2(M;\Z\pi)$.

The right hand side of the equation in \cref{lem:image-lambda} vanishes, since the intersection form vanishes on the $\coker d^2 = H^2(M^{(2)};\Z\pi)$ summand of $\pi_2(M)$.  Thus the left hand side of the equation in \cref{lem:image-lambda} vanishes, and then by injectivity of $\Phi_{H^2(M^{(2)};\Z\pi)} \circ \varpi_{\pi_2(M^{(2)})}$, we have that $f_*(\alpha)=0 \in \Z \otimes_{\Z\pi} \pi_3(M^{(2)})$; recall that $\alpha$ denotes the attaching map of the $4$-handle of~$M$. Here we also use the identification $S \colon \Z\otimes_{\Z\pi}\pi_3(M^{(2)}) \raiso \Z\otimes_{\Z\pi} \Gamma(\pi_2(M^{(2)}))$ from diagram~\ref{eq:master-sec-diagram}.

We claim that the images $x_n\in\Z\pi$ of $1\in\Z\pi$ under $\Z\pi\cong C_4\xrightarrow{d_4} C_3\cong \Z\pi^n\xrightarrow{p_n}\Z\pi$, where $p_n\colon \Z\pi^n\to \Z\pi$ are the projections, generate the augmentation ideal $I\pi\subseteq \Z\pi$. This can be seen as follows. By Poincar\'e duality the cokernel of the dual map $C^3\xrightarrow{d^4}C^4$ is $\Z$ and hence the image of $d^4$ is the augmentation ideal. The image is generated by $\{\overline{x_n}\}$ (when we view $C^4$ as a left module using the involution on $\Z\pi$).  By dualising, we obtain the claim.

By changing $f\colon M^{(3)}\to M^{(2)}$ on the $3$-cells, it follows from the claim that we can change $f_*(\alpha)\in \pi_3(M^{(2)})$ by elements of $I\pi\otimes_{\Z\pi}\pi_3(M^{(2)})$. Thus an extension $f'\colon M\to M^{(2)}$ exists if and only if $f_*(\alpha)=0\in \Z\otimes_{\Z\pi}\pi_3(M^{(2)})$. But that is exactly what we have, so such a map $f' \colon M \to M^{(2)}$ indeed exists.

Therefore $s$ is realised geometrically, not just in the sense that we have $s_2 = f^* \colon H^2(M^{(2)};\Z\pi) \to H^2(M^{(3)};\Z\pi)$, but now in the stronger sense that $s_2$ is induced from a map $f \colon M \to M^{(2)}$.

By \cref{lem:ter}, we therefore have $\tau_{M,s}$ is independent of $s$ and is a homomorphism, i.e.\ is an element of $\Hom(\ker \Sq^2,\Z/2)$.
Moreover $\tau_{M,s}$ maps to $\ter(M)$ under the composition in the statement of the theorem, completing the proof.
%
\end{proof}

Next, we want to express the injectivity condition of \cref{cor:conjter} as a property of the group $\pi$ alone.

\begin{lemma}
\label{lem:iotasum}
Let $N,N'$ be $\Z\pi$-modules that are free as $\Z$-modules. Then $\varpi_{N\oplus N'}$ is injective if and only if $\varpi_N$ and $\varpi_{N'}$ are both injective.
\end{lemma}

\begin{proof}
By \cref{lem:gammasum}, the inclusion $\Gamma(N\oplus N')\to (N\oplus N')\otimes_\Z(N\oplus N')$ is the direct sum of the inclusion $\Gamma(N)\to N\otimes_{\Z} N$, the diagonal map $(1+T)\colon N\otimes_\Z N'\to (N\otimes_\Z N')\oplus (N'\otimes_\Z N)$ and the inclusion $\Gamma(N')\to N'\otimes_\Z N'$. Thus, $\varpi_{N\oplus N'}$ is the direct sum of $\varpi_N$, $\varpi_{N'}$ and the diagonal map $(1+T)\colon N\otimes_{\Z\pi} N'\to (N\otimes_{\Z\pi} N')\oplus (N'\otimes_{\Z\pi} N)$. Since the diagonal map is always injective, the lemma follows.
\end{proof}

\begin{lemma}
\label{lem:inj-iota-proj}
Let $P$ be a finitely generated projective left $\Z\pi$-module. Then $\varpi_P$ is injective.
\end{lemma}

\begin{proof}
We begin with the case $P=\Z\pi$. The module $\Z\pi\otimes_\Z\Z\pi$ is a free $\Z\pi$-module with basis $\{1\otimes g\}_{g\in \pi}$. The action of $\Z\pi$ is the diagonal action. The module $\Gamma(\Z\pi)$ is generated as a subset of $\Z\pi\otimes_\Z\Z\pi$ by the elements $1\otimes 1$ and $1\otimes g+g\otimes 1$. Since $1\otimes g+g\otimes 1=1\otimes g+ g(1\otimes g^{-1})$ it follows that there is a homomorphism
\[\Z\pi\otimes_\Z\Z\pi \to \bigoplus_{g\in\pi, g^2\neq 1}\Z\pi\oplus\bigoplus_{g^2=1, g\neq 1}\Z\pi/(1+g).\]

The homomorphism above sends $1 \otimes g$ to the generator of the $g$ summand, with kernel the image of $\Gamma(\Z\pi)$.  If $g^2=1$, then $1 \otimes g + g \otimes 1 = (1+g)(1 \otimes g)$, whence the quotient by $(1+g)$.
Here we think of $(1+g)$ as a left ideal of the ring $\Z\pi$, and form the quotient \emph{ring}.  Then we consider the ring $\Z\pi/(1+g)$ as a left $\Z\pi$-module.
Therefore we have a short exact sequence
\[0\to \Gamma(\Z\pi)\to \Z\pi\otimes_\Z\Z\pi\to \bigoplus_{g\in\pi, g^2\neq 1}\Z\pi\oplus\bigoplus_{g^2=1, g\neq 1}\Z\pi/(1+g)\to 0.\]
Tensor this with $\Z$ over $\Z\pi$ and apply the $6$-term exact sequence to obtain
\begin{align*}
&\Tor_1^{\Z\pi}\Big(\Z,\bigoplus_{g\in\pi, g^2\neq 1}\Z\pi\oplus\bigoplus_{g^2=1, g\neq 1}\Z\pi/(1+g)\Big)\to \Z\otimes_{\Z\pi}\Gamma(\Z\pi) \\
 \xrightarrow{\varpi_{\Z\pi}}
& \Z\pi\otimes_{\Z\pi}\Z\pi\to \bigoplus_{g\in\pi, g^2\neq 1}\Z\oplus\bigoplus_{g^2=1, g\neq 1}\Z/2\to 0.
\end{align*}
A free $\Z\pi$-module resolution of $\Z\pi/(1+g)$, where $g$ has order two, is given by
\[\ldots \xrightarrow{1+g} \Z\pi \xrightarrow{1-g} \Z\pi \xrightarrow{1+g}\Z\pi\xrightarrow{} \Z\pi/(1+g) \to 0.\]
Tensor this with $\Z$ to obtain
\[\ldots \xrightarrow{2} \Z\xrightarrow{0}\Z\xrightarrow{2}\Z\xrightarrow{}\Z/2\to 0.\]
The first homology of the corresponding deleted resolution vanishes.  As $\Tor_1^{\Z\pi}(-,-)$  commutes with direct sums in the second factor, we obtain
\[\Tor_1^{\Z\pi}\Big(\Z,\bigoplus_{g\in\pi, g^2\neq 1}\Z\pi\oplus\bigoplus_{g^2=1, g\neq 1}\Z\pi/(1+g)\Big)=0.\]
This concludes the proof in the case $P=\Z\pi$. For a general finitely generated free $\Z\pi$-module the lemma follows from this and \cref{lem:iotasum} by induction on the rank.

If $P$ is finitely generated projective, there exists $P'$ such that $P\oplus P'$ is finitely generated free. So the injectivity of $\varpi_P$ follows from the injectivity of $\varpi_{P\oplus P'}$, again using \cref{lem:iotasum}.
\end{proof}

\begin{cor}
\label{cor:injiota}
Let $N$ be a left $\Z\pi$-module that is free as a $\Z$-module.
Then the map $\Phi_N\circ \varpi_{N^*}$ is injective if and only if $\Phi_{N\oplus\Z\pi}\circ \varpi_{N^*\oplus \Z\pi^*}$ is injective.
\end{cor}

\begin{proof}
By \cref{lem:phisum} and the definition of $\varpi$, the map $\Phi_{N\oplus \Z\pi}\circ \varpi_{N^*\oplus\Z\pi^*}$ is the direct sum of the maps $\Phi_N\circ \varpi_{N^*}$, $\Phi_{\Z\pi}\circ\varpi_{\Z\pi^*}$ and $(\Phi_{N,\Z\pi^*}\oplus\Phi_{\Z\pi,N^*})\circ(1+T)$.  It suffices to argue that the last two maps are injective.
 Note that $(1+T)\colon N^*\otimes_{\Z\pi} \Z\pi^*\to (N^*\otimes_{\Z\pi} \Z\pi^*)\oplus (\Z\pi^*\otimes_{\Z\pi} N^*)$ is injective, and $\varpi_{\Z\pi^*}$ is injective by \cref{lem:inj-iota-proj}. Moreover \cref{lem:injphi} says that $\Phi_{\Z\pi}$, $\Phi_{N,\Z\pi^*}$ and $\Phi_{\Z\pi,N^*}$ are injective.  It follows that $\Phi_{\Z\pi}\circ\varpi_{\Z\pi^*}$ and $(\Phi_{N,\Z\pi^*}\oplus\Phi_{\Z\pi,N^*})\circ(1+T)$ are injective, as required.
\end{proof}

We can now give the desired statement of a condition implying the \cref{tertiary-property} that depends only on the group~$\pi$.

\begin{theorem}\label{thm:ter-alg-condition}
Let $K$ be a finite 2-complex with fundamental group $\pi$.  Suppose that \[\Phi_{H^2(K;\Z\pi)}\circ\varpi_{\pi_2(K)}\colon \Z\otimes_{\Z\pi} \Gamma(\pi_2(K)) \to \Sesq(H^2(K;\Z\pi))\] is injective.
Then $\pi$ has the \cref{tertiary-property}.
\end{theorem}

\begin{proof}
The theorem with $K=M^{(2)}$ follows from \cref{cor:conjter}. Since any two choices of $K$ are homotopy equivalent after wedging with enough copies of $S^2$, it remains to check that the injectivity condition is preserved under wedging with $S^2$.
But this follows from \cref{cor:injiota}.
Therefore the injectivity condition is independent of the choice of $2$-complex $K$, as desired.
\end{proof}

\section{The Tertiary Property for finite groups}

We will prove the following corollary of \cref{thm:ter-alg-condition}, which will enable us to verify that many finite groups have the \cref{tertiary-property}.

\begin{cor}
\label{cor:ter-finite-group}
Let $\pi$ be a finite group and let $K$ be a finite $2$-dimensional CW complex with fundamental group $\pi$. If $\Z\otimes_{\Z\pi}\Gamma(\pi_2(K))$ is torsion-free, then $\pi$ has the \cref{tertiary-property}.
\end{cor}

For example, computations showing that $\Z\otimes_{\Z\pi}\Gamma(\pi_2(K))$ is torsion-free have been done for finite 2-generator abelian groups~\cite{KPR}, dihedral groups~\cite{KNR}, and groups with 4-periodic cohomology~\cite{Ham-kreck-finite}.  We will make use of some of these in \cref{chapter:examples}.

We will prove the corollary via a couple of lemmas.

Let $\pi$ be a finite group and let $M$ be a left $\Z\pi$-module.
Consider the map of abelian groups
\begin{align*}
  F_M\colon M^*\otimes_\Z M^* &\to \Hom_{\Z\pi}(M,\Hom_{\Z\pi}(M,\Z\pi\otimes_\Z\Z\pi)) \\
  f\otimes h &\mapsto (m\mapsto(n\mapsto f(n)\otimes \overline{h(m)})).
  \end{align*}
  Here we view $\Z\pi\otimes \Z\pi$ as a left $\Z\pi$-module via left multiplication in the first factor. The group $\Hom_{\Z\pi}(M,\Z\pi\otimes_\Z\Z\pi)$ is a right $\Z\pi$-module via right multiplication in the second factor of $\Z\pi\otimes\Z\pi$ and as usual we view it as a left $\Z\pi$-module using the involution.

\begin{lemma}\label{lem:F_M-injective}
    Let $\pi$ be a finite group and let $M$ be a finitely generated left $\Z\pi$-module. Then the map $F_M$ is injective.
\end{lemma}
\begin{proof}
    Consider the abelian group homomorphism $\ev_1\colon M^*\to \Hom_\Z(M,\Z)$ given by evaluation at the neutral group element. This map has an inverse given by sending $f$ to the map $(x\mapsto \sum_{g\in\pi}f(g^{-1}x)g)$. Consider the commutative diagram
    \[\xymatrix{
    M^*\otimes M^*\ar[d]^{\ev_1\otimes \ev_1}_{\cong}\ar[r]^-{F_M}& \Hom_{\Z\pi}(M,\Hom_{\Z\pi}(M,\Z\pi\otimes_\Z\Z\pi))\ar[d]^{\ev_{1\otimes 1}}\\
    \Hom_\Z(M,\Z)\otimes \Hom_\Z(M,\Z)\ar[r]&\Hom_\Z(M,\Hom_\Z(M,\Z)).
    }
    \]
    Since $\Hom_\Z(M,\Z)$ is a free $\Z$-module, the lower horizontal map is an isomorphism. Hence $F_M$ is injective.
\end{proof}

Recall that we have an inclusion map $\omega \colon \Gamma(M^*)\to  M^*\otimes_\Z M^*$, and that
$\varpi \colon \Z\otimes_{\Z\pi}\Gamma(M^*)\to M^*\otimes_{\Z\pi} M^*$ denotes the map induced after tensoring with $\Z$.
In the codomain, as per the conventions stipulated at the start of \cref{sec:sesq-forms-tate-cohomology}, the first $M^*$ in the tensor product is a right $\Z\pi$-module and the second $M^*$ is a left $\Z\pi$-module.
We obtain the following commutative diagram, where the right-hand vertical map is induced by the $(\Z\pi,\Z\pi)$-bimodule homomorphism $\Z\pi\otimes_{\Z} \Z\pi\to\Z\pi$, $x\otimes y\mapsto xy$, where as above the left $\Z\pi$ action is via left multiplication in the first factor.
\[\xymatrix{
\Gamma(M^*)\ar@{^(->}[r]^-{\omega}\ar[d]& M^* \otimes_\Z M^* \ar@{^(->}[r]^-{F_M}\ar[d]&\Hom_{\Z\pi}(M,\Hom_{\Z\pi}(M,\Z\pi\otimes_\Z\Z\pi))\ar[d]\\
\Z\otimes_{\Z\pi}\Gamma(M^*)\ar[r]^-{\varpi}& M^*\otimes_{\Z\pi} M^* \ar[r]^-{\Phi_M}&\Sesq(M)
}\]
Let $N_\pi\in \Z\pi$ be the norm element $N_\pi := \sum_{g\in\pi}g$. It is easy to see that the map
\[\tr_M\colon \Z\otimes_{\Z\pi}M\to M,\quad 1\otimes m\mapsto N_\pi \cdot m\]
is well-defined. We obtain a map
\[\tr_{\Sesq}\colon \Sesq(M)\to \Hom_{\Z\pi}(M,\Hom_{\Z\pi}(M,\Z\pi\otimes_\Z\Z\pi))\] by applying $\tr\colon \Z\pi\cong \Z\pi\otimes_{\Z\pi}\Z\pi\to \Z\pi\otimes_\Z \Z\pi$. Then we also have that
\[F_M\circ \omega \circ \tr_{\Gamma(M^*)}=\tr_{\Sesq} \circ\varpi \circ  \Phi_M.\]

\begin{lemma}\label{lemma:torsion-gamma-kernel}
	Let $K$ be a connected, finite 2-complex with finite fundamental group $\pi$. Then	the kernel of the map
	\[\Phi_{H^2(K;\Z\pi)}\circ\varpi \colon\Z\otimes_{\Z\pi}\Gamma(\pi_2(K))\to \Sesq(H^2(K;\Z\pi))\] is the torsion in $\Z\otimes_{\Z\pi}\Gamma(\pi_2(K))$.
\end{lemma}

\begin{proof}
We will apply the discussion above with  $M = H^2(K;\Z\pi) \cong \pi_2(K)^*$, so that $M^* = \pi_2(K)$ by \cref{lem:pi2-H2-dual}.
Since $\Sesq(H^2(K;\Z\pi))$ is torsion free, all torsion elements have to lie in the kernel.
	Now suppose $a\in \Z\otimes_{\Z\pi}\Gamma(\pi_2(K))$ lies in the kernel of $\Phi_{H^2(K;\Z\pi)}\circ \varpi$.
	Then
	\begin{align*}0&=\tr_{\Sesq}\circ\Phi_{H^2(K;\Z\pi)}\circ \varpi(a)\\
	&=F_{H^2(K;\Z\pi)}\circ \omega\circ \tr_{\Gamma(\pi_2(K))}(a).
	\end{align*}
	Since $F_{H^2(K;\Z\pi)}$ and $\omega$ are injective, this implies that $N_\pi a=0$ and thus $0=|\pi| a \in \Z\otimes_{\Z\pi}\Gamma(\pi_2K)$. In particular, $a$ is a torsion element.
\end{proof}

\cref{cor:ter-finite-group} now follows immediately from \cref{lemma:torsion-gamma-kernel} and \cref{thm:ter-alg-condition}.


\chapter{Inheritance results}
\label{sec:inheritance}

In this chapter we give two inheritance results for each of the secondary and tertiary properties.
The last two show that the properties are inherited from odd index subgroups, so it suffices to prove the properties on these.
The inheritance results from this chapter will be be used in the next chapter to prove the corresponding properties for families of groups.

\section{Secondary property inheritance}
We will actually prove an inheritance result for \cref{cond} instead of the \cref{secondary-property} itself.

As above, let $K$ be a finite, connected $2$-complex with $\pi_1(K) =\pi$.  Let $(C_*,d_*)$ be the cellular chain complex $C_*(K;\Z\pi)$, let $(C^*,d^*)$ be the dual complex, and let $H := H^2(K;\Z\pi)$.

\begin{defi}
If $\phi\colon \pi\to G$ is a group homomorphism and $K_G$ a finite connected $2$-complex with $\pi_1(K_G)=G$, we obtain a homomorphism
\[\phi_*\colon \widehat{H}^0(\Sesq(H))\to \widehat{H}^0(\Sesq(H^2(K_G;\Z G)))\]
as follows. First construct a homomorphism \[\begin{array}{rcl} \widehat{H}^0(\Sesq(H)) &\to& \widehat{H}^0(\Sesq(\Z G\otimes C^2/\im(\id_{\Z G}\otimes d^2)))\\ \lambda &\mapsto & [a\otimes x]\otimes [a'\otimes x']\mapsto a\phi(\lambda(x,x'))\overline{a'}.\end{array}\]
Then compose this with the map on $\widehat{H}^0(\Sesq(-))$ induced by $H^2(K_G;\Z G)\to H^2(K;\Z G)\cong \Z G\otimes C^2/\im(\id_{\Z G}\otimes d^2)$.
\end{defi}

Recall that for a finitely presented group $\pi$ we defined a map $A \colon H_3(\pi;\Z/2) \to \widehat{H}^0(\Sesq(H^2(K_{\pi};\Z\pi)))$ in \cref{defn:map-A}.

\begin{lemma}\label{lem:commuting-diagram-A-phi}
The following square commutes
\[\xymatrix{
H_3(\pi;\Z/2)\ar[d]^A\ar[r]^-{\phi_*}&H_3(G;\Z/2)\ar[d]^A\\
\widehat{H}^0(\Sesq(H))\ar[r]^-{\phi_*}&\widehat{H}^0(\Sesq(H^2(K_G;\Z G)))}\]
\end{lemma}

\begin{proof}
By tracing the definition of the map $A \colon H_3(\pi;\Z/2) \to \widehat{H}^0(\Sesq(H^2(K;\Z\pi)))$, via the maps $\iota$, $I$, $\Theta$ and $\Phi$, it is not too hard to see, for $c \in C_3(B\pi;\Z\pi)$, $m,m' \in C^2 = C^2(K;\Z\pi)$ representatives of cohomology classes in $H^2(K;\Z\pi)$, that
\[A([1\otimes c]) = \big(([m],[m']) \mapsto \overline{m(\partial_3(c))}m'(\partial_3(c)) \in \Z\pi \big).\]
Here $\partial_3 \colon C_3(B\pi;\Z\pi) \to C_2(B\pi;\Z\pi)$ is the boundary map.
 The map \[\phi_* \colon C_*(B\pi;\Z\pi) \to \Z G \otimes_{\Z\pi} C_*(B\pi;\Z\pi) = C_*(B\pi;\Z G) \to C_*(BG;\Z G)\] is given by reducing coefficients followed by the chain map induced by the map of spaces $B\pi \to BG$.
In the following formulae, we also use the induced map
\[\phi^* \colon C^2(BG; \Z G) \to C^2(B\pi;\Z G) = \Hom_{\Z\pi}(C_2(B\pi;\Z\pi),\Z G) \cong \Z G\otimes_{\Z\pi} C^2(B\pi;\Z\pi).\]
Then we compute, with $n,n' \in C^2(K_G;\Z G) = C^2(BG;\Z G)$, as follows.
\begin{align*}
A(\phi_*([1\otimes c])) &= A([1 \otimes \phi_*(c)]) \\
&= \big(([n],[n']) \mapsto \overline{n(\partial_3(\phi_*(c)))}n'(\partial_3(\phi_*(c))) \big) \\
&= \big(([n],[n']) \mapsto \overline{n(\phi_*(\partial_3(c)))}n'(\phi_*(\partial_3(c))) \big) \\
&= \big(([n],[n']) \mapsto \big(\overline{\phi^*(n)(\partial_3(c))} \phi^*(n')(\partial_3(c))\big) \big) \\
&= \phi_* \big(([m],[m']) \mapsto \overline{m(\partial_3(c))} m'(\partial_3(c))\big) \\
&= \phi_* (A([1 \otimes c ])).\qedhere
\end{align*}
\end{proof}

Here is our first secondary property inheritance result.

\begin{thm}
	\label{thm:semidirectprod}
Let $G,G'$ be finitely presented groups satisfying \cref{cond}.
\begin{enumerate}
	\item Let $p\colon \pi\to G'$ be a homomorphism, and let $x\in H_3(\pi;\Z/2)$ be such that $p_*(x)$ does not lie in the image of $Sq_2\circ \red_2$. Then $A(x)\neq 0$.
	\item Let $i\colon G\to \pi$ be a homomorphism and let $x\in H_3(G;\Z/2)$ be in the image of $Sq_2\circ\red_2$. Then $A(i_*x)=0$.
\end{enumerate}
Moreover, if $\pi\cong N \rtimes G$ for some normal subgroup $N$ of $\pi$ and $x\in H_3(G;\Z/2)$, then $i_*(x)\in H_3(\pi;\Z/2)$ is in the kernel of $A$ if and only if it is in the image of $\Sq_2\circ\red_2$. Here $i\colon G\to\pi$ denotes the inclusion of $G$ into $\pi$ that sends $g \mapsto (1_{N},g)$.
\end{thm}

\begin{proof}
Let $K,K_{G'}$ and $K_{G}$ be finite $2$-complexes with fundamental groups $\pi, G'$ and $G$ respectively. Consider the following diagram
	\[\xymatrix @C-0.3cm{
		H_5(G;\Z)\ar[d]^{\Sq_2\circ \red_2}\ar[r]^-{i_*}&H_5(\pi;\Z)\ar[d]^{\Sq_2\circ \red_2}\ar[r]^-{p_*}&H_5(G';\Z)\ar[d]^{\Sq_2\circ \red_2}\\
		H_3(G;\Z/2)\ar[d]^A\ar[r]^-{i_*}&H_3(\pi;\Z/2)\ar[d]^A\ar[r]^-{p_*}&H_3(G';\Z/2)\ar[d]^A\\
		\widehat{H}^0(\Sesq(H^2(K_G;\Z G)))\ar[r]^-{i_*}&\widehat{H}^0(\Sesq(H^2(K;\Z \pi)))\ar[r]^-{p_*}&\widehat{H}^0(\Sesq(H^2(K_{G'};\Z G')))}\]
The bottom two squares commute by \cref{lem:commuting-diagram-A-phi}, and the top two squares commute by naturality of $\Sq_2$ and $\red_2$.  The left and right columns are exact by assumption, but the middle column need not be.
The enumerated items in the theorem now follow from \cref{thm:mainsec} and a diagram chase.

To see the last part, note that the homomorphism $i \colon G \to \pi = N \rtimes G$  is a splitting for the quotient map $Q \colon \pi \to \pi/N \cong G$, that is $Q \circ i = \Id_G$.
The map $Q$ induces a diagram similar to the diagram above.
	\[\xymatrix{
		H_5(G;\Z)\ar[d]^{\Sq_2\circ \red_2}\ar@/^/[r]^{i_*} & H_5(\pi;\Z)\ar[d]^{\Sq_2\circ \red_2} \ar @/^/[l]^{Q_*} \\
		H_3(G;\Z/2)\ar[d]^A \ar @/^/[r]^{i_*} & H_3(\pi;\Z/2)\ar[d]^A \ar @/^/[l]^{Q_*}\\
		\widehat{H}^0(\Sesq(H^2(K_G;\Z G))) \ar @/^/[r]^{i_*} & \widehat{H}^0(\Sesq(H^2(K;\Z \pi))) \ar @/^/[l]^{Q_*} }\]
The squares commute with the horizontal arrows going in both directions.  The left hand column is exact.  Also note that $Q_* \circ i_* = \Id$, but be warned that we say nothing about $i_* \circ Q_*$.  Another diagram chase now shows that $i_*(x)\in H_3(\pi;\Z/2)$ is in the kernel of $A$ if and only if it is in the image of $\Sq_2\circ\red_2$, as required.
\end{proof}

\begin{cor}
	\label{cor:secabelian}
\cref{cond} holds for all finitely generated abelian groups if and only if it holds for all abelian groups with at most three generators.
\end{cor}

\begin{proof}
Let $E$ be an abelian group. For every $y\in H_3(E;\Z/2)$, there exists a decomposition $E\cong E'\oplus B$, such that $E'$ has at most $3$ generators and $y$ is in the image of $i_*$, where $i\colon E'\to E$ is the inclusion. Hence the sequence
\[H_5(E;\Z)\xrightarrow{\Sq_2\circ\red_2} H_3(E;\Z/2) \xrightarrow{A} \widehat{H}^0(\Sesq(H^2(K_E;\Z E)))\]
is ``exact at $y$,'' that is $y$ lies in the kernel of $A$ if and only if it lies in the image of $\Sq_2\circ\red_2$,  by \cref{thm:semidirectprod}. The corollary now follows from \cref{thm:mainsec}.
\end{proof}

\section{Tertiary property inheritance}

Before we state our tertiary inheritance theorem, we have a couple of lemmas. The first is a well known fact from group theory.

\begin{lemma}\label{lem:kernel-finitely-gen}
Let $p \colon \pi \to G$ be a surjective homomorphism between finitely generated groups $\pi$ and $G$. If $G$ is finitely presented, then $\ker(p)$ is finitely normally generated.
\end{lemma}

\begin{proof}
Let $\pi\cong \langle X=\{x_i\}_{i=1}^n\mid R\rangle$, where $R$ is some set of relations. Let $G\cong \langle X'=\{x'_i\}_{i=1}^{n'}\mid R'=\{w'_j(X')\}_{j=1}^{m'}\rangle$, where $w'_i$ is a word in $n'$ variables and $w'_i(X')$ denotes the word given by inserting $x'_i$ for the $i$th variable. Let $y_i\in \pi$ be a preimage of $x'_i$ for all $i=1,\ldots,n'$, $Y=\{y_i\}_{i=1}^{n'}$ and let $y_i=v_i(X)$ for some words $v_i$. Then
\[\pi\cong \langle X,Y\mid R,\{y_i^{-1}v_i(X)\}_{i=1}^{n'}\rangle.\]
Let $p(x_i)=v_i'(X')$ for some words $v_i'$. Then
\[G \cong \langle X',p(X)\mid R', \{p(x_i)^{-1}v_i'(X')\}_{i=1}^n\rangle.\]
Since $p(Y)=X'$ by definition and the relations of $\pi$ hold in $G$ we also have
\[G\cong \langle X,Y\mid R,\{y_i^{-1}v_i(X)\}_{i=1}^{n'}, \{w'_j(Y)\}_{j=1}^{m'}, \{x_i^{-1}v_i'(Y)\}_{i=1}^{n}\rangle.\]
To obtain this presentation for $G$, substitute $p(Y) =X'$ and then remove all the instances of $p$.
In particular, the kernel of $p$ is normally generated by the elements $\{w'_j(X')\}_{j=1}^{m'}$ and $\{x_i^{-1}v_i'(Y)\}_{i=1}^{n}$.
\end{proof}

\begin{lemma}
\label{lem:taucommute}
Let $N$ be a codimension zero submanifold $($possibly with boundary$)$ of a spin $4$-manifold $M$ and let $i\colon N\to M$ be the inclusion. Suppose $\alpha\in \pi_2(N)$ has $\mu(\alpha)=0$ and $i_*(\alpha)\in \pi_2(M)$ is $\RPT$-characteristic. Then $\tau(\alpha)$ and $\tau(i_*\alpha)$ are well-defined and equal.
\end{lemma}

\begin{proof}
From $\mu(i_*\alpha)=i_*(\mu(\alpha))=0$, it follows that $\tau(i_*\alpha)\in \Z/2$ is well-defined.
The number of transverse intersection points between $\alpha$ and a generically immersed $\RPT$ $\beta$ is the same as those between $i_*\alpha$ and $i_*\beta$. Therefore $\alpha$ is $\RPT$-characteristic and $\tau(\alpha)\in\Z/2$ is well-defined.

Pair up the self-intersections of $\alpha$ with Whitney discs in such a way that the intersection of these discs with $\alpha$ computes $\tau(\alpha)$. Then the same Whitney discs can be used to compute $\tau(i_*\alpha)$ inside the bigger manifold $M$. Hence $\tau(\alpha)=\tau(i_*\alpha)$.
\end{proof}

Here is our first tertiary property inheritance result.

\begin{thm}
	\label{thm:tersurj}
	Let $\pi$ be a finitely presented group and assume that there are surjections $p_i\colon \pi\to G_i$ for some finitely presented groups $G_1,\dots,G_n$, such that
	\[\prod_{i=1}^n(p_i)_*\colon H_2(\pi;\Z/2)/\im (d^2_{4,1}, d^3_{5,0}) \to \prod_{i=1}^n H_2(G_i;\Z/2)/\im (d^2_{4,1}, d^3_{5,0})\]
	is an injection.
		Furthermore, assume that all the $G_i$ have the \cref{tertiary-property}. Then $\pi$ has the \cref{tertiary-property}.	
\end{thm}

\begin{proof}
	Let $M$ be a spin manifold with $\pi_1(M)=\pi$ and $c_*([M])=0\in H_4(B\pi;\Z)$. Choose a model for $B\pi$ with finite $2$-skeleton. Let $f\colon M^{(3)}\to B\pi^{(2)}$ be a map that induces an isomorphism on fundamental groups, and such that the equivariant intersection form vanishes on the image of \[f^*\colon H^2(B\pi^{(2)};\Z\pi)\to H^2(M^{(3)};\Z\pi) \cong H^2(M;\Z\pi) \cong H_2(M;\Z\pi).\]
	In particular, we assume that the primary and secondary obstructions vanish for~$M$.

Recall that $\RC \subset \pi_2(M) \cong H_2(M;\Z\pi)$ denotes the subset of $\RPT$-characteristic elements on which $\mu$ vanishes. 	Consider the following diagram, in which $H^2(M;\Z\pi)_{\RC}$ denotes the subset $\RC \cap \ker \mu$ of $H^2(M;\Z\pi)$, $H^2(B\pi^{(2)};\Z\pi)_{\RC}$ denotes the preimage of $H^2(M;\Z\pi)_{\RC}$, and $H^2(B\pi^{(2)};\Z/2)_{\RC}$ denotes the image of $H^2(M;\Z\pi)_{\RC}$. By \cref{lem:char}, the maps land in the subsets claimed.
	\[\xymatrix{&H^2(B\pi^{(2)};\Z\pi)_{\RC} \ar[r]^-{f^*}\ar[d]& H^2(M;\Z\pi)_{\RC} \ar[r]^-{\tau}&\Z/2\\
		H^2(B\pi;\Z/2)\ar[r] \ar@{-->}[ru] &H^2(B\pi^{(2)};\Z/2)_{\RC}&&}\]
Make a choice of lift, as shown by the dashed arrow; that is, for each element of $H^2(B\pi^{(2)};\Z/2)_{\RC}$ choose an element in $H^2(B\pi^{(2)};\Z\pi)_{\RC}$ that maps to it.

	We want to show that \[\kappa(\tau_f|_{\ker (d^3_{5,0})^*})=\ter(M) \in H_2(B\pi;\Z/2)/(\im d^2_{4,1},\im d^3_{5,0})\]
where the isomorphism $\Hom(H^2(B\pi;\Z/2),\Z/2) \cong H_2(B\pi;\Z/2)$ induces an identification
\[\kappa \colon \Hom (\ker (d^3_{5,0})^*,\Z/2) \raiso H_2(B\pi;\Z/2)/(\im d^2_{4,1},\im d^3_{5,0}).\]
Note that this includes showing that $\tau_f|_{\ker (d^3_{5,0})^*}$ is a homomorphism.

We will also write $\kappa$ for the corresponding identification with $B\pi$ replaced by $BG$.  Observe that it suffices, by the injectivity hypothesis of \cref{thm:tersurj}, to show that for each surjective homomorphism $p \colon \pi \to G$, where $G$ has the Tertiary property, that $p_* \kappa(\tau_f|_{\ker (d^3_{5,0})^*}) = p_* \ter (M)$.
The remainder of the proof verifies this equality.

The element $p_*\ter(M)$ can be computed via the bordism class $[M \xrightarrow{c} B\pi \xrightarrow{p} BG] \in \Omega_4^{\spin}(BG)$.
We want to identify a class in $H_2(BG;\Z/2)$, from the $E_2$-page of the Atiyah-Hirzebruch spectral sequence for $\Omega_4^{\spin}(BG)$, with a $\tau$ invariant, using that $G$ has the Tertiary property.  For this, we need a manifold with fundamental group~$G$.

Perform surgery on $M$ along a normal generating set of curves $\coprod S^1 \subset M$ for $\ker (p \colon \pi \to G)$; this is a finite set by \cref{lem:kernel-finitely-gen}.  This removes, for each surgery, a copy of $D^2 \times S^2$.  Note that we will allow ourselves in the future to modify the representative circles normally generating $\ker p$ within their homotopy classes.  Let
\[M' := M \sm \Big(\coprod S^1 \times D^3\Big) \cup_{\coprod S^1 \times S^2} \coprod D^2 \times S^2.\]
The gluing map in the surgery arises from the framing of the normal bundle of  each $S^1$ i.e.\ the identification of a regular neighbourhood of $S^1$ with $S^1 \times D^3$.  We use the unique (up to homotopy) identification for which  the spin structure of $M$ extends over $M'$.

Note that we can arrange for the surgery data $\coprod S^1 \times D^3$ to be contained in $M^{(3)}$. Choose a model for $BG$ with finite $2$-skeleton. We obtain a commutative diagram
\[\xymatrix{M^{(3)} \ar[dr]_{p \circ f} & M^{(3)}\sm\Big(\coprod S^1 \times D^3\Big) \ar[d]^{p\circ f=f'}\ar[r]^-{i_{M'}}\ar[l]_-{i_M} & (M')^{(3)} \ar[dl]^{f'} \\ &BG^{(2)}. & }\]

We will view $M^{(3)}\sm\Big(\coprod S^1 \times D^3\Big)$ as a manifold with boundary, by only removing the interiors of the $S^1\times D^3$. In the following diagram, the boundary will be denoted by $\partial$.

The proof will be based on the diagram in \cref{big-diagram}, which we will explain in detail below.

\begin{figure}
\adjustbox{center, scale=0.9}{%
\begin{tikzcd}[column sep=scriptsize]
  	H^2(BG^{(2)};\Z/2)& H^2(BG;\Z/2)\arrow[l,"i^*"]&&\\
	H^2(BG^{(2)};\Z G)\arrow[u,twoheadrightarrow,"\red_2"] \arrow[rr,"p^*"] \arrow[dd,"(f')^*"] \arrow[ddr, "(f')^*"]&&H^2(B\pi^{(2)};\Z G)\arrow[d, "f^*"]&H^2(B\pi^{(2)};\Z\pi)\arrow[l, "\red"]\arrow[dd, "f^*"]\\
	&&H^2(M^{(3)};\Z G)&\\
	H^2((M')^{(3)};\Z G)&H^2(M^{(3)}\setminus\coprod S^1\times D^3,\partial;\Z G)\arrow[l, "\ex"]\arrow[ur, "\ex"]&H^2(M^{(3)}\setminus\coprod S^1\times D^3,\partial;\Z \pi)\arrow[l, "\red"]\arrow[r, "\ex"]&H^2(M^{(3)};\Z\pi)\arrow[ul, "\red"]\\
	H^2(M';\Z G)\arrow[u, "\cong","i^*"']\arrow[d, "PD"]&H^2(M\setminus\coprod S^1\times D^3,\partial;\Z G)\arrow[l, "\ex"]\arrow[u, "\cong","i^*"']\arrow[d, "PD"]&H^2(M\setminus\coprod S^1\times D^3,\partial;\Z \pi)\arrow[l, "\red"]\arrow[r, "\ex"]\arrow[u, "\cong","i^*"'] \arrow[d, "PD"] & H^2(M;\Z\pi)\arrow[u, "\cong","i^*"'] \arrow[d, "PD"]\\
	H_2(M';\Z G)\arrow[d, "h^{-1}"]&H_2(M\setminus\coprod S^1\times D^3;\Z G)\arrow[l, "i_*"]&H_2(M\setminus\coprod S^1\times D^3;\Z \pi)\arrow[l, "\red"]\arrow[r, "i_*"]\arrow[d, "h^{-1}"]&H_2(M;\Z\pi)\arrow[d, "h^{-1}"]\\
	\pi_2(M')\arrow[rrd, dashrightarrow, "\tau"']&&\pi_2(M\setminus\coprod S^1\times D^3)\arrow[ll, "i_*"']\arrow[d,dashrightarrow,"\tau"]\arrow[r, "i_*"]&\pi_2(M)\arrow[dl,dashrightarrow,"\tau"]\\
	&&\Z/2&
\end{tikzcd}}
\caption{The proof of tertiary inheritance.}\label{big-diagram}
\end{figure}
The strategy of the proof is as follows. The passage from $H^2(BG;\Z/2)$ (including a choice of lift for the top left $\red_2$ map), that goes down the left hand side of the diagram, computes $\tau_{f'}$.
Note that
\[\tau_{f'}|_{\ker (d^3_{5,0}(G))^*} = \kappa^{-1}(\ter(M')) = \kappa^{-1}(p_*(\ter(M))).\]
The passage along the right hand side of the diagram computes $\tau_f \circ p^*$.  If we show that $\tau_{f'}$ and $\tau_f \circ p^*$ are equal, this will imply that $\tau_f$ is a group homomorphism on the image of $p^*$. The injectivity statement in \cref{thm:tersurj} dualises to a surjectivity statement $\prod_{i=1}^n p_i^*\colon \prod_{i=1}^n \ker (d_{5,0}^3(G_i))^* \twoheadrightarrow \ker (d_{5,0}^3(\pi))^*$  and hence we will obtain that $\tau_f|_{\ker(d^3_{5,0})^*(\pi)}$ is a homomorphism. Here $d^3_{5,0}(\pi)$ and $d^3_{5,0}(G_i)$ denote the differentials of the Atiyah-Hirzebruch spectral sequences for $B\pi$ and $BG_i$ respectively. By naturality of the spectral sequences, $p^* \colon H^2(G;\Z/2) \to H^2(\pi;\Z/2)$ restricts to a homomorphism $p^*\colon \ker (d^3_{5,0}(G))^* \to \ker (d^3_{5,0}(\pi))^*$.

Now we know that $\tau_f\circ p^*|_{\ker (d^3_{5,0}(G))^*}$ is a homomorphism, and we can apply $\kappa$ to the equality
$\kappa^{-1}(p_*(\ter(M))) = \tau_{f'}|_{\ker (d^3_{5,0}(G))^*}=\tau_f\circ p^*|_{\ker (d^3_{5,0}(G))^*}$ to obtain
\begin{align*}
p_*(\ter(M)) &= \kappa(\tau_{f'}|_{\ker (d^3_{5,0}(G))^*}) = \kappa((\tau_f \circ p^*)|_{\ker (d^3_{5,0}(G))^*}) \\ &= \kappa\big((p^*)^*(\tau_f|_{\ker (d^3_{5,0}(\pi))^*})\big) = p_*(\kappa(\tau_f|_{\ker (d^3_{5,0}(\pi))^*})). \end{align*}
This equality for all the surjective group homomorphisms $p_i \colon \pi \to G_i$, together with the injectivity assumption in \cref{thm:tersurj}, will then imply that $\ter(M) = \kappa(\tau_f|_{\ker (d^3_{5,0})^*})$ as desired.
Once we have explained the diagram and shown that it commutes, it will follow easily that $\tau_{f'}$ and $\tau_f \circ p^*$ are equal.

Next we explain the maps in the diagram in \cref{big-diagram}.
\begin{itemize}
\item Arrows labelled with $f^*,p^*,(f')^*$ are the maps induced by $f\colon M^{(3)} \to B\pi^{(2)}$, $p \colon \pi \to G$ and $f' \colon (M')^{(2)} \to BG^{(2)}$ respectively.
\item Arrows labelled with $\red$ are reduction of the coefficients.
\item Arrows labelled with $\ex$ are given by the inverse of excision
\[H^2(M,\coprod S^1\times D^3)\xrightarrow{\cong}H^2(M\setminus \coprod S^1\times D^3,\partial)\]
composed with the map of from the long exact sequence of a pairs
\[H^2(M,\coprod S^1\times D^3)\to H^2(M)\]
with the stated coefficients or by the analogous maps for $M^{(3)}$ or $M'$ instead of $M$. The fact that $M\setminus \coprod S^1\times D^3 \cong M'\setminus \coprod S^1\times D^3$ is also used in the definition of the left-hand $\ex$ maps.
\item Arrows labelled with $i^*$ or $i_*$ are maps induced by inclusions of the $3$-skeleta.
\item Arrows labelled with $PD$ are the isomorphisms from Poincar\'e duality.
\item Arrows labelled with $h^{-1}$ are the inverses of the Hurewicz isomorphism.
\item Arrows labelled with $\tau$ are the $\tau$ invariant. These arrows are dashed since the $\tau$ invariant is only defined on a subset.
\end{itemize}

\noindent More precisely, the map \[(f')^*\colon H^2(BG^{(2)};\Z G)\to H^2(M^{(3)}\sm\coprod S^1 \times D^3,\partial;\Z G)\]
is given as follows.

Consider the cohomology long exact sequence of the pair $(M'^{(3)}, \coprod D^2\times S^2)$:
\[0\to H^2(M'^{(3)}, \coprod D^2\times S^2;\Z G)\to H^2(M'^{(3)};\Z G)\to H^2(\coprod D^2\times S^2;\Z G)\]
Since $f'$ restricted to the attached copies of $D^2\times S^2$ is null homotopic, the map $f'\colon H^2(BG^{(2)};\Z G)\to H^2(M'^{(3)};\Z G)$ lands in the image of $H^2(M'^{(3)}, \coprod D^2\times S^2;\Z G)$.
We can compose this $(f')^*$ further with the map induced by the map of pairs
\[H^2(M'^{(3)}, \coprod D^2\times S^2;\Z G)\to H^2(M^{(3)}\setminus \coprod S^1\times D^3,\partial;\Z G)\]
(by excision this is an isomorphism) to obtain the diagonal $(f')^*$ map in the big diagram
\[(f')^*\colon H^2(BG^{(2)};\Z G)\to H^2(M^{(3)}\sm\coprod S^1 \times D^3,\partial;\Z G).\]

All quadrilaterals in the diagram commute by naturality of the involved maps. Note that one of the quadrilaterals looks like a triangle at first glance.  The commutativity of the two triangles at the bottom (when they are defined) follows from \cref{lem:taucommute}.

Let $\alpha\in H^2(BG^{(2)};\Z G)$ be given and let $\wt{\alpha}\in H^2(B\pi^{(2)};\Z\pi)$ be a lift of $p^*\alpha$.
\begin{claim}
There exists $y\in H^2(M^{(3)}\setminus \coprod S^1\times D^3,\partial;\Z\pi)$ with $\red(y)=(f')^*\alpha$ and $\ex(y)=f^*\wt{\alpha}$.
\end{claim}

The claim can be seen as follows. Consider the following commutative diagram with exact rows coming from the long exact sequences of pairs.
\begin{equation}\label{diagram:tau-proof}
  \xymatrix{
H^1(\coprod S^1\times D^3;\Z\pi)\ar[r]\ar@{->>}[d]^{\red}&H^2(M^{(3)},\coprod S^1\times D^3;\Z\pi)\ar[r]^-{j}\ar[d]^{\red}&H^2(M^{(3)};\Z\pi)\ar[d]^{\red}\ar[r]&0\\
H^1(\coprod S^1\times D^3;\Z G)\ar[r]&H^2(M^{(3)},\coprod S^1\times D^3;\Z G)\ar[r]&H^2(M^{(3)};\Z G)\ar[r]&0\\
}\end{equation}

The left hand vertical map is surjective, as we now argue.
The circles in the left hand column represent elements of $\pi$ that normally generate $\ker(\pi \to G)$.  Let us denote these elements of $\pi$ by $g_1,\dots,g_n$. We have that $H^1(\coprod S^1\times D^3;\Z G) \cong \bigoplus^n \Z G$, with one summand per copy of $S^1 \times D^3$.  On the other hand, $H^1(\coprod S^1\times D^3;\Z \pi) \cong \bigoplus_{i=1}^n \Z\pi /(g_i-1)$.  The reduction map simply adds the relations $g_j-1$ to each summand, for $j \neq i$.  This shows that the left hand vertical map is surjective.

Using excision and commutativity of the big diagram from \cref{big-diagram}, it not too hard to see, by taking $a_1 =f^*(\wt{\alpha})$ and $a_2 = (f')^*(\alpha)$, that the claim follows if we can show the following: for every $a_1\in H^2(M^{(3)};\Z\pi)$ and for every  $a_2\in H^2(M^{(3)},\coprod S^1\times D^3;\Z G)$ with the same image in $H^2(M^{(3)};\Z G)$, there exists $y\in H^2(M^{(3)},\coprod S^1\times D^3;\Z \pi)$ with $j(y)=a_1$ and $\red(y)=a_2$.  The map $$j \colon H^2(M^{(3)}, \coprod S^1\times D^3;\Z \pi) \to H^2(M^{(3)};\Z\pi)$$ is identified with the map $$H^2(M^{(3)} \setminus \coprod S^1\times D^3,\partial ;\Z \pi) \to H^2(M^{(3)};\Z\pi)$$ using excision.
The existence of such an element $y$ follows from a diagram chase in the above diagram~\eqref{diagram:tau-proof}. Lift $a_1$ to $z  \in H^2(M^{(3)},\coprod S^1\times D^3;\Z \pi)$.  The element $z$ might not map to $a_2$ under $\red$. Take $\red(z)-a_2$, lift it to the top left corner and map it to $w \in H^2(M^{(3)},\coprod S^1\times D^3;\Z \pi)$.  Define $y=z-w$. It is now straightforward to see that $y$ has the desired properties.
This completes the proof of the claim.

It follows from the claim that for any two elements $\alpha,\alpha'\in H^2(BG^{(2)};\Z G)$, there are elements \[PD(y),PD(y')\in \pi_2(M\sm \coprod S^1\times D^3)\] that map to $PD((f')^*\alpha), PD((f')^*\alpha)$ and $PD(f^*\wt \alpha), PD(f^*\wt \alpha')$ under the respective inclusions. Since the equivariant intersection form vanishes on $PD\circ f^*$, we see that
\[0=\lambda(PD(f^*\wt \alpha),PD(f^*\wt \alpha'))=\lambda(PD(y),PD(y'))=p\big(\lambda(PD((f')^*\alpha), PD((f')^*\alpha))\big).\]
Hence the equivariant intersection form of $M'$ vanishes on $PD\circ (f')^*$. For every $x\in H^2(BG;\Z/2)$ and every lift $\wh x\in H^2(BG^{(2)};\Z G)$ of $i^*x$ we have
\[\tau_{f'}(x)=\tau(PD((f')^*\wh x)).\]
For every lift $\wt x\in H^2(B\pi^{(2)};\Z\pi)$ of $p^*\wh x$ and every $y\in H^2(M^{(3)}\setminus \coprod S^1\times D^3,\partial;\Z\pi)$ with $\red(y)=(f')^*\wh x$ and $\ex(y)=f^*\wt{x}$ we have
\[\tau_f(p^*x)=\tau(PD(f^*\wt x))=\tau(PD(y))=\tau(PD((f')^*\wh x))=\tau_{f'}(x).\]
This completes the proof of the assertion that $\tau_f\circ p^* =\tau_{f'}$ and therefore completes the proof of the inheritance theorem for the Tertiary property.
\end{proof}

\begin{cor}
	\label{cor:terabelian}
The \cref{tertiary-property} holds for all finitely generated abelian groups if and only if it holds for all abelian groups with at most two generators.
\end{cor}

\begin{proof}
Assume that the \cref{tertiary-property} holds for all abelian groups with at most two generators.

Let $\pi=\bigoplus_{i=1}^nC_i$ with $C_i$ cyclic. Let $p_i\colon \pi\to C_i$ be the projection homomorphism and let  $s_i\colon C_i\to \pi$ be the inclusion homomorphism. For $i\neq j$ let $G_{ij}:=C_i\oplus C_j$, let $p_{ij}\colon G_{ij}\to \pi$ be the projection homomorphism and let $s_{ij}\colon G_{ij}\to \pi$ be the inclusion homomorphism. By the  K\"{u}nneth theorem
\[H_2(\pi;\Z/2)\cong \bigoplus_{1\leq i\leq n}H_2(C_i;\Z/2)\oplus\bigoplus_{1\leq i<j\leq n}H_1(C_i;\Z/2)\otimes H_1(C_j;\Z/2).\] It follows, that \[p:=\bigoplus_{1\leq i<j\leq n}(p_{ij})_*\colon  H_2(\pi;\Z/2)\to\bigoplus_{1\leq i<j\leq n}H_2(G_{ij};\Z/2)\] is injective. Let
\[s:=\bigoplus_{1\leq i<j\leq n}(s_{ij})_*\colon \bigoplus_{1\leq i<j\leq n}H_2(G_{ij};\Z/2)\to H_2(\pi;\Z/2).\]
Depending on $n$, the composition $s\circ p$ might not be the identity since any nontrivial  elements from $H_2(C_i;\Z/2)$ (i.e.\ whenever $C_i$ has even order) appear $(n-1)$ times in the terms $H_2(G_{ij};\Z/2)$ and hence get multiplied by $(n-1)$. For $n$ even this requires no change, but for $n$ odd this needs a remedy.  We give a unified treatment.
Define the $\Z/2$-module
\[\Lambda_n := \begin{cases}
\bigoplus_{1\leq i<j\leq n}H_2(G_{ij};\Z/2) \oplus \bigoplus_{1\leq i\leq n} H_2(C_i;\Z/2) & n \text{ odd}\\
\bigoplus_{1\leq i<j\leq n}H_2(G_{ij};\Z/2) \oplus \{0\} & n \text{ even.}
\end{cases}
\]
For $n$ odd, define
\[p'_n:=\bigoplus_{1\leq i\leq n}(p_i)_*\colon H_2(\pi;\Z/2)\to \bigoplus_{1\leq i\leq n} H_2(C_i;\Z/2)\]
and
\[s'_n:=\bigoplus_{1\leq i\leq n}(s_i)_*\colon \bigoplus_{1\leq i\leq n} H_2(C_i;\Z/2)\to H_2(\pi;\Z/2).\]
For $n$ even, take $p'_n \colon H_2(\pi;\Z/2) \to \{0\}$ and $s'_n \colon \{0\} \to H_2(\pi;\Z/2)$ to be the zero maps.
Then
\[(s \oplus s'_n)\circ(p \oplus p'_n) \colon H_2(\pi;\Z/2) \to \Lambda_n \to H_2(\pi;\Z/2)\]
is the identity map.
If follows by naturality of the Atiyah Hirzebruch spectral sequence that this composition still induces the identity map after modding out the $d^2_{4,1}$ and $d^3_{5,0}$ differentials.

In particular, the induced map
\begin{align*}
p\oplus p'_n\colon H_2(\pi;\Z/2)/\im (d^2_{4,1}, d^3_{5,0})\to \Lambda_n/\im(d^2_{4,1}, d^3_{5,0})
\end{align*}
is injective.
We can therefore apply the inheritance result \cref{thm:tersurj} above, with $G_i$, in the notation of that theorem, as all the abelian groups $G_{ij}$ or $C_i$ with either one or two generators involved in the sum $\Lambda_n$.
We conclude that $\pi$ has the \cref{tertiary-property} as asserted.
\end{proof}

\section{Reduction to odd index subgroups}\label{sec:reduction-odd-index}

This section will show that if either the \cref{secondary-property} or the \cref{tertiary-property} holds for a finite index subgroup $P \leq \pi$ where the index of $P$ in $\pi$ is odd, then the corresponding property also holds for $\pi$.

Let $P \leq \pi$ be a subgroup of a finitely presented group $\pi$.  Let $X$ be a CW complex with fundamental group $\pi$.  There is a covering space $p\colon \wh{X} \to X$ with fundamental group $P$.  Both spaces have the same universal cover $\wt{X}$, and we can express the cohomology of the universal cover as $H^*(X;\Z\pi) \cong H^*(\wh{X};\Z P)$, where the isomorphism is of $\Z P$-modules. Similarly, there are isomorphisms of the homology groups  $H_*(X;\Z\pi) \cong H_*(\wh{X};\Z P)$, where again the $\Z\pi$ homology is thought of as a $\Z P$ module by restriction.  We will make use of this simple observation several times in this section.

Now suppose that $P \leq \pi$ is a finite index subgroup. If $M$ is a right $\Z\pi$-module we can restrict the action to $P$ and consider the projection
\[M\otimes_{\Z P}C_*\xrightarrow{p}M\otimes_{\Z\pi}C_*\]
inducing
\[p_*\colon H_*(\widehat{X};M)\to H_*(X;M).\]
On the chain level we obtain a \emph{transfer map} in the other direction by
\[\tr\colon M\otimes_{\Z\pi}C_*\to M\otimes_{\Z P}C_*,\quad m\otimes c\mapsto \sum_{Pg\in P\backslash \pi}mg^{-1}\otimes gc.\]
This map induces a map
\[\tr_*\colon H_*(X;M)\to H_*(\widehat{X};M).\]
Similarly one constructs a map
\[\tr^*\colon H^*(\widehat X;M)\to H^*(X;M).\]
The transfer maps have the key property that $p_* \circ \tr_*$ and $\tr^*\circ p^*$ are equal to multiplication by the index $[\pi:P]$ of $P$ in $\pi$.

Let $K$ be a finite $2$-complex with $\pi_1(K)=\pi$. Then there is a finite cover $K_P$ of $K$ with fundamental group $P$.

Start with a map $f\colon H^2(K;\Z\pi)\to \Hom_{\Z\pi}(H^2(K;\Z\pi),\Z\pi)$, and define \[\tr(f)\in \Sesq(H^2(K_P;\Z\pi))\] to be the composition
\begin{align*} H^2(K_P;\Z\pi) & \xrightarrow{\tr^*}H^2(K;\Z\pi)\xrightarrow{f} \Hom_{\Z\pi}(H^2(K;\Z\pi),\Z\pi) \\ &\xrightarrow{(\tr^*)^*}\Hom_{\Z\pi}(H^2(K_P;\Z\pi),\Z\pi).\end{align*}
Post-compose with the isomorphism
\[\Hom_{\Z\pi}(H^2(K_P;\Z\pi),\Z\pi)\cong \Hom_{\Z P}(H^2(K_P;\Z\pi),\Z P),\]
given by post-composition with the map
\[\begin{array}{rcl}
  \ev_{\Z P}\colon \Z\pi &\to & \Z P \\
\sum_{g\in\pi}n_g g &\mapsto & \sum_{g\in P}n_g g,
\end{array}\] and pre-compose with the map $H^2(K_P;\Z P)\to H^2(K_P;\Z\pi)$ given by the inclusion $\Z P\to \Z\pi$, to restrict $\tr(f)$ to an element of $\Sesq(H^2(K_P;\Z P))$ i.e.\ a map $H^2(K_P;\Z P) \to \Hom_{\Z P}(H^2(K_P;\Z P),\Z P)$. This construction commutes with the transposition $T$ and so defines a transfer map
\[\tr_*\colon \wh{H}^0\big(\Sesq(H^2(K;\Z\pi))\big) \to \wh{H}^0\big(\Sesq(H^2(K_P;\Z P))\big).\]

\begin{lemma}
\label{lem:A-tr-commute}
The following diagram commutes:
\[\xymatrix{
H_3(\pi;\Z/2)\ar[d]^{\tr_*}\ar[r]^-A&\wh{H}^0(\Sesq(H^2(K;\Z\pi)))\ar[d]^{\tr_*}\\
H_3(P;\Z/2)\ar[r]^-A&\wh{H}^0(\Sesq(H^2(K_P;\Z P)))
}\]
\end{lemma}

\begin{proof}
By tracing the definition of the map $A \colon H_3(\pi;\Z/2) \to \widehat{H}^0(\Sesq(H^2(K;\Z\pi)))$, via the maps $\iota$, $I$, $\Theta$ and $\Phi$, it is not too hard to see, for $c \in C_3(B\pi;\Z\pi)$, $m,m' \in C^2 = C^2(K;\Z\pi)$ representatives of cohomology classes in $H^2(K;\Z\pi)$, that
\[A([1\otimes c]) = \big(([m],[m']) \mapsto \ol{m(\partial_3(c))}m'(\partial_3(c)) \in \Z\pi \big).\]
Here $\partial_3 \colon C_3(B\pi;\Z\pi) \to C_2(B\pi;\Z\pi)$ is the boundary map.

Similarly to the proof that the above definition of $A([1\otimes c])$ does not depend on the preimage $c\in C_3$ of $[1\otimes c]$, one checks that for $[c']\in H_3(P;\Z/2)$ and for any choice of a preimage $\wt{c'}\in \Z\pi\otimes_{\Z P}C_3(B P;\Z P)$, the $\wh{H}^0(\Sesq H^2(K_P;\Z P))$ class of the map
\[([m],[m']) \mapsto \ev_{\Z P}\big(\ol{(1\otimes m)(\partial_3(\wt{c'}))}(1\otimes m')(\partial_3(\wt{c'})) \in \Z\pi \big)\]
agrees with $A([c'])$.

For convenience we provide the details. For all $g,g'\in\pi$ and all $c\in C_3(BP;\Z P)$ we have the following equality in $\widehat{H}^0(\Sesq (H^2(K_P;\Z P)))$.
\begin{eqnarray*}
	&&\left( (m,m')\mapsto \ev_{\Z P}(\overline{(1\otimes m)((g\pm g')\otimes \partial_3(c))}(1\otimes m')((g\pm g')\otimes \partial_3(c)))  \right)\\
	&=&\left( (m,m')\mapsto \ev_{\Z P}(\overline{m(\partial_3(c))}((\overline{g\pm g'})(g\pm g'))m'(\partial_3(c)))  \right)\\
	&=&\left( (m,m')\mapsto \ev_{\Z P}(\overline{m(\partial_3(c))}(2\pm(g^{-1}g'+\overline{g^{-1}g'}))m'(\partial_3(c)))  \right)\\
	&=&\left( (m,m')\mapsto \ev_{\Z P}(\overline{m(\partial_3(c))}(g^{-1}g'+\overline{g^{-1}g'})m'(\partial_3(c))) \right)\\
	&=&0
\end{eqnarray*}

Since any two choices of preimage of $[c'] \in H_3(P;\Z/2)$ in $\Z\pi \otimes_{\Z P} C_3(BP;\Z P)$ differ by a sum of elements of the form $(g \pm g') \otimes c$, this shows that the element in $\widehat{H}^0(\Sesq (H^2(K_P;\Z P)))$ is independent of the choice of preimage.
Hence,
\begin{align*}
A(\tr_*[1\otimes c])&=\ev_{\Z P}\big(([m],[m']) \mapsto \ol{(1\otimes m)(\partial_3(\tr_* c))}(1\otimes m')(\partial_3(\tr_* c)) \big)\\
&=\ev_{\Z P}\big(([m],[m']) \mapsto \ol{\tr^*(1\otimes m)(\partial_3(c))}\tr^*(1\otimes m')(\partial_3(c))\big)\\
&=\tr_* A([1\otimes c]).\qedhere
\end{align*}
\end{proof}

\subsection{Secondary Property inheritance from odd index subgroups}

Like the Secondary inheritance from Section~\ref{sec:inheritance}, we will not actually show that the \cref{secondary-property} is inherited from finite odd index subgroups, but instead show inheritance for \cref{cond}.  A group can have the Secondary Property without necessarily satisfying \cref{cond}, although \cref{cond} implies the \cref{secondary-property} by \cref{thm:mainsec}.

\begin{theorem}
\label{lem:2sylow-sec}
Let $\pi$ be a finitely presented group and let $P$ be a finite odd index subgroup of~$\pi$. Then \cref{cond} holds for $\pi$ if it holds for~$P$.
\end{theorem}

\begin{proof}
Let $K$ be a finite $2$-complex with $\pi_1(K)=\pi$. Then there is a finite cover $K_P$ of $K$ with fundamental group $P$. In particular, $H^2(K_P;\Z P)\cong H^2(K;\Z\pi)$, by the observation made at the beginning of this section.

Consider the following diagram. Since $i_*\circ \tr_*$ is multiplication by the index $[\pi:P]$, which is odd, the middle vertical composition is the identity.
\[\xymatrix @C+0.2cm{
H_5(\pi;\Z)\ar[r]^-{Sq_2\circ \red_2}\ar[d]^{\tr_*}&H_3(\pi;\Z/2)\ar[d]^{\tr_*}\ar[r]^-A&\wh{H}^0(\Sesq(H^2(K;\Z\pi)))\ar[d]^{\tr_*}\\
H_5(P;\Z)\ar[r]^-{Sq_2\circ \red_2}\ar[d]^{i_*}&H_3(P;\Z/2)\ar[r]^-A\ar[d]^{i_*}&\wh{H}^0(\Sesq(H^2(K_P;\Z P)))\ar[d]^{i_*}\\
H_5(\pi;\Z)\ar[r]^-{Sq_2\circ \red_2}&H_3(\pi;\Z/2)\ar[r]^-A&\wh{H}^0(\Sesq(H^2(K;\Z\pi)))
}\]
Exactness of the top and bottom rows (which are equal) now follows from the exactness of the middle row by a diagram chase using the fact that $i_* \circ \tr_*=\id_{H_3(\pi;\Z/2)}$.
\end{proof}

\subsection{Tertiary property inheritance from odd index subgroups}

\begin{theorem}
\label{lem:2sylow-ter}
Let $\pi$ be a finitely presented group, and let $P$ be a finite index subgroup of odd index. Then the \cref{tertiary-property} holds for $\pi$ if it holds for $P$.
\end{theorem}

\noindent For the proof we first need a couple of lemmas.

\begin{lemma}
\label{lem:red-tr-commute}
Let $X$ be a finite CW complex with fundamental group $\pi$.
The composition
\[H^2(\wh{X};\Z P)\xrightarrow{\red}H^2(\wh{X};\Z/2)\xrightarrow{\tr^*}H^2(X;\Z/2)\]
agrees with the reduction of coefficients
\[H^2(X;\Z\pi)\xrightarrow{\red}H^2(X;\Z/2)\]
under the identification of $H^2(\wh{X};\Z P)$ and $H^2(X;\Z\pi)$ given at the beginning of the subsection.
\end{lemma}

\begin{proof}
Let $C_*$ denote the cellular $\Z\pi$-chain complex of $X$. On the chain level, the identification of $H^2(X;\Z\pi)$ with $H^2(\wh X;\Z P)$ is given by
\[\ev_{\Z P}\colon \Hom_{\Z\pi}(C_2,\Z\pi)\to \Hom_{\Z P}(C_2,\Z P).\]
Here $\ev_{\Z P} \colon \Z\pi \to \Z P$, sending $\sum_{g \in \pi} n_g g$ to $\sum_{g \in P} n_g g$, induces the map on $\Hom$ modules above by post-composition; we abuse notation and also denote the map on $\Hom$ modules by $\ev_{\Z P}$.

The reduction of coefficients is given by post-composing with the augmentations $\epsilon\colon \Z\pi\to \Z/2$ and $\epsilon\colon \Z P\to \Z/2$ respectively. The lemma now follows from the straightforward computation that the following diagram commutes.
\[\xymatrix{
\Hom_{\Z\pi}(C_2,\Z\pi)\ar[r]^-{\ev_{\Z P}}\ar[d]^\epsilon&\Hom_{\Z P}(C_2;\Z P)\ar[d]^\epsilon\\
\Hom_{\Z\pi}(C_2,\Z/2)&\Hom_{\Z P}(C_2,\Z/2)\ar[l]_-{\tr^*}
}\]
\end{proof}

\begin{lemma}
\label{lem:naturality-ter}
Let $i\colon P\to \pi$ be the inclusion. Let $M$ be a spin manifold with fundamental group $\pi$ and such that $\pri(M) = \msec(M) =0$. Then \[i_*\colon H_2(P;\Z/2)/\im (d_{5,0}^3, d_{4,1}^2) \to H_2(\pi;\Z/2)/\im (d_{5,0}^3, d_{4,1}^2)\] maps $\ter(\wh M)$ to $\ter(M)$, where $p \colon \wh{M} \to M$ is the cover associated to $P \leq \pi$.
\end{lemma}

\begin{proof}
Stabilisation of $M$ by a single $S^2 \times S^2$ corresponds to stabilisation of $\wh M$ with $[\pi:P]$ copies of $S^2\times S^2$. Since stabilisation does not change the $\ter$ invariant, and since stably the map $c\colon M\to B\pi$ factors through $B\pi^{(2)}$ up to homotopy (because $\pri(M) = \msec(M) =0$), we will assume that we are in the situation that $c\colon M\to B\pi$ maps to the $2$-skeleton $B \pi^{(2)} \subset B\pi$, and therefore also that $\wh{c} = c\circ p \colon \wh{M} \to M \to B\pi^{(2)}$ factors as $M \to BP^{(2)} \to B\pi^{(2)}$, through the 2-skeleton of $BP$. By \cref{lem:arf}, $\ter(M)$ is given as follows. Choose a point $e_i$ in each $2$-cell $c_i$ of $B\pi^{(2)}$, and let $F_i$ be the regular preimage of $e_i$ under $c$. The spin structure of $M$ induces a spin structure on $F_i$ and we have $\ter(M,c)=\big[\sum_{i} \Arf(F_i)c_i\big]$. A regular preimage of $e_i$ under $c\circ p \colon \wh{M} \to B\pi^{(2)}$ consists of $[\pi:P]$ copies of $F_i$, each with the same spin structure. Therefore, when we view $\wh M$ as an element of $\Omega_4^{\spin}(B\pi)$, we get
\[\ter(\wh M, c\circ p)=[\pi:P]\Big[\sum_{i} \Arf(F_i)c_i\Big]=[\pi:P]\ter(M,c)=\ter(M,c),\]
where the last equality uses that the index $[\pi:P]$ is odd. Since $\ter(\wh M, c\circ p)=i_*\ter(\wh M,\wh c)$, the lemma follows.
\end{proof}

\begin{proof}[Proof of \cref{lem:2sylow-ter}]
Let $K$ be a finite $2$-complex with fundamental group $\pi$, and let $\wh{K}$ be the covering space corresponding to $P$.

Every splitting $s$ of $H^2(M;\Z\pi)\to H^2(K;\Z\pi)$ also is a splitting $s_P$ of $H^2(\wh M;\Z P)\to H^2(\wh K;\Z P)$, because $H^2(\wh M;\Z P) = H^2(M;\Z\pi)$ and $H^2(\wh K;\Z P) = H^2(K;\Z\pi)$, with the latter considered as $\Z P$-modules.

 By \cref{lem:naturality-ter} it suffices to show that $i_*\tau_{s_P}=\tau_s$, since then $\tau_s=i_*\tau_{s_P}=i_*\ter(\wh M)=\ter(M)$. Here for the middle equality we used the assumption that $P$ has the \cref{tertiary-property}.

Recall that the map $\tau_{s_P}\colon H^2(B\pi;\Z/2)\to \Z/2$ is defined as follows: restrict $x\in H^2(BP;\Z/2)$ to an element of $H^2(\wh K;\Z/2)$, choose a preimage $x'\in H^2(\wh K;\Z P)$ and then apply $\tau\circ PD\circ s_P$. By \cref{lem:red-tr-commute} with $X=K$, $x'\in H^2(\wh K;\Z P)\cong H^2(K;\Z\pi)$ is a preimage of the restriction of $\tr^*(x) \in H^2(B\pi;\Z/2)$ to $H^2(K;\Z/2)$. Hence $\tau_s(\tr^*(x))=\tau(PD(s(x')))=\tau_{s_P}(x)$ and we obtain that $\tau_s\circ \tr^*=\tau_{s_P}$. Equivalently, we have $\tr_*\tau_s=\tau_{s_P}$, when we view $\tau_s$ and $\tau_{s_P}$ as elements of $H_2(\pi;\Z/2)\cong \Hom(H^2(\pi;\Z/2),\Z/2)$ and $H_2(P;\Z/2)\cong \Hom(H^2(P;\Z/2),\Z/2)$ respectively. Apply $i_*$ to obtain
\[i_*\tau_{s_P}=i_*\tr_*\tau_s=[\pi:P]\tau_s=\tau_s,\]
where the last equality follows from the assumption that the index of $P$ in $\pi$ is odd.
\end{proof}


\chapter{Examples}\label{chapter:examples}

In this chapter we prove that several families of groups possess the Secondary and Tertiary properties, thus classifying, at least modulo the limitations discussed in the introduction, spin $4$-manifolds with these fundamental groups up to stable diffeomorphism.
In particular we prove that both properties hold for all finitely generated abelian groups.

\section{The secondary and tertiary properties for three families of groups}
\label{sec:ex}

We begin by considering the following three families of groups: cohomologically 3-dimensional groups, right-angled Artin groups, and generalised quaternion groups.

\subsection{Cohomologically 3-dimensional groups}

A group $G$ is said to have cohomological dimension at most $n$ if $\Z$, thought of as a $\Z G$-module via the augmentation map $\Z G \to \Z$, admits a $\Z G$-module projective resolution of length $n$. A group $G$ is said to be of type $FP_n$ if $\Z$ has a projective resolution such that the first $n$ terms are finitely generated.
A group $G$ has type $FP_3$ and cohomological dimension at most three if and only if there exists a finite $2$-complex $K$ with $\pi_1(K)\cong G$ such that $\pi_2(K)$ is finitely generated projective as a $\Z G$-module.

\begin{theorem}\label{thm:conjs-cohom-dim-three}
The \cref{secondary-property} and the \cref{tertiary-property} hold for groups of type $FP_3$ that have cohomological dimension at most three.
\end{theorem}

Under the assumption that the equivariant intersection form is even, Hambleton-Hildum~\cite{Hamb-Hild} give an alternative classification to ours, using the stable quadratic $2$-type instead of the $\tau$ invariant.  By \cref{thm:conjs-cohom-dim-three}, the intersection form being even is equivalent to $\msec(M)=0$. As discussed in the introduction, Hambleton and Hildum use the whole intersection form. This means that their results also apply to the unstable homeomorphism problem, unlike ours.  On the other hand, the $\tau$-invariant on a summand of $\pi_2(M)$ gives an alternative formulation. In favourable circumstances the $\tau$-invariant requires fewer computations, because one considers only a subset of $\pi_2(M)$. In addition deciding whether two $\Z/2$ valued functions coincide could be easier than deciding whether two sesquilinear, hermitian forms over the group ring are stably isometric.

The proof of \cref{thm:conjs-cohom-dim-three} is split into two parts, one for each property.  Each property is dealt with in its own subsection below.

\subsubsection*{Three dimensional groups have the secondary property}\label{section:sec-3-dim-groups}

For the convenience of the reader, we recall that \cref{cond} requires the sequence
\[H_5(\pi;\Z)\xrightarrow{\Sq_2\circ \red_2}H_3(\pi;\Z/2)\xrightarrow{A}\widehat{H}^0(\Sesq(H))\]
be exact at $H_3(\pi;\Z/2)$.

For a cohomologically 3-dimensional group $\pi$ of type $FP_3$, choose a $2$-dimensional complex $K$ with $\pi_1(K)=\pi$ and such that $\pi_2 (K)$ is a finitely generated projective $\Z\pi$-module. Let $H:=H^2(K;\Z\pi)$.  Recall from \cref{lem:pi2-H2-dual} that $\pi_2(K) \cong H^*$.
Then the double dual $H^{**}$ is the dual of $\pi_2(K)$ and hence is a finitely generated projective $\Z\pi$-module. We can consider the natural evaluation map $e_H \colon H \to H^{**}$. Since $H^*\cong \pi_2(K)$ is finitely generated projective, the map $e_{H^*}\colon H^* \to H^{***}$ is an isomorphism. It is straightforward to verify that $e_H^*\circ e_{H^*} \colon H^* \to H^*$ is the identity on $H^*$.  Thus the map $e_H^* \colon H^{***} \to H^*$ is an isomorphism.

For a $\Z\pi$-module $M$ and a projective $\Z \pi$-module $P$, the canonical map $$\Hom_{\Z\pi}(M,\Z\pi)\otimes_{\Z\pi}P\to \Hom_{\Z\pi}(M,P)$$ is an isomorphism. Hence the map $\Hom_{\Z\pi}(H^{**},H^{*})\to \Hom_{\Z\pi}(H,H^*)$ given by precomposing with $e_H$ is an isomorphism, because $e_H^*$ is an isomorphism, and the given map can be identified with $e_H^* \otimes \Id_{H^*}$. This implies that
\begin{align*}
\Sesq(e_H) \colon \Sesq(H^{**}) \cong & \Hom_{\Z\pi}(H^{**},H^{***}) \xrightarrow{\Hom(-,e_H^*)} \\ & \Hom_{\Z\pi}(H^{**},H^*) \xrightarrow{\cong} \Hom_{\Z\pi}(H,H^*)\cong \Sesq(H)
\end{align*}
is also an isomorphism.

Consider the diagram
\[\xymatrix @R+0.5cm @C+0.5cm{
H^{***}\otimes_{\Z\pi} H^{***} \ar[d]^-{\Phi_{H^{**}}}\ar[r]_-{\cong}^-{e_H^*\otimes e_H^*}&H^*\otimes_{\Z\pi}H^*\ar[d]^{\Phi_{H}}\\
\Sesq(H^{**})\ar[r]_-{\cong}^-{\Sesq(e_H)} & \Sesq(H)}\]
Now the map $\Phi_{H^{**}}$ is an isomorphism by \cref{lem:injphi}, and hence $\Phi_H$ is also an isomorphism.
It follows that $\Psi_H = \widehat{H}^0(\Phi_H) \circ \Theta_H$ is an isomorphism as well, since $\Theta_H$ is an isomorphism by \cref{lem:Theta-Iso}. The map $A \colon H_3(\pi;\Z/2)\to \widehat{H}^0(\Sesq(H))$, which by definition is a composite $\Psi_H\circ (\Id \times I) \circ \iota$, where $I$ is an isomorphism by \cref{lem:pi2-H2-dual} and $\iota$ is injective by \cref{lem:inj-iota-proj}, is therefore injective.  Since $H_5(\pi;\Z)=0$, this implies that \cref{cond} holds, and therefore that $\pi$ has the \cref{secondary-property}.

\subsubsection*{Three dimensional groups have the tertiary property}

Since $\Phi_H$ is injective by the argument in \cref{section:sec-3-dim-groups}, and $\iota_{H^*}$ is injective  by \cref{lem:inj-iota-proj}, it follows from \cref{cor:conjter} that $\pi$ has the \cref{tertiary-property}.

\subsection{Right angled Artin Groups}

Let $(V,E)$ be a finite graph, with $V=\{v_i\}_{i\in I}$  the set of vertices and $E=\{e_j=\{v_j^1,v_j^2\}\}_{j\in J}$ the set of edges, where $e_j$ is an edge between $v_j^1, v_j^2\in V$. The corresponding right-angled Artin group (or RAAG for short) is $\langle V\mid \{[v_j^1,v_j^2]\}_{j\in J}\rangle$.
An \emph{$n$-clique} is a subset of $n$ vertices $v_1,\dots,v_n \in V$, such that for every pair of vertices $v_i,v_j$, with $i \neq j \in \{1,\dots,n\}$, we have $\{v_i,v_j\} \in E$.

\begin{theorem}
The \cref{secondary-property} and the \cref{tertiary-property} hold for right-angled Artin groups.
\end{theorem}

We will use the following well-known theorem.
\begin{thm}[{\cite[Corollary 3.2.2]{charney-davis}}]
\label{thm:homraags}
Let $\pi$ be the RAAG associated to a graph $(V,E)$. The dimension of $H_n(\pi;\Z/2)$ is the number of $n$-cliques in $(V,E)$.
\end{thm}

\subsubsection*{RAAGs have the secondary property}\label{sec:RAAG-sec}

For each $n$-clique in $(V,E)$ with vertices $v_k'$, $k=1,\ldots,n$, there is an inclusion $\Z^n=\langle x_k\mid [x_k,x_{k'}]\rangle\to \pi$ given by $x_k\mapsto v_k'$, and a projection $\pi\to \Z^n$ given by $v_k'\mapsto x_k$ and $v_i\mapsto 0$ for all $v_i\notin \{v_k'\}_{k=1,\ldots, n}$. In particular, $\pi$ is a semi-direct product $N\rtimes \Z^n$ of $\Z^n$ and a group $N$.

We apply this with $n=3$, for each $3$-clique.  By \cref{thm:homraags}, the rank of $H_3(\pi;\Z/2)$ is the same as the number of $3$-cliques.  We can find a basis for $H_3(\pi;\Z)$ such that each basis element is the image $i_*(x)$ for some $x \in H_3(\Z^3;\Z/2) =\Z/2$.  We know that $\Z^3$ has the Secondary property by \cref{thm:conjs-cohom-dim-three}.  Then the ``moreover'' part of \cref{thm:semidirectprod} says that $i_*(x)$ lies in the kernel of $A$ if and only if it lies in the image of $\Sq_2 \circ \red_2$.  Apply this argument at each $3$-clique, and therefore for each basis element of $H_3(\pi;\Z/2)$, to see that \cref{cond} holds for $\pi$. Therefore $\pi$ has the \cref{secondary-property}.

\subsubsection*{RAAGs have the tertiary property}

By \cref{thm:homraags} and the discussion at the beginning of the proof of the Secondary Property in \cref{sec:RAAG-sec}, it follows that there are surjections $p_j \colon \pi \to \Z^2$ such that $H_2(\pi;\Z/2)\to H_2(\oplus_{j \in J}\Z^2;\Z/2)$ is injective, where as above $J$ indexes the edges of the graph defining the RAAG $\pi$. Since all differentials in the spectral sequence for $\Omega_4^{\spin}(B\Z^2)$ are trivial, this map is still injective after dividing out the image of the differentials. Hence the \cref{tertiary-property} for $\pi$ follows from \cref{thm:tersurj} and the fact that we know the property holds for $\Z^2$ (since $\Z^2$ has cohomological dimension two).

\subsection{Generalised quaternion groups}\label{sec:gen-quaternion-groups}

By \cite[Theorem VI 9.3]{brown}, every $2$-group with periodic cohomology is either a cyclic group or a generalised quaternion group. In the next section, we will show that all finitely generated abelian groups have the \cref{secondary-property} and the \cref{tertiary-property}. Thus by \cref{lem:2sylow-sec} and \cref{lem:2sylow-ter}, every finite group whose 2-Sylow subgroup has periodic cohomology has both properties once we have shown them for all generalised quaternion groups.  This is the purpose of this subsection.  In combination with our other results, \cref{theorem:quarterion-groups} is therefore part of the stable diffeomorphism classification of spin $4$-manifolds with finite fundamental group, whose 2-Sylow subgroup has periodic cohomology.

Note that the \cref{secondary-property} for finite groups with quaternion 2-Sylow subgroup was already proved in \cite[Theorem 6.4.1]{teichnerthesis}. Since the statement there is slightly different, we reprove it here for completeness, as an instructive illustration of the use of \cref{cond}.

Let $n$ be a power of two. A presentation of the generalised quaternion group with $8n$ elements is given by
\[Q_{8n}=\langle x,y\mid x^{2n}y^{-2},xyxy^{-1}\rangle.\]
(The quaternion group with $4$ elements is omitted because it is cyclic.)  Note that $xyx=y$ implies $x^{2n}yx^{2n}=y$, so that $y^4=1$ and therefore $x^{4n}=1$.
In the case $n=1$, sending $x \mapsto i$, $y \mapsto j$ gives an isomorphism with the usual presentation of $Q_8$ given by $\langle i,j,k \mid i^2 = j^2 = k^2 = ijk \rangle$.

\begin{theorem}\label{theorem:quarterion-groups}
The generalised quaternion group $Q_{8n}$ has the \cref{secondary-property} and the \cref{tertiary-property}.
\end{theorem}

\begin{proof}
By \cite[Proposition 4.2.1]{teichnerthesis}, the differential
\[d^3_{5,0}\colon H_5(Q_{8n};\Z)\to H_2(Q_{8n};\Z/2)\]
is an isomorphism. In particular, the $\ter$ invariant lies in the trivial group.
 Moreover, this implies that $d^2_{5,0}\colon H_5(Q_{8n};\Z)\to H_3(Q_{8n};\Z/2)$ is the zero map. It follows immediately from the fact that $H_2(Q_{8n};\Z/2)/\im(d^3_{5,0})=0$ that $Q_{8n}$ has the \cref{tertiary-property}.

Now we work on showing that the Secondary Property holds for $Q_{8n}$.
Let $N\in \Z Q_{8n}$ be the norm element $N=\sum_{g\in Q_{8n}}g$ and let $\epsilon\colon\Z Q_{8n}\to \Z$ be the augmentation. By \cite[page 253]{cartan-eilenberg} the beginning of a free resolution of $\Z$ as a $\Z Q_{8n}$-module is given as follows:
\[\xymatrix{
C_4\ar[d]_{d_4}&=&&\Z Q_{8n}\ar[d]^N&\\
C_3\ar[d]_{d_3}&=&&\Z Q_{8n}\ar[dr]^{1-xy}\ar[dl]_{x-1}&\\
C_2\ar[d]_{d_2}&=&\Z Q_{8n}\ar[d]_{\sum_{i=0}^{2n-1}x^i}\ar[drr]_(.25){-y-1}& \oplus &\Z Q_{8n}\ar[d]^{x-1}\ar[dll]^(.25){xy+1}\\
C_1\ar[d]_{d_1}&=&\Z Q_{8n}\ar[dr]_{x-1}& \oplus &\Z Q_{8n}\ar[dl]^{y-1}\\
C_0&=&&\Z Q_{8n}\ar[d]^\epsilon&\\
&&&\Z&
}\]

Let $K$ denote the $2$-dimensional CW complex determined by the start of this resolution (which is the same as the presentation complex for the above presentation of $Q_{8n}$).

One easily computes that $H_3(Q_{8n};\Z/2)\cong \Z/2\otimes_{\Z Q_{8n}}\pi_2(K)\cong\Z/2$.  Recall from the proof of the Tertiary Property for these groups, that $d^2_{5,0}$ is the zero map.  To verify \cref{cond} and therefore show that the \cref{secondary-property} holds, we therefore have to show that the map $A\colon H_3(Q_{8n};\Z/2)\to \wh{H}^{0}(\Sesq(H^2(K;\Z\pi)))$ is injective.

The nontrivial element in $\Z/2\otimes_{\Z Q_8} \pi_2(K)\cong \Z/2\otimes_{\Z Q_8} \ker d_2$ is represented by $1\otimes d_3(1)$. By \cref{lem:pi2-H2-dual} and the definition of $A$ (\cref{defn:map-A}), the form $\lambda:=A(1\otimes d_3(1))$ on
\[H^2(K;\Z Q_{8n})\cong \coker d^2=(\Z Q_{8n})^2/\bigg\langle \Big(\sum_{i=0}^{2n-1}x^{-i},(xy)^{-1}+1\Big),\Big(-y^{-1}-1,x^{-1}-1\Big)\bigg\rangle\] is given by the matrix
\[L:= \begin{pmatrix}
(x^{-1}-1)(x-1)&(x^{-1}-1)(1-xy)\\
(1-(xy)^{-1})(x-1)&(1-(xy)^{-1})(1-xy)
\end{pmatrix}.\]
Note that we used the involution on $\Z Q_{8n}$ to view $\coker d^2$ as a left module.
The entries of the matrix are obtained from the third boundary map $d_3$ in the above resolution.   This matrix a priori defines a pairing on the free module $\Z Q_{8n}^2$.  However the matrix $L$ determines a well-defined pairing on $H^2(K;\Z Q_{8n})\cong \coker d^2$, due to the fact that $d_2 \circ d_3 =0$ in the resolution above.  Throughout our verification of the Secondary Property, we will often write forms on quotient modules as matrices defining forms on free modules, and it will always be necessary that these descend to well defined maps on the quotient modules.

\begin{claim}
The form $\lambda:=A(1\otimes d_3(1))$ is odd.
\end{claim}

To see the claim, we investigate possible forms $q$ that might possibly exhibit $\lambda$ as even, and show that no such $q$ can exist.
A possible $q$ with $q+q^*=\lambda$ can be written as
\[\begin{pmatrix}
1-x+z_1&(x^{-1}-1)(1-xy)+z_2\\
-\overline{z_2}&1-xy+z_3
\end{pmatrix}\]
for some $z_1,z_2,z_3\in \Z Q_{8n}$ satisfying $z_1=-\overline{z_1}, z_3=-\overline{z_3}$. Since \[\big[-y^{-1}-1,x^{-1}-1\big],\Big[\sum x^{-i},(xy)^{-1}+1\Big]=0\in\coker d^2,\]
the form $q$ has to satisfy the relations
\begin{align*}
0 &=q((1,0),(-y^{-1}-1,x^{-1}-1))\\ &= (1-x+z_1)(-y-1)+(x^{-1}-1)(1-xy)(x-1)+z_2(x-1)
\end{align*}
and
\[0=q\big((0,1),\big(\sum x^{-i},(xy)^{-1}+1\big)\big)=-\overline{z_2}\big(\sum x^i\big)+(1-xy+z_3)(xy+1).\]

To derive a contradiction, we will show that these equations cannot be solved after passing to the abelianisation $(Q_{8n})_{ab}\cong(\Z/2)^2$.
We get that $z_1 = -\ol{z_1}$ implies $2z_1 =0$ in $\Z[(\Z/2)^2]$, so $z_1=0$. Similarly $z_3=0$.  Also $x^2=1$ implies that $\sum_{i=1}^{2n-1}x^i=n(1+x)$.
The first relation, with $z_1=0$ substituted, gives
\[0=(1-x)(-y-1)+2(1-x)(1-xy)+z_2(x-1)=(1-x)(1+y)+z_2(x-1) \in \Z[(\Z/2)^2].\]

It follows that $z_2=(1+y)-a(1+x)$ for some $a\in\Z (Q_{8n})_{ab}$, since $(1+x)$ divides any element of $\Z[\Z/2]$ that annihilates $(1-x)$.
Insert $z_2=(1+y)-a(1+x)$ into the second relation to obtain
\begin{align*}
0 &=-n(1+x)(1+y-a-ax)+(1-xy)(1+xy) \\ &=-n(1+x)(1+y-a-ax)=2na(1+x)-nN,
\end{align*}
where $N$ now denotes the norm element in $\Z[(\Z/2)^2]$ of $Q_{8n}^{ab}=(\Z/2)^2$.
But the first summand has all coefficients divisible by $2n$, while the second summand has all coefficients only divisible by $n$. Hence this element does not vanish in $\Z[(\Z/2)^2]$. Hence there are no solutions for $z_1, z_2, z_3$ that define a $q$ as desired and so $\lambda$ cannot be even as claimed.  This completes the proof of the claim, so that \cref{cond} holds for $Q_{8n}$, completing the proof of \cref{theorem:quarterion-groups}.
\end{proof}

\section{The secondary and tertiary properties for abelian groups}\label{sec:abelian-groups}

\noindent In this section we will focus on finitely generated abelian groups.  The goal is to prove the next theorem, classifying spin $4$-manifolds with abelian fundamental groups up to stable diffeomorphism.

\begin{theorem}\label{thm:conjectures-abelian}
The \cref{secondary-property} and the \cref{tertiary-property} hold for all finitely generated abelian groups.
\end{theorem}

The proof breaks up into stages, each of which is considered in one of the next few subsections.
First we consider cyclic groups, then we consider abelian groups with at most two generators, then at most three generators, and finally we use the inheritance properties from \cref{sec:inheritance} to deduce that the properties hold for any finitely generated abelian group.

We will need the following structure of certain cohomology rings of finite cyclic groups.

\begin{thm}[{\cite[Example 3.41]{hatcher}}]
Let $C$ be a cyclic group of order $n=2k$. Then the cohomology ring with $\Z/n$ coefficients is a quotient of the polynomial ring $(\Z/n)[\alpha,\beta]$, as follows:
\[H^*(C;\Z/n)\cong (\Z/n)[\alpha,\beta]/(\alpha^2-k\beta),\]
where $|\alpha|=1$ and $|\beta|=2$.
\end{thm}

\begin{cor}
\label{cor:cup}
Let $C$ be a cyclic group of order $n=2k$. Then the cohomology ring with $\Z/2$ coefficients is a quotient of the polynomial ring $(\Z/2)[\alpha,\beta]$, as follows:
\[H^*(C;\Z/2)\cong\left\{\begin{matrix} \Z/2[\alpha,\beta]/\alpha^2&\text{if }k=2m\\
\Z/2[\alpha]&\text{if }k=2m+1\end{matrix}\right.\]
where $|\alpha|=1$ and $|\beta|=2$.
\end{cor}

Note that the corollary is not obtained simply by setting $n=2$.   Rather, $\alpha^2 = k\beta$ becomes either $\alpha^2=0$ or $\alpha^2=\beta$, when $k=2m$ or $k=2m+1$ respectively.

\subsection{Cyclic groups}
\label{sec:cyclic-gps}
\begin{lemma}\label{lem:cyclic-gps-sec-ter}
The \cref{secondary-property} and the \cref{tertiary-property} hold for all cyclic groups.
\end{lemma}

\subsubsection*{Cyclic groups have the secondary property}\label{section:cyclic-gps-sec-prop}

Let $G$ be a cyclic group.
If the cyclic group $G$ is of odd order or infinite order, then $H_3(G;\Z/2)=0$, and there is nothing further to show.

If $G=\langle T\rangle$ is cyclic of even order $2k$, then we need to check that \cref{cond}.  Note that
\[\red_2 \colon \Z/|G|\cong H_5(G;\Z)\to H_5(G;\Z/2)\cong \Z/2\]
is the projection.  The generator is dual to $\alpha\beta^2$ in the notation of \cref{cor:cup}, if $k=2m$ for some $m$, and the generator is dual to $\alpha^5$ if $k=2m+1$. Similarly the generator of $H_3(G;\Z/2)$ is dual to $\alpha\beta$ when $k=2m$ and $\alpha^3$ when $k=2m+1$ In both cases
\[Sq^2\colon H^3(G;\Z/2)\to H^5(G;\Z/2)\]
is an isomorphism by \cref{cor:cup}, together with a straightforward computation using the axioms of the Steenrod operations. In particular recall the Cartan formula $\Sq^n(xy) = \sum_{p+q=n} \Sq^p(x)\Sq^q(y)$. Thus
\[Sq_2\circ \red_2\colon H_5(G;\Z)\to H_3(G;\Z/2)\]
is surjective.

Therefore, to verify \cref{cond}, we have to show that
\[A\colon H_3(G;\Z/2)\to \widehat{H}^0(\Sesq(H^2(K;\Z G)))\]
is trivial. A free resolution of $\Z$ as a $\Z G$ module is given by
\[\ldots\to \Z G \xrightarrow{N_T}\Z G \xrightarrow{1-T}\Z G \xrightarrow{N_T}\Z G \xrightarrow{1-T}\Z G \xrightarrow{\epsilon}\Z,\]
where $N_T=\sum_{i=0}^{n-1}T^i$. Let $K$ be the corresponding $2$-complex with $\pi_1(K)\cong G$. Then $\Z/2 \otimes_{\Z G} \pi_2(K) \cong \Z/2 \otimes_{\Z G} \ker(N_T)$ is generated by the image of $1$ under $$d_3=(1-T)\colon \Z G\to \ker(N_T) \to \Z/2 \otimes_{\Z G} \ker(N_T).$$
Hence by \cref{lem:pi2-H2-dual}, the form $\lambda$ on $H^2(K;\Z G)$ that we have to consider is given by \[\lambda(x,y)=x(1-T^{-1})(1-T)\overline{y},\] for $x,y\in \Z G /(\ol{N_T}) = \Z G/ (N_T)$. For $q$ given by $q(x,y)=x(1-T)\overline{y}$, we have $\lambda=q+q^*$. Thus $\lambda$ is even and so the map $A$ is trivial as required.  It follows that the sequence of \cref{cond} is exact as required, so that the cyclic group $G$ has the \cref{secondary-property}.

\subsubsection*{Cyclic groups have the tertiary property}\label{sec:cyclic-groups-ter}

By \cref{cor:cup}, we have that $d^2_{4,1}=Sq_2\colon H_4(G;\Z/2)\to H_2(G;\Z/2)$ is an isomorphism. Hence the group in which the $\ter$ invariant resides is trivial, and the \cref{tertiary-property} trivially holds for cyclic groups.

\subsection{Abelian groups with at most two generators}

This section proves the following lemma.

\begin{lemma}\label{lem:abelian-gps-two-gens-sec-ter}
The \cref{secondary-property} and the \cref{tertiary-property} hold for all abelian groups with at most two generators.
\end{lemma}

The proof of this lemma will require about eight pages.
We consider the secondary property first.

\subsubsection*{Two generator abelian groups have the secondary property}\label{section:secondary-ab-gps-two-gens}

\begin{claim}
To prove the secondary property for two generator abelian groups, it suffices to consider the groups $\pi\cong \Z/2^{k_1}\times \Z/2^{k_2}=\langle a,b\mid a^{2^{k_1}},b^{2^{k_2}},[a,b]\rangle$, for some $k_1, k_2\geq 1$, and the groups $\pi\cong \Z \times \Z/2^k$ for some $k \geq 1$.
\end{claim}

The claim follows from the fact that we can pass to a finite odd index subgroup by \cref{lem:2sylow-sec}, and the fact that $\Z \times \Z$ has cohomological dimension two, which we already know to have the Secondary property by \cref{thm:conjs-cohom-dim-three}.  Every abelian group with two generators (as the minimal number of generators) has a finite odd index subgroup that belongs to the list in the claim.  We already proved that cyclic groups have the Secondary property in \cref{section:cyclic-gps-sec-prop}.
This completes the proof of the claim.\qed

For each of the groups in the claim, we need to check \cref{cond}.  This will occupy the next few pages.
\medskip

\noindent\textit{The groups $\Z/2^{k_1}\times \Z/2^{k_2}$.}\\ \smallskip
We begin by considering the groups $\Z/2^{k_1}\times \Z/2^{k_2}$, with $k_1,k_2 \geq 1$.
In this case $H_3(\pi;\Z/2)\cong (\Z/2)^4$. If $\alpha\in H_3(\pi;\Z/2)$ is in the image of the inclusion of one of the two factors, then we know exactness at $\alpha$ by the last subsection. Thus by symmetry we only have to consider the nontrivial element $\gamma\in H_3(\pi;\Z/2)$ coming from the generator of $H_1(\Z/2^{k_1};\Z/2)\otimes H_2(\Z/2^{k_2};\Z/2)\cong \Z/2$.

\begin{claim}
The element $\gamma$ lies in $\im(Sq_2\circ \red_2)$ if and only if $k_1\leq k_2$.
\end{claim}

The proof of this claim will take the next page and a half.
We will show, in the case that $k_1 < k_2$, that $\gamma$ lies in the image of $Sq_2\circ \red_2$, and that the nontrivial element $\ol{\gamma} \in H_2(\Z/2^{k_1};\Z/2)\otimes H_1(\Z/2^{k_2};\Z/2) \subseteq H_3(\pi;\Z/2)$ does not lie in the image.  In the case that $k_1=k_2$, we will show that one, and hence both, of $\gamma$ and $\ol{\gamma}$ lies in $\im(Sq_2\circ \red_2)$.

First we compute the image of the $\Sq^2 \colon H^3(\pi;\Z/2) \to H^5(\pi;\Z/2)$.
A preliminary computation that we will soon need is that $\Sq^1(\beta_i)=0$, with $\beta_i\in H^2(\Z/2^{k_i};\Z/2)$ as in \cref{cor:cup}, and $k_i>1$.  To see this, we use the fact that $\Sq^1$ coincides with the Bockstein homomorphism $BS = \Sq^1 \colon H^2(\Z/2^{k_i};\Z/2) \to H^3(\Z/2^{k_i};\Z/2)$ associated to the coefficient sequence $0 \to \Z/2 \to \Z/4 \to \Z/2 \to 0$ (see \cite[Chapter 3, Theorem 1]{mosher-tangora}).
We have that $H^2(\Z/2^{k_i};\Z/4) =\Z/4 = H^3(\Z/2^{k_i};\Z/4)$.  The differentials in the tensored down resolutions all vanish, so maps in the Bockstein long exact sequence coincide with the maps in the short exact sequence $0 \to \Z/2 \to \Z/4 \to \Z/2 \to 0$ that induces it.  Thus the connecting homomorphism vanishes.    This completes the proof that $\Sq^1(\beta_i)=0$.

To compute the image of $\Sq^2$, we start with the case that $1 < k_1 \leq k_2$. Then $H^3(\Z/2^{k_1}\times \Z/2^{k_2};\Z/2)$ is generated by elements of the form $\alpha_i\beta_j$, in the notation of \cref{cor:cup} with $i,j \in \{1,2\}$. We compute:
\[\Sq^2(\alpha_i\beta_j) = \Sq^2(\alpha_i)\beta_j + \Sq^1(\alpha_i)\Sq^1(\beta_j) + \alpha_i \Sq^2(\beta_j) = \alpha_i\beta_j^2.\]

Next, we consider the case that $k_1=1 < k_2$.  Now $H^3(\Z/2^{k_1}\times \Z/2^{k_2};\Z/2)$ is generated by elements  $\alpha_1\beta_2$, $\alpha_2\beta_2$, $\alpha_1^3$, and $\alpha_1^2 \alpha_2$.
The image of the first two cases are $\alpha_1\beta_2^2$ and $\alpha_2\beta_2^2$ by the computation above.  We also have $\Sq^2(\alpha_1^3) = \alpha_1^5$ and $\Sq^2(\alpha_1^2\alpha_2) = \alpha_1^4\alpha_2$.

Finally we consider the case $k_1=k_2=1$.  In this case $H^3(\Z/2\times \Z/2;\Z/2)$ is generated by elements $\alpha_1^{i}\alpha_2^j$ with $i+j=3$.  An inductive argument using the Cartan formula shows that $\Sq^2(\alpha^i)$ is equal to $\alpha^{i+2}$ if $i\equiv 2,3 \mod{4}$ and $\Sq^2(\alpha^i)=0$ if $i \equiv 0,1 \mod{4}$.  Thus we compute
\[\Sq^2(\alpha_1^2\alpha_2) = \alpha_1^4\alpha_2 ;\; \Sq^2(\alpha_i^3) = \alpha^5_i  \text{; and } \Sq^2(\alpha_1 \alpha_2^2) = \alpha_1 \alpha_2^4.\]
It follows that the map $\Sq_2\colon H_5(\pi;\Z/2) \to H_3(\pi;\Z/2)$ is onto in all cases.
A resolution of $\Z$ by free $\Z\pi$-modules is given below.  To complete the proof of the claim, one has to check whether the elements of $H_5(\pi;\Z/2) = H_5(\Z/2^{k_1}\times \Z/2^{k_2};\Z/2)$, that hit the generators of $H_3(\Z/2^{k_1}\times \Z/2^{k_2};\Z/2)$ under the dual of the Steenrod square, are in the image of the reduction modulo two map $\red_2$.

First we consider the case that $k_1<k_2$.
$(\alpha_1\beta_2^2)^* \in H_5(\Z/2^{k_1}\times \Z/2^{k_2};\Z/2)$ maps to $\gamma$ and $(\alpha_2\beta_1^2)^*$ maps to $\ol{\gamma}$.
The terms in degree $4$, $5$ and $6$ of the $\Z\pi$-module resolution of $\Z$, tensored down over $\Z$, are shown in the diagram below.  Let $k := k_1$ and let $\ell := k_2$, so that $k \leq \ell$ and $\pi = \Z/2^{k} \times \Z/2^{\ell}$. Also denote $C_* := C_*(\pi;\Z\pi)$.  The next diagram shows \[\Z \otimes C_6 \cong \Z^7 \to \Z \otimes C_5 \cong \Z^6 \to \Z \otimes C_4 \cong \Z^5.\]
Ignore the underlining of two $\Z$ summands for now.
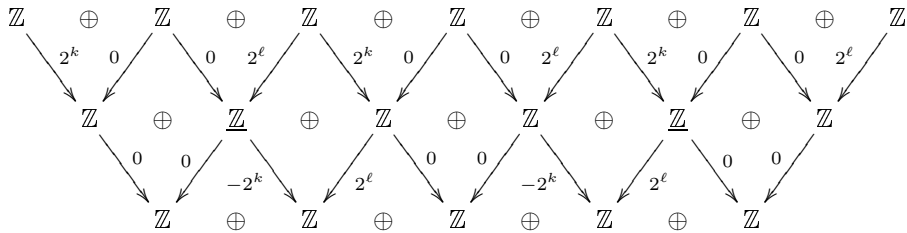
\[\xymatrix @C-0.35cm {\Z \ar[dr]^-{2^{k}} & \oplus & \Z \ar[dl]_-{0} \ar[dr]^-{0} & \oplus & \Z\ar[dl]_-{2^{\ell}} \ar[dr]^-{2^{k}} & \oplus & \Z\ar[dl]_-{0} \ar[dr]^-{0} & \oplus & \Z\ar[dl]_-{2^{\ell}} \ar[dr]^-{2^{k}} & \oplus & \Z\ar[dl]_-{0} \ar[dr]^-{0} & \oplus & \Z\ar[dl]_-{2^{\ell}} \\
 & \Z \ar[dr]^-{0} & \oplus & \underline{\Z} \ar[dl]_-{0} \ar[dr]_>>>>>>{-2^{k}}& \oplus & \Z\ar[dl]^>>>>>>{2^{\ell}} \ar[dr]^-{0}& \oplus & \Z\ar[dl]_-{0} \ar[dr]_>>>>>>{-2^{k}}& \oplus & \underline{\Z}\ar[dl]^>>>>>>{2^{\ell}} \ar[dr]^-{0}& \oplus & \Z\ar[dl]_-{0} & \\
  & & \Z & \oplus & \Z & \oplus & \Z & \oplus & \Z & \oplus & \Z & & &}\]
In $\Z/2 \otimes C_5$, we have
$e_2 = (\alpha_2\beta_1^2)^*$, which maps to $\ol{\gamma} = (\alpha_2\beta_1)^* \in H_3(\pi;\Z/2)$.  Similarly $e_5 = (\alpha_1\beta_2^2)^*$, which maps to $\gamma = (\alpha_1\beta_2)^* \in H_3(\pi;\Z/2)$.
The summands generated by $e_2$ and $e_5$ are underlined in the diagram above.
In the case that one or both of $k,\ell$ are equal to $1$, replace $\beta_i$ by $\alpha_i^2$ in the above statements.  Otherwise the computation is the same.

The relevant generator of $H_5(\pi;\Z)$ for $\ol{\gamma}$ is $2^{k-\ell}e_2 + e_3$.  If $k < \ell$ then this maps to $e_3$ in $\Z/2 \otimes C_5$, and so $\ol{\gamma}$ is not hit by $\Sq_2 \circ \red_2$.  On the other hand, the relevant generator for $H_5(\pi;\Z)$ for $\gamma$ is $e_5 + 2^{k-\ell}e_4$, which maps to $e_5$ in $\Z/2 \otimes C_5$.  Then $e_5$ maps to $\gamma$ under $\Sq_2 \colon H_5(\pi;\Z/2) \to H_3(\pi;\Z/2)$.  Thus when $k<\ell$, we see that $\gamma$ lies in the image of $\Sq_2 \circ \red_2$.
This proves the claim in the case that $k \neq \ell$, that is $k_1 \neq k_2$.

Now consider the case that $k_1=k_2$, that is $k =\ell$.   We just need to show that $\gamma$ does lie in the image of $\Sq_2 \circ \red_2 \colon H_5(\pi;\Z) \to H_3(\pi;\Z/2)$.  The relevant generator of the $\Z$-homology  is $e_4 +e_5$.  This maps to $e_4+e_5$ in $\Z/2 \otimes C_5$.  Now, $e_4$ is dual to $\beta_1\beta_2\alpha_2$ if $k=\ell \geq 2$, or $\alpha_1^2\alpha_2^3$ if $k_1=k_2=1$.  From the computation of the Steenrod square $\Sq^2 \colon H^3(\pi;\Z/2) \to H^5(\pi;\Z/2)$ that we made not long ago, it follows that $e_4 + e_5 \mapsto (\alpha_1\beta_2)^* = \gamma$ under the dual map $\Sq_2 \colon H_5(\pi;\Z/2) \to H_3(\pi;\Z/2)$, as required. Thus the element $\gamma$ lies in $\im(Sq_2\circ \red_2)$ if and only if $k_1\leq k_2$. This completes the proof of the claim.\qed

Now we need to show that the cases in which $\gamma \in \im(\Sq_2 \circ \red_2)$ correspond to the cases in which $A(\gamma)$ is even.
Let $n := 2^{k_1} = 2^k$ and $m := 2^{k_2}=2^{\ell}$.
Let $N_a:=\sum_{i=0}^{n-1}a^i$ and $N_b:=\sum_{i=0}^{m-1}b^i$. From the $\Z\pi$-chain complex
\[\xymatrix{
	\Z\pi\ar[dr]_{1-a}&\oplus&\Z\pi\ar[dl]_{1-b}\ar[dr]^{-N_a}&\oplus&\Z\pi\ar[dl]_{N_b}\ar[dr]^{1-a}&\oplus&\Z\pi\ar[dl]^{1-b}\\
	&\Z\pi\ar[dr]_{N_a}&\oplus&\Z\pi\ar[dl]_{1-b}\ar[dr]^{a-1}&\oplus&\Z\pi\ar[dl]^{N_b}&\\
	&&\Z\pi\ar[dr]_{1-a}&\oplus&\Z\pi\ar[dl]^{1-b}&&\\
	&&&\Z\pi&&&}\]
one computes $H^2(K;\Z\pi)\cong (\Z\pi)^3/\langle(N_a,1-b^{-1},0),(0,a^{-1}-1,N_b)\rangle$, and that $A(\gamma)$ is the form represented by:
\[\begin{pmatrix}
0&0&0\\
0&mN_b&(1-a)N_b\\
0&(1-a^{-1})N_b&2-a-a^{-1}\end{pmatrix}.\]

\begin{claim}
For $k_1 > k_2$, that is $n>m$, the form $A(\gamma)$ is not even.
\end{claim}

A possible $q$ with $A(\gamma)=q+q^*$ has to have the following shape:
\[\begin{pmatrix}
u&v&w\\
-\overline{v}&\tfrac{m}{2}N_b+y&-aN_b+x\\
-\ol{w}&N_b-\overline x&1-a+z\end{pmatrix}\]
with $u,v,w,x,y,z\in\Z\pi$, $u+\ol u=y+\ol y=z+\ol z=0$.
In order to determine a form on $H^2(K;\Z\pi)$ it has to satisfy:
\[0=(N_b-\ol{x})(a-1)+(1-a+z)N_b=(1-a)\ol x+N_bz\]
and
\[0=\tfrac{m}{2}(a-1)N_b+(a-1)y-amN_b+xN_b=-\tfrac{m}{2}(a+1)N_b+(a-1)y+xN_b.\]
Apply the involution to the second equation and then multiply by $(1-a)$ to yield
\begin{align*}
0 &=-\tfrac{m}{2}(1-a)(a^{-1}+1)N_b-(1-a)(a^{-1}-1)y+(1-a)\ol xN_b\\
&=\tfrac{m}{2}(a-a^{-1})N_b+(2-a-a^{-1})y+(1-a)\ol xN_b\\
&= \tfrac{m}{2}(a-a^{-1})N_b+(2-a-a^{-1})y-mzN_b.
\end{align*}
Here we used $\ol y=-y$, while for the last equality we used $(1-a)\ol x=-N_bz$ as well as $N_b^2=mN_b$.
Reduce coefficients modulo $m$ and take the coefficient of $a$ to obtain:
\[0=\frac{m}{2}+1y_a-y_0-y_{a^2}\mod m.\]
Since $y=-\ol y$ we have $y_0=0$ and thus $y_{a^2}=2y_a+\tfrac{m}{2}$.
Take the coefficient of $a^k$ for $1<k\leq \tfrac{n}{2}$ to obtain:
\[0=2y_{a^k}-y_{a^{k-1}}-y_{a^{k+1}}\mod m.\]
Thus by induction we have:
\begin{equation}\label{eq:y}
y_{a^k}=ky_a+(-1)^{\delta_{4|k}}\delta_{2|k}\tfrac{m}{2}\mod m,\end{equation}
where $\delta_{2|k}$ and $\delta_{4|k}$ are the characteristic functions of the Boolean statements $2|k$ and $4|k$ respectively. That is, for example, $\delta_{2|k} =1$ if $2$ divides $k$, and $\delta_{2|k} =0$ otherwise.
If $n=2^km$ for $k\geq 1$, then $m|\tfrac{n}{2}$ and
\[y_{a^{n/2}}=\pm\tfrac{m}{2}\mod m\]
But since $y=-\ol y$, the coefficient $y_{a^{n/2}}$ has to be zero. This shows the claim that the form $A(\gamma)$ cannot be even if $n=2^k m$ for some $k\geq 1$, i.e.\ if $k_1>k_2$.\qed

\begin{claim}
The form $A(\gamma)$ is even if $k_1\leq k_2$.
\end{claim}
First we suppose that $k_1=k_2$, that is $n=m$.
Let $y$ be such that $y_{a^kb^l}=k-\tfrac{n}{2}$ for $k\neq0$. First note that $y_{a^{-k}b^{-l}}=(n-k)-\tfrac{n}{2}=-(k-\tfrac{n}{2})=-y_{a^kb^l}$. Therefore, $\ol y=-y$. Consider the form $q$ given by
\[\begin{pmatrix}
0&0&-(1-b^{-1})\\
0&\tfrac{m}{2}N_b+y&-aN_b+N_a\\
(1-b)&N_b-N_a&1-a\end{pmatrix}\]
It follows from $N_a=\ol{N_a}$ and $\ol y=-y$ that $q+q^*=A(\gamma)$. It remains to show that $q$ is indeed a form on $H^2(K;\Z\pi)$. For this it has to satisfy the equations:
\begin{align}
0&=0\cdot N_a+0\cdot(1-b)\nonumber\\
0&=0\cdot N_a+(\tfrac{m}{2}N_b+y)(1-b)\label{eq:twogen2}\\
0&=(1-b)\cdot N_a+(N_b-N_a)(1-b)\nonumber\\
0&=0\cdot(a-1)-(1-b^{-1})N_b\nonumber\\
0&=(\tfrac{m}{2}N_b+y)(a-1)-(aN_b+N_a)N_b\label{eq:twogen5}\\
0&=N_b(a-1)+(1-a)N_b\nonumber
\end{align}
Equation \eqref{eq:twogen2} is true since $y$ is a multiple of $N_b$ by definition. All but Equation \eqref{eq:twogen5} are satisfied, using $(1-b)N_b=(1-b^{-1})N_b=0$. To check \eqref{eq:twogen5}, we have to show that
\[0=(\tfrac{m}{2}N_b+y)(a-1)-(aN_b+N_a)N_b=-\tfrac{m}{2}(a+1)N_b+(a-1)y-N_aN_b.\]
Since each term is a multiple of $N_b$, it suffices to show that this holds for the coefficients of $a^k$, for every $k\geq 0$. For $a^0$ we have
\[0=-\tfrac{m}{2}+y_{a^{-1}}-y_0+1=-\tfrac{m}{2}+m-1-\tfrac{m}{2}-0+1=0.\]
For $a^1$ we get
\[0=-\tfrac{m}{2}+y_{0}-y_a+1=-\tfrac{m}{2}+0-1+\tfrac{m}{2}+1=0.\]
Finally for $a^k$ with $k>1$, we have
\[0=0+y_{a^{k-1}}-y_{a^k}+1=k-1-\tfrac{m}{2}-k+\tfrac{m}{2}+1=0.\]
This completes the verification of Equation \eqref{eq:twogen5}. This completes the proof of the claim that $A(\gamma)$ is even when $k_1=k_2$.

If $k_1 < k_2$, consider the projection $\Z/2^{k_2}\times\Z/2^{k_2}\to\Z/2^{k_1}\times \Z/2^{k_2}$ applied to the form $q$, to see that the form $A(\gamma)$ also has to be even in this case.
This completes the proof of the claim that $A(\gamma)$ is even when $k_1 \leq k_2$, completing the proof that the groups $\Z/s^{k_1} \times \Z/2^{k_2}$ have the \cref{secondary-property}.\qed

\medskip

\noindent\textit{The groups $\Z \times \Z/2^{k}$.}\\ \smallskip
Next, we consider the group $\pi = \Z \times \Z/2^k = \langle a,b \mid [a,b], b^{2^k}\rangle$.
Once again we consider the element $\gamma \in H_3(\pi;\Z/2)$ coming from $H_1(\Z;\Z/2) \otimes H_2(\Z/2^{k};\Z/2)$.

\begin{claim}
The form $A(\gamma)$ is not even.
\end{claim}
A resolution of $\Z$ by $\Z\pi$ modules can be obtained from the tensor product
\[(\Z\pi \xrightarrow{1-a} \Z\pi\big) \otimes \big(\Z\pi \otimes_{\Z[\Z/2]} C_*(B\Z/2;\Z[\Z/2])\big).\]
This is a subcomplex of the resolution depicted above in the case $\Z/2^{k_1} \times \Z/2^{k_2}$, obtained by deleting everything apart from the two right-most summands of each chain group.
We compute that
$H^2(K;\Z\pi)\cong (\Z\pi)^2/\langle(1-b^{-1},0),(a^{-1}-1,N_b)\rangle$ and that $A(\gamma)$ is the form
\[\begin{pmatrix}
2^k N_b&(1-a)N_b\\
(1-a^{-1})N_b&2-a-a^{-1}\end{pmatrix}.\]
Now let $m=2^k$.  A possible $q$ with $A(\gamma)=q+q^*$ has to have the following shape:
\[\begin{pmatrix}
\tfrac{m}{2}N_b+y&-aN_b+x\\
N_b-\overline x&1-a+z\end{pmatrix}\]
with $x,y,z\in\Z\pi$, $y+\ol y=z+\ol z=0$.

A verbatim repetition of the argument for $\Z/2^{k_1}\times \Z/2^{k_2}$ in the case $k_1>k_2$, starting at ``In order to determine a form on $H^2(K;\Z\pi)$ it has to satisfy,'' shows that \eqref{eq:y} also has to hold in the case $\Z\times \Z/m = \Z \times \Z/2^k$.  Since in this case $a$ has infinite order, \eqref{eq:y} implies that infinitely many coefficients of $y$ are non-zero. This is a contradiction to $y$ being an element of the group ring. Hence the form $A(\gamma)$ cannot be even for $\Z\times \Z/2^k$.   This completes the proof of the claim that $A(\gamma)$ is not even.\qed

\begin{claim}
The element $\gamma$ is not in the image of $\Sq_2\circ\red_2$.
\end{claim}
It is straightforward to compute that the map
\[\red_2 \colon H_5(\Z \times \Z/2^k;\Z) \cong \Z/2^k \to H_5(\Z \times \Z/2^k;\Z/2) \cong \Z/2 \oplus \Z/2\]
is given by $(0,\pr)$, where $\pr \colon \Z/2^{k} \to \Z/2$ is the reduction modulo two. The generator of the image of this map is dual to $\alpha_2\beta_2^2$.
For $k>1$, the cohomology group $H^3(\Z/2^{k};\Z/2)\cong \Z/2 \oplus \Z/2$ is generated by $\alpha_2\beta_2$ and $\alpha_1\beta_2$; recall that $H^*(\Z/2^k;\Z/2) = (\Z/2)[\alpha_2,\beta_2]/\alpha_2^2$, and note that $H^*(\Z;\Z/2) = \Z[\alpha_1]/\alpha_1^2$, where $\alpha_1$ has degree $1$.
But $\Sq^2(\alpha_2\beta_2) = \alpha_2\beta_2^2$ and $\Sq^2(\alpha_1 \beta_2) = \alpha_1\beta_2^2$. It follows that \[\Sq_2((\alpha_2\beta_2^2)^*) = \alpha_2\beta_2 \neq \alpha_1\beta_2 = \gamma.\]
The same argument goes through in the case that $k=1$ if we replace $\beta_2$ by $\alpha_2^2$ in the above computations.
This completes proof of the claim that $\gamma$ is not in the image of $\Sq_2\circ\red_2$. \qed

We have therefore completed the proof that $\Z\times\Z/2^k$ satisfies \cref{cond} and therefore has the Secondary Property. Since this is the last case we had to check, this completes the proof that abelian groups with at most two generators have the \cref{secondary-property}.\qed

\subsubsection*{Two generator abelian groups have the tertiary property}
\label{sec:twogenter}
Let us now consider the \cref{tertiary-property}. From the discussion of the Tertiary property for cyclic groups (\cref{sec:cyclic-groups-ter}), it follows that the $\ter$ invariant of an abelian group with two generators either lives in the trivial group or in $\Z/2$.
In more detail, consider a group $\pi := C_1 \times C_2$, where $C_1$ and $C_2$ are cyclic groups.
We have that \[H_2(\pi;\Z/2) \cong \bigoplus_{i=0}^2 H_i(C_1;\Z/2) \otimes H_{2-i}(C_2;\Z/2).\]
However the parts with $i=0$ and $i=2$ lie in the image of $\Sq_2 =d^2_{4,0} \colon H_4(\pi;\Z/2) \to H_2(\pi;\Z/2)$ by the structure of the cohomology ring of cyclic groups (\cref{cor:cup}).  Therefore if the $\ter$ invariant is nontrivial, it lives in $H_1(C_1;\Z/2) \otimes H_1(C_2;\Z/2) \cong \Z/2$.
In the case that the $\ter$ invariant lives in $\Z/2$, the nontrivial element corresponds to the commutator relation (i.e.\ to the 2-cell in a model for the classifying space inducing that relation).

By \cite[Proposition~6.1]{KPR}, $\Z \otimes_{\Z\pi} \Gamma(\pi_2(K))$ is torsion-free if $\pi$ is a finite abelian group with two generators. Hence the \cref{tertiary-property} holds for these groups by \cref{cor:ter-finite-group}.

Any infinite abelian group with two generators admits a surjection onto $\Z/8\times \Z/2$, mapping the commutator relation to the commutator relation. By the discussion at the beginning of this subsection, this surjection induces an injection on the $\ter$ invariant if the $\ter$ invariant for $\Z/8\times \Z/2$ is nontrivial. Hence by the Inheritance \cref{thm:tersurj}, it suffices to show that the $\ter$ invariant for $\Z/8\times \Z/2$ is nontrivial.

The commutator relation is dual to the product $\alpha_1 \alpha_2\in H^2(\Z/8 \times \Z/2;\Z/2)$ where $\alpha_1$ and $\alpha_2$ are the generators of $H^1$ of the two cyclic subgroups. Since $\Z/8$ has order higher than two, $\alpha_1^2$ is trivial. In particular, $Sq^2(\alpha_1 \alpha_2)=\alpha_1^2\alpha_2^2=0$, and the commutator relation is not in the image of the dual of~$Sq^2$.

Now consider the commutative square
\begin{equation}
\label{eqn:diagram-two-claims}
\xymatrix{
	\Z/2\cong H_2(\Z/8 \times \Z/2;\Z/2)/\im d^2_{4,1} \ar[r]^-{\cong}_-{p_*} & \Z/2\cong H_2(\Z/4\times\Z/2;\Z/2)/\im d^2_{4,1}\\
	\ker d^2_{5,0}\subseteq H_5(\Z/8 \times \Z/2;\Z)\ar[r]_-{p_*}\ar[u]^{d^3_{5,0}}&\ker d^2_{5,0}\subseteq H_5(\Z/4\times \Z/2;\Z), \ar[u]^{d^3_{5,0}}}
\end{equation}
where $p\colon \Z/8\times\Z/2\to \Z/4\times\Z/2$ is the projection.

\begin{claim}
	The map $p_*$ in the bottom row of the diagram \eqref{eqn:diagram-two-claims} is trivial on $\ker d^2_{5,0}\subseteq H_5(\Z/8\times\Z/2;\Z)$.
\end{claim}

Using the claim, it follows by commutativity that $d^3_{5,0}\colon \ker d^2_{5,0}\subseteq H_5(\Z/8\times\Z/2;\Z)\to \Z/2 \cong H_2(\Z/8\times\Z/2;\Z/2)/\im d^2_{4,1}$ is trivial.  Thus the $\ter$ invariant for $\Z/8\times\Z/2$ is nontrivial.  This completes the proof that abelian groups with at most two generators have the \cref{tertiary-property}, modulo the proof of the claim.

Before we prove the claim we give a computation of the maps between group homology that is needed.
Consider the projection map $p \colon \Z/2^n \to \Z/2^{n-1}$, and view $\Z[\Z/2^{n-1}]$ as a module over $\Z[\Z/2^n]$ via this projection.  Write $T$ for the generator of $\Z/2^n$ and $U$ for the generator of $\Z/2^{n-1}$, and note that $p_*(T)=U$.  The induced map on resolutions, shown in degrees up to four, is as follows:
\[\xymatrix{ \Z[\Z/2^n] \ar[r]^{N_{2^n}} \ar[d]^{4\cdot p_*} & \Z[\Z/2^n] \ar[r]^{1-T} \ar[d]^{2 \cdot p_*} & \Z[\Z/2^n] \ar[r]^{N_{2^n}} \ar[d]^{2 \cdot p_*} & \Z[\Z/2^n] \ar[r]^{1-T} \ar[d]^{p_*} & \Z[\Z/2^n] \ar[d]^{p_*} &  \\
	\Z[\Z/2^{n-1}] \ar[r]^-{N_{2^{n-1}}} & \Z[\Z/2^{n-1}] \ar[r]^{1-U} & \Z[\Z/2^{n-1}] \ar[r]^-{N_{2^{n-1}}} & \Z[\Z/2^{n-1}] \ar[r]^{1-U} & \Z[\Z/2^{n-1}],
}\]
where $N_{2^j} = 1+ S + \cdots + S^{2^{j} - 1}$ is the norm element of $\Z/2^{j}$ for $(S,j) \in \{(T,n),(U,n-1)\}$.
This can be seen by starting with the identity map on $C_0$, and then extending in order to make the diagram commute.  After tensoring with $\Z$, the induced vertical maps are unchanged: multiplication by $2^{i}$ in degrees $2i$ and $2i+1$.

\medskip

\noindent\textit{Proof of the claim.}\\ \smallskip
Now we embark upon the proof of the claim, which is elementary but giving a detailed computation will require most of the next two pages. We need to show that $p_*\colon H_5(\Z/8 \times \Z/2;\Z) \to H_5(\Z/4 \times \Z/2;\Z)$ is trivial on $\ker d^2_{5,0}\subseteq H_5(\Z/8 \times \Z/2;\Z)$.
Here we recall the part of the chain complex
\[C_6(\Z/8 \times \Z/2;\Z) \cong \Z^7 \to C_5(\Z/8 \times \Z/2;\Z) \cong \Z^6 \to C_4(\Z/8 \times \Z/2;\Z) \cong \Z^5\]
relevant to computing the fifth homology:
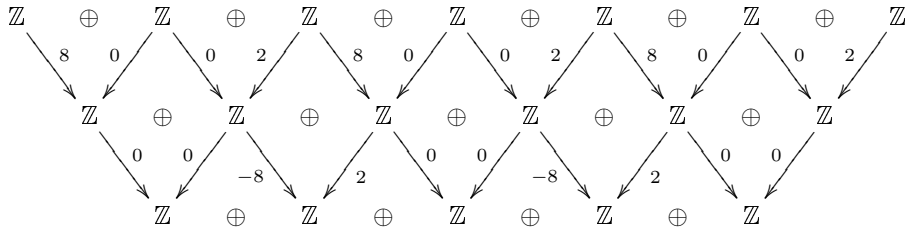
\[\xymatrix @C-0.35cm {\Z \ar[dr]^-{8} & \oplus & \Z \ar[dl]_-{0} \ar[dr]^-{0} & \oplus & \Z\ar[dl]_-{2} \ar[dr]^-{8} & \oplus & \Z\ar[dl]_-{0} \ar[dr]^-{0} & \oplus & \Z\ar[dl]_-{2} \ar[dr]^-{8} & \oplus & \Z\ar[dl]_-{0} \ar[dr]^-{0} & \oplus & \Z\ar[dl]_-{2} \\
	& \Z \ar[dr]^-{0} & \oplus & \Z \ar[dl]_-{0} \ar[dr]_>>>>>>{-8}& \oplus & \Z\ar[dl]^>>>>>>{2} \ar[dr]^-{0}& \oplus & \Z\ar[dl]_-{0} \ar[dr]_>>>>>>{-8}& \oplus & \Z\ar[dl]^>>>>>>{2} \ar[dr]^-{0}& \oplus & \Z\ar[dl]_-{0} & \\
	& & \Z & \oplus & \Z & \oplus & \Z & \oplus & \Z & \oplus & \Z & & &}\]
Let $e_0,\dots,e_5$ be the basis for $C_5(\Z/8 \times \Z/2;\Z)$ generating the summands shown in order.  So $e_i \in C_{5-i}(\Z/8;\Z) \otimes C_i(\Z/2;\Z)$. We see that
\[\ker(C_5(\Z/8 \times \Z/2;\Z) \to C_4(\Z/8 \times \Z/2;\Z)) = \langle e_0,e_1+4e_2,e_3 + 4e_4, e_5 \rangle.\]
Meanwhile
\[\im(C_6(\Z/8 \times \Z/2;\Z) \to C_5(\Z/8 \times \Z/2;\Z)) = \langle 8e_0,2e_1+8e_2,2e_3 + 8e_4, 2e_5 \rangle.\]
Therefore
\[H_5(\Z/8 \times \Z/2;\Z) \cong \Z/8 \oplus (\Z/2)^3 = \langle e_0,e_1+4e_2,e_3 + 4e_4, e_5 \rangle.\]
On the other hand part of the chain complex
\[C_6(\Z/4 \times \Z/2;\Z) \cong \Z^7 \to C_5(\Z/4 \times \Z/2;\Z) \cong \Z^6 \to C_4(\Z/4 \times \Z/2;\Z) \cong \Z^5\]
relevant to computing the fifth homology is:
\[\xymatrix @C-0.35cm {\Z \ar[dr]^-{4} & \oplus & \Z \ar[dl]_-{0} \ar[dr]^-{0} & \oplus & \Z\ar[dl]_-{2} \ar[dr]^-{4} & \oplus & \Z\ar[dl]_-{0} \ar[dr]^-{0} & \oplus & \Z\ar[dl]_-{2} \ar[dr]^-{4} & \oplus & \Z\ar[dl]_-{0} \ar[dr]^-{0} & \oplus & \Z\ar[dl]_-{2} \\
	& \Z \ar[dr]^-{0} & \oplus & \Z \ar[dl]_-{0} \ar[dr]_>>>>>>{-4}& \oplus & \Z\ar[dl]^>>>>>>{2} \ar[dr]^-{0}& \oplus & \Z\ar[dl]_-{0} \ar[dr]_>>>>>>{-4}& \oplus & \Z\ar[dl]^>>>>>>{2} \ar[dr]^-{0}& \oplus & \Z\ar[dl]_-{0} & \\
	& & \Z & \oplus & \Z & \oplus & \Z & \oplus & \Z & \oplus & \Z & & &}\]
Let $f_0,\dots,f_5$ be the basis for $C_5(\Z/4 \times \Z/2;\Z)$ generating the summands shown in order. We see that
\[\ker(C_5(\Z/4 \times \Z/2;\Z) \to C_4(\Z/4 \times \Z/2;\Z)) = \langle f_0,f_1+2f_2,f_3 + 2f_4, f_5 \rangle.\]
We also have
\[\im(C_6(\Z/4 \times \Z/2;\Z) \to C_5(\Z/4 \times \Z/2;\Z)) = \langle 4f_0,2f_1+4f_2,2f_3 + 4f_4, 2f_5 \rangle.\]
Therefore
\[H_5(\Z/4 \times \Z/2;\Z) \cong \Z/4 \oplus (\Z/2)^3 = \langle f_0,f_1+2f_2,f_3 + 2f_4, f_5 \rangle.\]
The induced map on group homology sends $e_i \mapsto 4 f_i$ for $i=0,1$, $e_i \mapsto 2 f_i$ for $i=2,3$, and $e_i \mapsto f_i$ for $i=4,5$.  We therefore compute that the induced map \[p_* \colon H_5(\Z/8 \times \Z/2 ;\Z) \cong \Z/8 \oplus (\Z/2)^3 \to H_5(\Z/4 \times \Z/2;\Z) \cong \Z/4 \oplus (\Z/2)^3\] has kernel $\langle e_0,e_1+4e_2,e_3 + 4e_4\rangle \cong \Z/8 \oplus (\Z/2)^2$, and is an isomorphism between the last $\Z/2$ summands of each, so $e_5 \mapsto f_5$.

Next we need to compute the kernel of $d^2_{5,0} = \Sq_2 \circ \red_2$.  We start by looking at $\red_2 \colon H_5(\Z/8\times\Z/2;\Z) \to H_5(\Z/8\times\Z/2;\Z/2)$.  After tensoring with $\Z/2$, all the maps in the chain complex above become zero maps.  Therefore $H_5(\Z/8 \times \Z/2;\Z/2) \cong (\Z/2)^6$.  Denote the basis elements by $\wt{e}_i$, for $i=0,\dots,6$.  The map $\Z/8 \oplus (\Z/2)^3  \to (\Z/2)^6$ is given by projection $\Z/8 \to \Z/2$ on the first summand, and maps the second, third and fourth summands isomorphically to the second, fourth and sixth $\Z/2$ summands of $H_5(\Z/8 \times \Z/2;\Z/2) \cong (\Z/2)^6$ respectively.  The kernel is thus $2e_0$.   It remains to compute the dual of the Steedrod square.  First we look at the map
\[\Sq^2 \colon H^3(\Z/8 \times \Z/2;\Z/2) \cong (\Z/2)^4 \to H^5(\Z/8 \times \Z/2 ;\Z/2),\]
then we take its dual. Recall that $H^*(\Z/8;\Z/2) \cong (\Z/2)[\alpha_1,\beta_1]/(\alpha_1^2)$ and $H^*(\Z/2;\Z/2) \cong (\Z/2)[\alpha_2]$, where $\alpha_1$ and $\alpha_2$ have degree one and $\beta_1$ has degree two.
We computed earlier that:
\begin{align*}
\Sq^2(\alpha_2^3) &= \alpha_2^5 = \wt{e}_5^* \\
\Sq^2(\alpha_2^2\alpha_1) &= \alpha_2^4\alpha_1 = \wt{e}_4^* \\
\Sq^2(\alpha_2\beta_1) &= \alpha_2\beta_1^2 = \wt{e}_1^* \\
\Sq^2(\alpha_1\beta_1) &= \alpha_1\beta_1^2 = \wt{e}_0^*.
\end{align*}
Therefore \[\ker(\Sq_2 \colon (\Z/2)^6 \to (\Z/2)^4) = \langle \wt{e}_2, \wt{e}_3\rangle.\]
It follows from combining the maps $\Sq_2$ and $\red_2$ that \[\ker d^2_{4,0} = \langle 2e_0, e_3 + 4e_4 \rangle \cong \Z/4 \oplus \Z/2 \subset \Z/8 \oplus (\Z/2)^3,\]
where the $\Z/2$ in the domain includes into the second $\Z/2$ summand in the codomain.  This is contain in $\ker p_*$ as claimed.  This completes the proof of the claim: the map $p_*$ in the bottom row of diagram \eqref{eqn:diagram-two-claims} is trivial on $\ker d^2_{5,0}\subseteq H_5(\Z/8\times\Z/2;\Z)$.\qed

We have therefore completed the proof that abelian groups with at most two generators have the \cref{tertiary-property}, which completes the proof of \cref{lem:abelian-gps-two-gens-sec-ter}.\qed

\subsection{Abelian groups with at most 3 generators}

\begin{lemma}\label{lem:ab-groups-3-gens-sec-ter}
The \cref{secondary-property} and the \cref{tertiary-property} hold for all abelian groups with at most three generators.
\end{lemma}

For the Tertiary property, \cref{cor:terabelian} implies that knowing the property for abelian groups with at most two generators is sufficient.  Therefore we do not need to consider the Tertiary property in this subsection.
The proof for the Secondary property will take about five pages.

\subsubsection*{Three generator abelian groups have the secondary property}

Once again, we can pass to finite odd index subgroups by \cref{lem:2sylow-sec}.  Thus it suffices to consider groups $G$ of the form $\Z^3$, $\Z \times \Z  \times \Z/2^k$, $\Z \times \Z/2^{k_1} \times \Z/2^{k_2}$ or $\Z/2^{k_1} \times \Z/2^{k_2} \times \Z/2^{k_3}$.
We already know the Secondary property for $\Z^3$ by \cref{thm:conjs-cohom-dim-three}. We will compute that the \cref{secondary-property} holds for the groups $\pi = \Z/2\times \Z/2 \times \Z/2$ and $\pi =  \Z/2\times \Z/2 \times \Z/4$, in a fashion similar to the computation in \cref{section:secondary-ab-gps-two-gens}, and then we will deduce that this implies the Secondary property for all the groups $G$ just listed.
\medskip

\noindent\textit{The group $\pi = \Z/2 \times \Z/2 \times \Z/2$.}\\ \smallskip
If $\pi\cong \Z/2\times\Z/2\times \Z/2$, in order to check \cref{cond}, we only have to consider the element $\gamma\in H_3(\pi;\Z/2)$ coming from the generator of $H_1(\Z/2;\Z/2)\otimes H_1(\Z/2;\Z/2)\otimes H_1(\Z/2;\Z/2)$, since all other elements come from a direct summand with at most $2$ generators.

\begin{claim}
The element $\gamma$ is in the image of $\Sq_2\circ\red_2$.
\end{claim}

The $\Z/2$-coefficient cohomology of $\pi = \Z/2 \times \Z/2 \times \Z/2$ is
\begin{align*}
H^*(\Z/2 \times \Z/2 \times \Z/2;\Z/2) &\cong H^*(\Z/2;\Z/2) \otimes H^*(\Z/2;\Z/2) \otimes H^*(\Z/2;\Z/2) \\ &\cong \Z/2[\alpha_1] \otimes \Z/2[\alpha_2] \otimes \Z/2[\alpha_3],\end{align*}
where $\alpha_i$ has degree one for $i=1,2,3$.
Computing with Steenrod squares, we obtain:
\begin{align*}
\Sq^2(\alpha_1 \alpha_2 \alpha_3) &= \Sq^2(\alpha_1 \alpha_2) \alpha_3 + \Sq^1(\alpha_1 \alpha_2) \Sq^1(\alpha_3) \\ &= \alpha_1^2\alpha_2^2\alpha_3 + \Sq^1(\alpha_1) \alpha_2 \alpha_3^2 + \alpha_1 \Sq^1(\alpha_2) \alpha_3^2 \\  &= \alpha_1^2\alpha_2^2\alpha_3 + \alpha_1^2 \alpha_2 \alpha_3^2 + \alpha_1 \alpha_2^2 \alpha_3^2 \in H^5(\pi;\Z/2). \end{align*}
Let $a_i := \alpha_i^* \in H_1(\Z/2;\Z/2)$ be the dual to $\alpha_i$ and let $a_i^2 = (\alpha_i^2)^* \in H_2(\Z/2;\Z/2)$ be the dual to $\alpha_i^2$.
From the computation of $\Sq^2$ above, we see that the element $\gamma = a_1 \otimes a_2 \otimes a_3$ lies in the image of $\Sq_2\colon H_5(\pi;\Z/2)\to H_3(\pi;\Z/2)$.  We therefore if suffices to show that
\[(\alpha_1^2\alpha_2^2\alpha_3 + \alpha_1^2 \alpha_2 \alpha_3^2 + \alpha_1 \alpha_2^2 \alpha_3^2)^* = a_1^2 \otimes a_2^2 \otimes a_3 + a_1^2 \otimes a_2 \otimes a_3^2 + a_1 \otimes a_2^2 \otimes a_3^2 \]
lies in the image of $\red_2 \colon H_5(\pi;\Z) \to H_5(\pi;\Z/2)$.
Recall that a $\Z[\Z/2]$-module resolution of $\Z$, tensored down over $\Z$, is given by
\[\Z(6) \xrightarrow{2} \Z(5) \xrightarrow{0} \Z(4) \xrightarrow{2}\Z(3) \xrightarrow{0}\Z(2) \xrightarrow{2}\Z(1) \xrightarrow{0} \Z(0),\]
where the parenthetical numbers indicate the grading.
It is then straightforward to compute that this element is a cycle with integer coefficients.  Here we will also use the notation $a_i$ and $a_i^2$ as above, but now this notation represents the corresponding integral chains in $C_*(\Z/2;\Z)$  i.e.\ generators of $\Z(1)$ and $\Z(2)$ respectively.  We have the following computation in $C_*(\pi;\Z) = \bigotimes^3 C_*(\Z/2;\Z)$.
\begin{align*}
\partial(a_1^2 \otimes a_2^2 \otimes a_3 + a_1^2 \otimes a_2 \otimes a_3^2 + a_1 \otimes a_2^2 \otimes a_3^2) & = 2 a_1 \otimes a_2^2 \otimes a_3 +2 a_1^2 \otimes a_2 \otimes a_3 \\  &+ 2 a_1 \otimes a_2 \otimes a_3^2 - 2 a_1^2 \otimes a_2 \otimes a_3 \\ & -  2 a_1 \otimes a_2 \otimes a_3^2 - 2a_1 \otimes a_2^2 \otimes a_3 \\ &= 0.
\end{align*}
If $a_1^2 \otimes a_2^2 \otimes a_3 + a_1^2 \otimes a_2 \otimes a_3^2 + a_1 \otimes a_2^2 \otimes a_3^2$ were a boundary over $\Z$, then it would be a boundary over $\Z/2$ as well.  Therefore it represents an element of $H_5(\pi;\Z)$, that maps to $(\alpha_1^2\alpha_2^2\alpha_3 + \alpha_1^2 \alpha_2 \alpha_3^2 + \alpha_1 \alpha_2^2 \alpha_3^2)^*$ under the reduction modulo two.  As we have shown above that this element maps to $\gamma$ under $\Sq_2$, this completes the proof of the claim that $\gamma$ is in the image of $\Sq_2\circ\red_2$.\qed

\begin{claim}
The form $A(\gamma)$ is even.
\end{claim}
We have
\begin{align*}
H^2(K;\Z\pi)\cong (\Z\pi)^6/\langle (1-b,c-1,0,1+a,0,0),&(a-1,0,1-c,0,1+b,0),\\ &(0,1-b,a-1,0,0,1+c)\rangle
\end{align*}
and the form $A(\gamma)$ is given by
\[\begin{pmatrix}
2(1-c)&(1-c)(1-b)&(1-c)(1-a)&0&0&0\\
(1-b)(1-c)&2(1-b)&(1-a)(1-c)&0&0&0\\
(1-a)(1-c)&(1-a)(1-b)&2(1-a)&0&0&0\\
0&0&0&0&0&0\\
0&0&0&0&0&0\\
0&0&0&0&0&0\end{pmatrix}\]
For $q$ given by
\[\begin{pmatrix}
1-c&-(c+b)&-(c+a)&0&0&0\\
1+cb&1-b&-(b+a)&0&0&0\\
1+ca&1+ab&1-a&0&0&0\\
0&0&0&0&0&0\\
0&0&0&0&0&0\\
0&0&0&0&0&0
\end{pmatrix}\]
we have $\lambda=q+q^*$. This completes the proof of the claim that $A(\gamma)$ is even.\qed

This therefore completes the proof of the \cref{secondary-property} for $\pi = \Z/2\times\Z/2\times\Z/2$.
\medskip

\noindent\textit{The group $\pi = \Z/2 \times \Z/2 \times \Z/4$.}\\ \smallskip
Now we consider $\pi:=\Z/2\times \Z/2\times \Z/4$.  Verifying \cref{cond} for this group will occupy the next two and a half pages.

\begin{claim}
The element $\gamma\in H_3(\pi;\Z/2)$ coming from the generator of $$H_1(\Z/2;\Z/2)\otimes H_1(\Z/2;\Z/2)\otimes H_1(\Z/4;\Z/2)\cong\Z/2$$ is not in the image of $\Sq_2\circ\red_2$.
\end{claim}

The $\Z/2$-coefficient cohomology of $\pi = \Z/2 \times \Z/2 \times \Z/4$ is
\begin{align*}
H^*(\Z/2 \times \Z/2 \times \Z/4;\Z/2) &\cong H^*(\Z/2;\Z/2) \otimes H^*(\Z/2;\Z/2) \otimes H^*(\Z/4;\Z/2) \\ &\cong \Z/2[\alpha_1] \otimes \Z/2[\alpha_2] \otimes \Z/2[\alpha_3,\beta_3]/\alpha_3^2,\end{align*}
where $\alpha_i$ has degree one for $i=1,2,3$ and $\beta_3$ has degree two.
Note that $\gamma = (\alpha_1\alpha_2\alpha_3)^*$.
We first compute the image of the Steenrod square map $\Sq^2 \colon H^3(\pi;\Z/2) \to H^5(\pi;\Z/2)$.
First we have
\begin{align*} \Sq^2(\alpha_1\alpha_2\alpha_3) &= \Sq^2(\alpha_1) \alpha_2 \alpha_3 + \alpha_1^2 \Sq^1(\alpha_2 \alpha_3) + \alpha_1 \Sq^2(\alpha_2\alpha_3) \\
& = \alpha_1^2\alpha_2^2 \alpha_3 + \alpha_1 \alpha_2\alpha_3^2 + \alpha_1 \alpha_2^2\alpha_3^2 \\ &=  \alpha_1^2\alpha_2^2 \alpha_3 \end{align*}
since $\alpha_3^2=0$.
We also compute:
\[\Sq^2(\alpha_1^2 \alpha_2) = \alpha_1^4 \alpha_2; \; \Sq^2(\alpha_1 \alpha_2^2) = \alpha_1 \alpha_2^4; \; \Sq^2(\alpha_i^2\alpha_3) = \alpha_i^4\alpha_3\]
and
\[\Sq^2(\alpha_i\beta_3) = \alpha_i\beta_3^2;\; \Sq^2(\alpha_i^3) = \alpha_i^5.\]
We can thus compute the dual map $\Sq_2 \colon H_5(\pi;\Z/2) \to H_3(\pi;\Z/2)$.
We need to show that the element $(\alpha_1^2\alpha_2^2 \alpha_3)^* \in H_5(\pi;\Z/2)$ does not lie in the image of the reduction modulo two.

For the chain complex $C_*(\Z/2;\Z)$ we will again use $a_i^j$ to denote a generator of $C_j(\Z/2;\Z))$ in the $i$th copy, for $i=1,2$.
Recall that a $\Z[\Z/4]$-module resolution of $\Z$, tensored down over $\Z$, is given by
\[\Z(6) \xrightarrow{4} \Z(5) \xrightarrow{0} \Z(4) \xrightarrow{4}\Z(3) \xrightarrow{0}\Z(2) \xrightarrow{4}\Z(1) \xrightarrow{0} \Z(0),\]
where the parenthetical numbers indicate the grading.
We will use the notation $a_3^jb_3^k$ for the generator of $\Z(j+2k)$, an integral lift of the dual element to $\alpha_3^j\beta_3^k$. In $C_*(\pi;\Z) =  C_*(\Z/2;\Z) \otimes C_*(\Z/2;\Z) \otimes C_*(\Z/4;\Z)$, we compute:
\[\partial (a_1^2 \otimes a_2^2 \otimes a_3) = 2a_1 \otimes a_2^2 \otimes a_3 + 2 a_1^2 \otimes a_2 \otimes a_3.\]
So the obvious dual element to $\alpha_1^2\alpha_2^2 \alpha_3$ is not a cycle with $\Z$ coefficients.  To show that $\gamma$ does not lie in the image of $\Sq_2 \circ \red_2$ however, we need to argue that there is no way to add other elements of $C_5(\pi;\Z)$ to make $a_1^2 \otimes a_2^2 \otimes a_3$ into a cycle, in such a way that we preserve the correct image $\gamma$ of the reduction modulo two in $H_3(\pi;\Z/2)$.  We can try adding linear combinations of the chains
\[a_1^i \otimes a_2^j \otimes a_3;\; a_k \otimes b_3^2;\; a_s^3\otimes b_3;\; a_t^5 \text{ or } a_{\ell} \otimes a_m^2\otimes b_3,\]
for some $i$,$j$ with $i+j=4$, for some $k,s,t =1,2$, and for some nonempty $\{\ell,m\} \subseteq \{1,2\}$.
Then compute
\[\partial(a_1^3 \otimes a_2 \otimes a_3) = \partial(a_1 \otimes a_2^3 \otimes a_3) = \partial(a_1^5) = \partial(a_2^5)=0; \]
\[\partial(a_i^4 \otimes a_3) = 2a_i^3\otimes a_3;\; \partial(a_j^3 \otimes b_3) = -4a_j^3 \otimes a_3\]
for $i,j = 1,2$;
\[\partial(a_1 \otimes a_2^2 \otimes b_3) = -2a_1 \otimes a_2 \otimes b_3 -4a_1 \otimes a_2^2 \otimes a_3;\]
\[\partial(a_1^2 \otimes a_2 \otimes b_3) = 2a_1 \otimes a_2 \otimes b_3 -4a_1^2 \otimes a_2 \otimes a_3;\]
 and
\[\partial(a_k \otimes b_3^2) = -4a_k \otimes a_3 b_3.\]
Since $\alpha_i\beta_3^2$ and $\alpha_j^4\alpha_3$ are in the image of $\Sq^2$, for any $i,j=1,2$, any occurrence of their dual terms $a_i \otimes b_3^2$ or $a_j^4 \otimes a_3$ in a putative linear combination must have even coefficient, or this occurrence would alter the image under $\Sq_2$, causing it to deviate from being $\gamma = (\alpha_1\alpha_2\alpha_3)^* = a_1 \otimes a_2\otimes a_3$. Therefore the boundary of every term in our linear combination must have a term with coefficient divisible by $4$.  On the other hand the boundary we are trying to cancel is \[\partial (a_1^2 \otimes a_2^2 \otimes a_3) = 2a_1 \otimes a_2^2 \otimes a_3 + 2 a_1^2 \otimes a_2 \otimes a_3,\]
in which both terms are only divisible by $2$.  We see that there is no way to cancel the terms $-4a_1 \otimes a_2^2 \otimes b_3$ and
$-4a_1^2 \otimes a_2 \otimes b_3$ while still having the necessary odd coefficient of $a_1^2 \otimes a_2^2 \otimes a_3$.  It follows that $\gamma = a_1 \otimes a_2 \otimes a_3$ does not lie in the image of $\Sq_2 \circ \red_2 \colon H_5(\pi;\Z) \to H_3(\pi;\Z/2)$ as claimed.  For interest, we remark that $2a_1^2 \otimes a_2^2 \otimes a_3 + a_1\otimes a_2^2 \otimes b_3 + a_1^2 \otimes a_2 \otimes b_3$ represents a homology class in $H_5(\pi;\Z)$, but of course this does not map to $\gamma$. This completes the proof of the claim.\qed

We write $\pi:=\Z/2\times \Z/2\times \Z/4=\langle a,b,c\mid [a,b],[b,c],[c,a],a^2,b^2,c^4\rangle$.

\begin{claim}
The form $A(\gamma)$ is not even.
\end{claim}

\begin{align*}
H^2(K;\Z\pi)\cong (\Z\pi)^6/\langle & (1-b,c^{-1}-1,0,1+a,0,0),(a-1,0,1-c^{-1},0,1+b,0), \\ &(0,1-b,a-1,0,0,1+c+c^2+c^3)\rangle
\end{align*}
and the form $A(\gamma)$ is given by
\[\begin{pmatrix}
(1-c)(1-c^{-1})&(1-c^{-1})(1-b)&(1-c^{-1})(1-a)&0&0&0\\
(1-b)(1-c)&2(1-b)&(1-a)(1-b)&0&0&0\\
(1-a)(1-c)&(1-a)(1-b)&2(1-a)&0&0&0\\
0&0&0&0&0&0\\
0&0&0&0&0&0\\
0&0&0&0&0&0\end{pmatrix}\]

Assume there exists $q$ with $q+q^*=A(\gamma)$. Then $q$ has to be of the form
\[\begin{pmatrix}
(1-c)+z&(1-b)+x&*\\
-c(1-b)-\ol x&(1-b)+y&*\\
*&*&*
\end{pmatrix}\]
with $z+\ol z=y+\ol y=0$. The stars mean that we only specify the four entries in the upper left. The proof will show that this matrix cannot represent any  form on $H^2(K;\Z\pi)$, no matter how we fill in the rest of the matrix.

Now work modulo $2$ and under the projection $a=1$. To represent a form on $H^2(K;\Z\pi)$, the following equations have to be satisfied:
\[0=(1-c+z)(1-b)+(1-b+x)(c-1)=(1-b)z+(c-1)x\]
\[0=(-c(1-b)-\ol x)(1-b)+(1-b+y)(c-1).\]
Since $z+\ol z=0$, there are $z_1,z_2\in\Z$ with $z=(z_1+z_2b)(c-c^3)$. Thus, for $k:=z_1-z_2$, we have $(1-b)z=(k-kb)(c-c^3)$. It follows from the first equation that $x=(k-kb)(c+c^2)+x'N_c$ for some element $x'\in \Z\pi$. Since $(1-b)(1-b)=2(1-b)$ and we are working modulo 2, if we replace $x$ in the second equation above by $(k-kb)(c+c^2)+x'N_c$, we see that
\[0=-\ol{x'}N_c(1-b)+(1-b+y)(c-1).\]
Multiply by $(1+c)$, to obtain
\[0=(1-b+y)(1+c^2).\]
Evaluate at the neutral group element to yield
\[0=1+y_0+y_{c^2}=1,\]
since $y_0=y_{c^2}=0$. This is a contradiction, and it follows that $A(\gamma)$ cannot be even as claimed.\qed

Therefore, the \cref{secondary-property} holds for $\pi = \Z/2 \times \Z/2 \times \Z/4$.
\medskip

\noindent\textit{Any $3$-generator abelian group.}\\ \smallskip
Now we deduce the Secondary Property for $G$ one of the groups on the list from the start of this subsection.
For $G\cong G_1\times G_2\times G_3$, with $G_i$ cyclic and $|G_3|\geq 4$, the element $\gamma' \in H_3(G;\Z/2)$ coming from the generator of $$H_1(G_1;\Z/2)\otimes H_1(G_2;\Z/2)\otimes H_1(G_3;\Z/2)\cong \Z/2$$ is also not in the image of $\Sq_2\circ\red_2$. This can be seen by considering the projection $\phi \colon G\to \Z/2\times \Z/2\times \Z/4 = \pi$, under which $\gamma'$ is mapped to $\gamma$.
By \cref{lem:commuting-diagram-A-phi}, we have a commutative diagram
\[\xymatrix{
H_3(G;\Z/2) \ar[d]^A\ar[r]^-{\phi_*} & H_3(\pi;\Z/2) \ar[d]^A\\
\widehat{H}^0(\Sesq(H^2(K_G;\Z G))) \ar[r]^-{\phi_*}& \widehat{H}^0(\Sesq(H^2(K_{\pi};\Z\pi))}\]
Recall that $\gamma' \in H_3(G;\Z/2)$ is such that $\phi_*(\gamma') = \gamma$, and we just showed that $A(\gamma) \neq 0$.
It follows that $A(\gamma') \neq 0$.  Therefore the \cref{secondary-property} also holds for~$G$.  This completes the proof of \cref{lem:ab-groups-3-gens-sec-ter}.\qed

\subsection{Abelian groups with any finite number of generators}

\cref{thm:conjectures-abelian} now follows by combining the work we have done above.

\begin{proof}[Proof of \cref{thm:conjectures-abelian}]
The \cref{secondary-property} holds for all finitely generated abelian groups,
since by the previous subsections it holds for all abelian groups with at most three generators, and \cref{cor:secabelian} tells us that this suffices to prove the conjecture for all finitely generated abelian groups.

The \cref{tertiary-property} holds for all finitely generated abelian groups, since by the previous subsections it holds for all abelian groups with at most two generators, and \cref{cor:terabelian} applies to show that the conjecture holds for all abelian groups.
\end{proof}

\noindent We also have the following result on the Tertiary Property.

\begin{cor}
	If a finitely presented group~$G$ has the \cref{tertiary-property}, then so does $G\times \Z$.
\end{cor}

\begin{proof}
Let $G^{ab}$ be the abelianisation of $G$ and let $G^{ab} \cong \bigoplus_{i=1}^m C_i$ be a decomposition into cyclic groups~$C_i$.
Consider the projection $p_0\colon G \times \Z \to G =: G_0$ and the surjections \[p_i\colon G\times \Z\to G^{ab}\times \Z \to C_i \times \Z =:G_i,\] for $i=1,\dots,m$.
To apply \cref{thm:tersurj}, we need to see that the induced map
	\[\prod_{i=0}^m(p_i)_*\colon H_2(G \times \Z;\Z/2)/\im (d^2_{4,1}, d^3_{5,0}) \to \prod_{i=0}^m H_2(G_i;\Z/2)/\im (d^2_{4,1}, d^3_{5,0})\]
	is an injection.
This is obvious for elements of $H_2(G \times \Z;\Z/2)$ coming from $H_2(G;\Z/2)$, by considering the image under $(p_0)_*$. For elements of $H_2(G \times \Z;\Z/2)$ coming from $H_1(G;\Z/2)\otimes_{\Z} H_1(\Z;\Z/2)$ in the K\"{u}nneth theorem, we consider the appropriate map $(p_i)_*$, $i=1,\dots,m$, and use that for every cyclic group $C$ (of even or infinite order, so that $H_1(C;\Z/2)$ is nontrivial), the nontrivial element in $H_2(C \times \Z;\Z/2)$ coming from $H_1(C;\Z/2)\otimes_{\Z} H_1(\Z;\Z/2)$  does not lie in the image of the boundary maps $d^2_{4,1}, d^3_{5,0}$, as we showed in \cref{sec:twogenter}.
Thus the corollary indeed follows from \cref{thm:tersurj}.
\end{proof}

\begin{cor}
Let $\pi$ be a finite group whose 2-Sylow subgroup is abelian or has periodic cohomology. Then the \cref{secondary-property} and the \cref{tertiary-property} hold for $\pi$. In particular, both properties hold for every finite group with periodic cohomology.
\end{cor}

\begin{proof}
First, we can reduce the verification of the properties to $2$-Sylow subgroups, by \cref{lem:2sylow-sec,lem:2sylow-ter}.  By \cite[Theorem VI 9.3]{brown}, every finite $2$-group with periodic cohomology is either abelian or generalised quaternion.
The corollary therefore follows from the statement for abelian groups, \cref{thm:conjectures-abelian}, together with \cref{theorem:quarterion-groups} on the generalised quaternion groups.

By \cite[Theorem VI 9.5]{brown}, for every finite group with periodic cohomology its 2-Sylow subgroup also has periodic cohomology. This implies the last statement.
\end{proof}

\chapter{The stable classification for fundamental group \texorpdfstring{$\Z\times \Z/2$}{ZxZ/2}}\label{section:stable-classn-ZxZ2}

We apply our results to the case of smooth, closed, connected, spin $4$-manifolds with fundamental group $\Z\times\Z/2$, and prove \cref{thm:ZxZ2-intro}.
We will recall the statement later, after some preliminary lemmas.


		For the following, let $\pi:=\Z\times\Z/2=\langle T,t\mid [T,t],T^2\rangle$.
Let $I:=I\pi$ denote the augmentation ideal of $\Z\pi$, generated by $1-t$ and $1-T$.
We start by collecting some basic results about group homology and cohomology that we shall need.

\begin{lemma}\label{lem:group-homology-ZxZ2}
	Let $\pi := \Z \times \Z/2 = \langle t,T \mid [t,T], T^2 \rangle$.
	We have
	\[H_i(\pi;\Z/2) \cong \begin{cases} \Z/2 & i=0 \\ \Z/2 \oplus \Z/2 & i \geq 1 \end{cases}\]
	\[H_i(\pi;\Z) \cong \begin{cases} \Z & i=0 \\ \Z \oplus \Z/2 & i = 1 \\ \Z/2 & i \geq 2. \end{cases}\]
	The cohomology ring with $\Z/2$ coefficients is
	\[H^*(\pi;\Z/2) \cong \Z/2[t,T]/(t^2),\]
	where we abuse notation and use $t,T$ for elements of $H^1(\pi;\Z/2)$ that are dual to the corresponding generators of $\pi$.
	Moreover, we have the Steenrod operations:
	\[\xymatrix @C+0.2cm {H^3(\pi;\Z/2) \ar[d]_{\cong} \ar[r]^{\Sq^2} & H^5(\pi;\Z/2) \ar[d]_{\cong} \\
		\Z/2\langle T^3 \rangle \oplus \Z/2 \langle T^2t\rangle \ar[r]^{\begin{pmatrix}1 & 0 \\ 0 & 1 \end{pmatrix}} & \Z/2\langle T^5 \rangle \oplus \Z/2 \langle T^4t\rangle}\]
	\[\xymatrix @C+0.2cm {H^2(\pi;\Z/2) \ar[d]_{\cong} \ar[r]^{\Sq^2} & H^4(\pi;\Z/2) \ar[d]_{\cong} \\
		\Z/2\langle T^2 \rangle \oplus \Z/2 \langle Tt\rangle \ar[r]^{\begin{pmatrix}1 & 0 \\ 0 & 0 \end{pmatrix}} & \Z/2\langle T^4 \rangle \oplus \Z/2 \langle T^3t\rangle.}\]
	%
	%
	The reduction modulo 2 maps are given by:
	\[\operatorname{red}_2 \colon H_5(\pi;\Z) \cong \Z/2 \xrightarrow{(1,0)} H_5(\pi;\Z/2) \cong \Z/2\langle (T^5)^* \rangle \oplus \Z/2 \langle (T^4t)^*\rangle,\]
	\[\operatorname{red}_2 \colon H_4(\pi;\Z) \cong \Z/2 \xrightarrow{(0,1)} H_4(\pi;\Z/2) \cong \Z/2\langle (T^4)^* \rangle \oplus \Z/2 \langle (T^3t)^* \rangle.\]
	It follows that \[\Sq_2 \circ \red_2 \colon H_5(\pi;\Z) \cong \Z/2 \to H_3(\pi;\Z/2) \cong \Z/2 \oplus \Z/2\]
	is injective, while
	\[\Sq_2 \circ \red_2 \colon H_4(\pi;\Z) \cong \Z/2 \to H_2(\pi;\Z/2) \cong \Z/2 \oplus \Z/2\]
	is the zero map.
\end{lemma}

\begin{proof}
	The computations of the homology and the cohomology ring with $\Z/2$-coefficients follow from the K\"unneth theorem.
	The Steenrod square computations then follow from the axioms, as in \cref{sec:abelian-groups}.
\end{proof}

We can now calculate the bordism group $\Omega_4^{\spin}(B\pi)$.
\begin{lemma}\label{lem:Omega-4-ZxZ2}
	We have an isomorphism  $\Omega_4^{\spin}(B\pi)\cong 16\Z \oplus \Z/8$ given by the signature and the Brown invariant sending $f\colon M\to S^1\times\RP^3\subseteq B\pi$ to $f^{-1}(\{\pt\}\times \RP^2)\in\Omega_2^{\mathrm{Pin}^-}\cong \Z/8$.
\end{lemma}

\begin{proof}
	Recall that $\Omega_4^{\spin}\cong 16\Z$, where the isomorphism is given by the signature. As usual, there is a splitting $\Omega^{\spin}_4(X)\cong \wt \Omega_4^{\spin}(X)\oplus 16\Z$, as explained in \cref{prop:spin-AHSS}.  	
	So it suffices to show that $\wt \Omega_4^{\spin}(B(\Z\times\Z/2)) \cong \Z/8$.
	
	We have the following chain of isomorphisms
	\[\wt \Omega_4^{\spin}(B\pi)\cong \wt \Omega_4^{\spin}(B\Z/2) \oplus \Omega_3^{\spin}(B\Z/2) \cong  \Omega_3^{\spin}(B\Z/2)\cong 
	\Omega_2^{\pin^-}\cong \Z/8.\]
	The second isomorphism uses that $\wt \Omega_4^{\spin}(B\Z/2)=0$, as the calculations in \cref{sec:cyclic-gps} showed.
	The fourth isomorphism $\Omega_2^{\pin^-}\cong \Z/8$ was computed by Kirby and Taylor \cite{KT90}.
	Recall that $B\pi \simeq S^1 \times \RP^\infty$. The composite of the first two isomorphisms is given on $[M\xrightarrow{f} S^1\times \RP^\infty]$ by making $f$ transverse to $\{s_0\}\times\RP^\infty$ and then sending $[M,f]$ to $[f^{-1}(\{s_0\}\times\RP^\infty)\to \RP^\infty]$. The third isomorphism, the Smith map, is given by homotoping the map to $\RP^{\infty}$ to the 3-skeleton, making $f$ transverse to $\RP^2$, and then sending $[M\xrightarrow{f} \RP^3]$ to $[f^{-1}(\RP^2)]$. The spin structure on $M$ induces a $\pin^-$ structure on $f^{-1}(\RP^2)$.  Generators of the bordism groups are represented by $S^1 \times \RP^3 \in \wt \Omega_4^{\spin}(B(\Z\times\Z/2))$, $\RP^3 \in \Omega_3^{\spin}(B\Z/2)$ and $\RP^2 \in \Omega_2^{\pin^-}$.
	
	We argue that the third map $\Omega_3^{\spin}(B\Z/2)\cong \Omega_2^{\pin^-}$ is indeed an isomorphism, as we were unable to find an explicit reference. As described above~\cite[pp.~212--3]{KT90}, the map arises by representing a class $[M,c] \in \Omega_3^{\spin}(B\Z/2)$ by a map $c \colon M \to \RP^3$, taking the inverse image $c^{-1}(\RP^2)$.  A standard transversality argument shows that this is well-defined, and a homomorphism.  The map is surjective since $\RP^2$ with one of its $\pin^-$ structures represents a generator of $\Omega_2^{\pin^-}\cong \Z/8$~\cite[Section~3]{KT90}, and the map applied to $\id \colon \RP^3 \to \RP^3$ yields $\RP^2$.  Therefore $\Omega_3^{\spin}(B\Z/2)$  surjects to $\Omega_2^{\pin^-}\cong\Z/8$. 	
	We will now show that $\Omega^{\spin}_3(B\Z/2)\cong \wt{\Omega}^{\spin}_4(B(\Z\times\Z/2))$ has at most eight elements and thus this surjection is an isomorphism.
	Write $\pi\cong \Z\times \Z/2$. By \cref{prop:spin-AHSS} there is a filtration $F_2\subseteq F_1\subseteq \wt\Omega_4^{\spin}(B\pi)$, arising from the Atiyah-Hirzebruch spectral sequence, with quotients:
	\begin{enumerate}
		\item $\wt\Omega_4^{\spin}(B\pi)/F_1$, a subgroup of $H_4(\pi;\Z)\cong \Z/2$;
		\item $F_1/F_2$, isomorphic to the cokernel of $\Sq_2\circ \red_2\colon H_5(\pi;\Z)\to H_3(\pi;\Z/2)$, which is $\Z/2$ using the computation from \cref{lem:group-homology-ZxZ2};
		\item $F_2$ a quotient of the cokernel of $\Sq_2\colon H_4(\pi;\Z/2)\to H_2(\pi;\Z/2)$. The latter is isomorphic to $\Z/2$ using the computation from \cref{lem:group-homology-ZxZ2}.
	\end{enumerate}
	So $\wt\Omega_4^{\spin}(B\pi)$ has at most eight elements as claimed.
	It follows that in fact these three copies of $\Z/2$ survive to the $E^\infty$ page and assemble to give $\wt \Omega_4^{\spin}(B(\Z\times\Z/2))\cong \Z/8$, i.e.\ the filtration is $0 \subseteq \Z/2 \subseteq \Z/4 \subseteq \Z/8$.
\end{proof}

We remark that $(n,m) \in \Z \oplus \Z/8 \cong \Omega_4^{\spin}(B(\Z\times\Z/2))$ is realised by the $4$-manifold $\#^{n} K3 \# \big((\#^m \RP^3) \times S^1 \big)$.
We also make the following observations for later use:
\begin{enumerate}[(i)]
  \item $\pri(M) = 0$ if and only if $\beta(M)$ is even;
  \item if $\pri(M)=0$, then $\msec(M) =0$ if and only if $\beta(M) \in \{0,4\}$.
\end{enumerate}

\begin{lemma}
	\label{lem:action}	
	Changing the spin structure by $T\in H^1(\pi;\Z/2)$ acts on $\wt \Omega_4^{\spin}(B\pi)\cong \Z/8$ by multiplication with $-1$.
\end{lemma}

\begin{proof}
	The isomorphism $\wt\Omega^{\spin}_4(B\pi)\cong \Omega_3^{\spin}(B\Z/2)$ is equivariant with respect to changing the spin structure by $T$. Pick a fixed orientation and spin structure for $\RP^3$. Since $\RP^3$ generates $\Omega_3^{\spin}(B\Z/2)$, it suffices to see that changing the spin structure of $\RP^3$ yields $-[\RP^3] \in\Omega_3^{\spin}(B\Z/2)$. Use that $\RP^3$ admits an orientation reversing diffeomorphism to see that  $-[\RP^3]$ equals either $[\RP^3]$ or $[T^*\RP^3]$, where $T^*\RP^3$ denotes $\RP^3$ with the other spin structure. But $[\RP^3]$ is a generator of $\Omega_3^{\spin}(B\Z/2)\cong \Z/8$, so $[\RP^3]\neq -[\RP^3]$. It follows that $-[\RP^3]=[T^*\RP^3]$ as claimed.
\end{proof}

Using \cref{lem:Omega-4-ZxZ2,lem:action}, Kreck's modified surgery implies that smooth, closed, connected, spin $4$-manifolds with fundamental group $\pi = \Z\times\Z/2$ are classified by their signature and the Brown invariant in $(\Z/8)/\pm1$. As promised in the introduction,  we now apply our classification results to obtain an intrinsic description of the Brown invariant in this case.   For convenience we repeat the statement of \cref{thm:ZxZ2-intro}. 

\begin{thm}
	\label{thm:ZxZ2}
	For a smooth, closed, connected, spin $4$-manifold with fundamental group $\pi$, the Brown invariant $\beta(M)\in \Z/8$ is detected as follows:
	\begin{enumerate}
		\item\label{item:ZxZ2-1} $\beta(M)=\pm 1$ if and only if $\pi_2(M)$ is stably free and $\lambda_M$ is stably isometric to a form induced up to $\Z\pi$ from a form over $\Z$;
		\item\label{item:ZxZ2-2} $\beta(M)=\pm 3$ if and only if $\pi_2(M)$ is stably free but $\lambda_M$ is not stably induced from a form over $\Z$;
		\item\label{item:ZxZ2-3} $\beta(M)=\pm 2$ if and only if $\pi_2(M)$ is not stably free and $\lambda_M$ is odd;
		\item\label{item:ZxZ2-4} $\beta(M)\in\{0,4\}$ if and only if $\pi_2(M)$ is not stably free and $\lambda_M$ is even.
	\end{enumerate}
	In the last case, $\pi_2(M)$ is stably isomorphic to $I\oplus I^*$, where $I:=I\pi$ is the augmentation ideal, and there exists such a stable isomorphism for which the restriction of $\lambda_M$ to $I^*$ is trivial. Let $x\in\pi_2(M)$ be the image of $\phi\in I^*$ under such an isomorphism, where $\phi\in I^*$ is determined by $\phi(1-t)=1+T$ and $\phi(1-T)=0$. Then:
	\begin{enumerate}
		\setcounter{enumi}{4}
		\item\label{item:ZxZ2-5} $\beta(M)=0$ if and only if $\tau_M(x)=0$, where $\tau_M$ is the Kervaire--Milnor invariant.
	\end{enumerate}
\end{thm}
The rest of this section consists of the proof of this theorem.  We will assume that all 4-manifolds are smooth, closed, connected, and spin without further comment.
First we show that the stable isomorphism class of $\pi_2(M)$ determines the parity of $\beta(M)$.
Since the map from $\Omega_4^{\spin}(B\pi)\to H_4(\pi;\Z/2)\cong\Z/2$ is given by reducing the Brown invariant modulo two, this is implied by the following lemma.

\begin{lemma}\label{lem:stably-free-if-cM-nonzero}
	Let $M$ be a spin $4$-manifold with $\pi_1(M)=\Z\times\Z/2$. The image $c_*[M] \in H_4(\pi;\Z) \cong \Z/2$ of the fundamental class is nontrivial if and only if $\pi_2(M)$ is stably free as a $\Z\pi$-module.
\end{lemma}


\begin{proof}
	By \cite{KPT18}, the image $c_*[M]\in H_4(B(\Z\times\Z/2);\Z) \cong \Z/2$ of the fundamental class of $M$, stably determines $\pi_2(M)$ as a $\Z\pi$-module. Since $\pi_2(S^1\times \RP^3)=0$, this implies that $\pi_2(M)$ is stably free if $c_*[M]$ is nontrivial. Let $(C_*,d_*)$ be a resolution of $\Z$ as a $\Z\pi$-module by finitely generated free modules. If $c_*[M]$ is trivial, then $\pi_2(M)$ is stably isomorphic to $\ker d_2\oplus \coker d^2$, which is not stably free. To see this, note that if $\ker d_2$ were projective, then there would be a finite length projective resolution of $\Z$ as a $\Z \times \Z/2$-module, which cannot happen since $\Z[\Z \times \Z/2]$ has homology in all nonnegative degrees.
\end{proof}

We have shown that $\pri(M)$, i.e.\ whether $\beta(M)\in\{\pm 1,\pm 3\}$ or $\beta(M)\in \{0,\pm 2,4\}$, is detected by whether or not $\pi_2(M)$ is stably free.

\section{The case \texorpdfstring{$\beta(M)\in\{\pm1,\pm3\}$}{of odd Brown invariant}}

This section distinguishes between the cases \eqref{item:ZxZ2-1} and \eqref{item:ZxZ2-2} of \cref{thm:ZxZ2}.
A generator of $\wt\Omega_4^{\spin}(B\pi)\cong\Z/8$ is given by $S^1\times\RP^3$, which has trivial $\pi_2(M)$.  Applying \cite[Theorem C]{surgeryandduality} (\cref{thm:stablediffeoclasses}), it follows that $\beta(M)=\pm 1$ implies that $\lambda_M$ is stably isometric to a form over $\Z\pi$ on a f.g.\ free module induced from a form over $\Z$.

It therefore remains to show that, for a fixed signature, the stable isometry class of the equivariant intersection form distinguishes between the two stable diffeomorphism classes with nontrivial image of the fundamental class.
Since the signature determines the stable isometry class of forms over $\Z$, we deduce that if $\lambda_M$ is stably isometric to a form induced from $\Z$, then $\beta(M)=\pm 1$.

Let $M_1$ and $M_2$ be two spin $4$-manifolds with nontrivial $c_*[M_1]=c_*[M_2]\in H_4(\pi;\Z)\cong \Z/2$. Fix a spin structure on each of $M_1$ and $M_2$.

\begin{lemma}
 After changing the spin structure on $M_1$ if necessary, we can assume that the secondary invariant of the 1-skeleton sum $M_1\#_1 M_2$ in $H_3(\pi;\Z/2)/d_2$ is trivial.
\end{lemma}

\begin{proof}
The signature is unaffected by the change in spin structure and changing the signature by adding a simply-connected manifold does not change the secondary invariant. Hence we can consider the image of $M_1\#_1 M_2$ in reduced bordism $\wt{\Omega}_4^{\spin}(B\pi)$.

By \cref{lem:action}, we can change the Brown invariant of $M_1$ by a sign by changing the spin structure on $M_1$.  Since $\pri(M_i) \neq 0$, we know that $\beta(M_i) \in \{\pm 1, \pm 3\}$ for $i=1,2$. It follows that $\beta(M_1 \#_1 M_2) = \beta(M_1) + \beta(M_2)$ is even, so $\pri(M_1 \#_1 M_2)=0$. But in addition, for each $x \in \{1,3\}$, and for every possible $\beta(M_2) \in \{\pm 1,\pm 3\}$, there is a choice of sign $\pm$ such that $\pm x + \beta(M_2) \in \{0,4\}$. The lemma follows.
\end{proof}

Next we will show that the tertiary invariant of $M_1\#_1M_2$ in $H_2(\pi;\Z/2)/d_2\cong \Z/2$ is determined by $[\lambda_{M_1}], [\lambda_{M_2}]\in L_4(\Z\pi)$.
Let $\mathbb{L}\langle 1\rangle$ denote the 1-connective cover of the quadratic $L$-theory spectrum $\mathbb{L}$.  Recall that
\[\mathbb{L}\langle 1\rangle_4(X) \cong H_4(X;\mathbb{L}\langle 1 \rangle_{\bullet}) \cong \mathcal{N}_{\operatorname{TOP}}(X) \]
for a 4-dimensional Poincar\'{e} complex~$X$.

\begin{lemma}\label{lem:assembly}
	The assembly map $\alpha\colon \mathbb{L}\langle 1\rangle_4(S^1 \times \RP^3)\to L_4(\Z\pi)$ has domain isomorphic to $\Z\oplus \Z/2\oplus\Z/2$, codomain isomorphic to $\Z \oplus \Z \oplus \Z/2$, and the torsion element $(0,0,1)$ in the codomain lies in the image of $\alpha$.
\end{lemma}

\begin{proof}
	From a spectral sequence computation using the Atiyah-Hirzebruch sequence $H_p(S^1 \times \RP^3;\mathbb{L}\langle 1\rangle_q) \Rightarrow \mathbb{L}\langle 1\rangle_{p+q}(S^1 \times \RP^3)$, and noting that $\mathbb{L}\langle 1\rangle_q =0$ for $q \leq 1$ and $q=3$, $\mathbb{L}\langle 1\rangle_4 \cong \Z$ and $\mathbb{L}\langle 1\rangle_2 \cong \Z/2$,  it follows from $H_0(S^1 \times \RP^3;\Z) \cong \Z$ and $H_2(S^1 \times \RP^3;\Z/2) \cong \Z/2 \oplus \Z/2$ that \[\mathbb{L}\langle 1\rangle_4(S^1 \times \RP^3)\cong \Z\oplus \Z/2\oplus\Z/2\]
since there are no nontrivial differentials going in or out of the 4-line. The only potentially nonzero map is $H_3(S^1 \times \RP^3;\Z/2) \to H_0(S^1 \times \RP^3;\Z) \cong \Z$, but there can be no nontrivial homomorphism since the domain is torsion.
	
Since $\Wh(\Z \times \Z/2)= 0 = \Wh(\Z/2)$ by \cite[Theorem~3.2~(d-iii)~ and~(e)]{Luck-Stamm}  we have that
	\[L_4^h(\Z[\Z \times \Z/2]) \cong L_4^s(\Z[\Z \times \Z/2]) \text{ and } L^h_k(\Z[\Z/2]) \cong L^s_k(\Z[\Z/2])\]
	by \cite[Proposition~4.1]{Shaneson}, for $k \geq 0$.
	Therefore combining this with Shaneson splitting~\cite[Theorem~5.1]{Shaneson} we have
	\[L_4^h(\Z[\Z \times \Z/2]) \cong L_4^s(\Z[\Z \times \Z/2]) \cong L_4^s(\Z[\Z/2]) \oplus L_3^h(\Z[\Z/2]).\]
	Now since $L^h_4(\Z[\Z/2]) \cong L^s_4(\Z[\Z/2])$, by~\cite[Chapter~13A.1]{Wall} we have
	\[L_4^s(\Z[\Z/2]) \oplus L_3^h(\Z[\Z/2]) \cong L_4^h(\Z[\Z/2]) \oplus L_3^h(\Z[\Z/2]) \cong (\Z \oplus \Z) \oplus \Z/2.\]
	
	We have therefore identified the assembly map $\alpha\colon \mathbb{L}\langle 1\rangle_4(S^1 \times \RP^3)\to L_4(\Z\pi)$ with a map $\Z \oplus \Z/2 \oplus \Z/2 \to \Z \oplus \Z \oplus \Z/2$.   We show that $(0,0,1)$ lies in the image.   We start with the commuting diagram which says that Shaneson splitting of quadratic $L$ groups arises from taking the product of a degree one normal map with the given surgery obstruction with $S^1$ (see also~\cite{Ranicki-splittings}):
	\[\xymatrix{\mathbb{L}\langle 1 \rangle_3(\RP^3) \ar[r]^-{\times S^1} \ar[d]^{\alpha_{\RP^3}} & \mathbb{L}\langle 1 \rangle_3(S^1 \times \RP^3) \ar[d]^{\alpha_{S^1 \times \RP^3}} \\
		L_3(\Z[\Z/2]) \ar@{^{(}->}[r] \ar[d]^{\cong} & L_4(\Z[\Z \times \Z/2]) \ar[d]^{\cong}  \\
		\Z/2 \ar@{^{(}->}[r] & \Z \oplus \Z \oplus \Z/2 .}\]
	In order to show that $(0,0,1)$ lies in the image of the right hand surgery obstruction map $\alpha_{S^1 \times \RP^3}$, it suffices to show that the left hand $\alpha_{\RP^3}$  is surjective. To see this we consider the following commutative diagram of $\mathbb{L}$-homology groups and assembly maps. The key input is that the vertical map $L_3(\Z[\Z]) \to L_3(\Z[\Z/2])$ induced by the group homomorphism $\Z \to \Z/2$ is an isomorphism by \cite[Lemma~13A.9]{Wall}.  We also use that the assembly map  $\alpha \colon \mathbb{L}_3(S^1) \xrightarrow{} L_3(\Z[\Z])$ is an isomorphism, which is shown in~\cite[p.~126, Proof of Lemma~6.2(1)]{Connolly-Davis-Khan-rigidity}.	
\[\xymatrix @C-0.3cm {\Z/2 \ar[r]^-{\cong} \ar[d]^{=} & \mathbb{L}\langle 1 \rangle_3(S^1) \ar[rr]^-{\cong}  \ar[d]^{\cong} && \mathbb{L}_3(S^1) \ar[r]_-{\alpha}^-{\cong} \ar[d] & L_3(\Z[\Z]) \ar[r]^-{\cong} \ar[d]^{\cong} & \Z/2 \ar[d]^{=} \\
		\Z/2 \ar[r]^-{\cong}  & \mathbb{L}\langle 1 \rangle_3(\RP^3) \ar[r] & \mathbb{L}\langle 1 \rangle_3(\RP^\infty) \ar[r] & \mathbb{L}_3(\RP^{\infty}) \ar[r]^-{\alpha}  & L_3(\Z[\Z/2]) \ar[r]^-{\cong}  & \Z/2
	}\]
	It now follows from commutativity that $\alpha_{\RP^3} \colon \mathbb{L}\langle 1 \rangle_3(\RP^3) \to L_3(\Z[\Z/2])$, which in the diagram we factored through $\mathbb{L}\langle 1 \rangle_3(\RP^\infty) \to \mathbb{L}_3(\RP^{\infty})$,  is an isomorphism and so is in particular surjective.
\end{proof}

This assembly map can be identified with the surgery obstruction map and we can interpret elements of $\mathbb{L}\langle 1\rangle_4(S^1 \times \RP^3)$ as degree 1 normal maps over $S^1 \times \RP^3$.

\begin{lemma}\label{lem:spin-structures-ZxZ2}
	Suppose that $f \colon M\xrightarrow{} S^1 \times \RP^3$ is a degree 1 normal map. Then $M$ is spin.
\end{lemma}

\begin{proof}
	The existence of a degree 1 normal map implies that there is a stable vector bundle $\xi$ over $N:= S^1\times \RP^3$, which reduces to the  Spivak normal fibration $\eta$ of $N$. Moreover $\xi$ pulls back to the stable normal bundle $\nu_M$ of $M$. Let $\nu_N$ be the stable normal bundle of $N$, which also reduces to the Spivak normal fibration of $N$.    Since $M$ and $N$ are orientable, $w_1(M) = 0 =w_1(N)$. By the Whitney sum formula, we therefore have $w_2(\nu_N) = w_2(N)$ and  $w_2(\nu_M) = w_2(M)$.  Since $N$ is spin $w_2(N)=0$.
 Since there is an analogue of the Thom isomorphism theorem, Stiefel-Whitney class are defined for spherical fibrations, and we have that $w_2(\xi) = w_2(\eta) = w_2(\nu_N)=0$. Thus $0 = f^*(0) =  f^*(w_2(\nu_N)) = w_2(\nu_M) = w_2(M)$. Therefore $M$ is spin as asserted.
\end{proof}

By Shapiro's lemma we have $H^2(\pi;\Z\pi)=H^2(\Z;\Z[\Z])=0$ and $H^3(\pi;\Z\pi)=H^3(\Z;\Z[\Z])=0$.
It then follows from \cref{prop:int-form-singular} that $\lambda_M$ is nonsingular for closed $M$ with fundamental group $\pi$.

\begin{lemma}\label{lemma:stably-diffeo-iff-lambdas}
	Let $M_1$ and $M_2$ be spin $4$-manifolds with fundamental group $\Z \times \Z/2$ and $c_*[M_i]\neq 0$ for $i=1,2$. Then $M_1$ and $M_2$ are stably diffeomorphic if and only if $[\lambda_{M_1}]=[\lambda_{M_2}]\in L_4(\Z\pi)$.
\end{lemma}


\begin{proof}
	Let $\Omega^{\spin}_{4,\deg 1}(S^1 \times \RP^3)$ be the subset of $\Omega_4^{\spin}(S^1 \times \RP^3)$ of elements of degree 1. On both sets, $H^1(\pi;\Z/2)$ acts by changing the spin structure.
	
	Similarly to \cref{lem:spin-structures-ZxZ2}, a degree 1 normal bordism gives rise to a spin bordism. Hence we can define a map (of sets)
	\[u\colon \mathbb{L}\langle 1\rangle_4(S^1 \times \RP^3)\to \Omega^{\spin}_{4,\deg 1}(S^1 \times \RP^3)/H^1(\pi;\Z/2).\]
	Also let $\Omega^{\spin}_{4,\deg 1}(B\pi)$ denote the subset of $\Omega_4^\spin(B\pi)$ of elements $([M,f \colon M \to S^1 \times \RP^3])$ with $c_*[M]\neq 0$. Define a map
	\[\mathcal{O}\colon \Omega^{\spin}_{4,\deg 1}(B\pi)/H^1(\pi;\Z/2)\to L_4(\Z\pi)\]
	by sending $[M\xrightarrow{f} S^1 \times \RP^3]$, represented by a map $f$ that is an isomorphism on fundamental groups, to $(\pi_2(M),\lambda_M)\in L_4(\Z\pi)$. This map is defined, since $\pi_2(M)$ is stably free by~\cref{lem:stably-free-if-cM-nonzero}, $\lambda_M$ admits a quadratic refinement because $M$ is spin, and since $\lambda_M$ is nonsingular. The map $\mathcal{O}$ is well-defined because two $4$-manifolds with fundamental group $\pi$ representing the same spin bordism class are stably diffeomorphic and the action of $H^1(B\pi;\Z/2)$ does not change the intersection form.
	%
	Now consider the following diagram, obtained from restricting the above maps to reduced (generalised) homology.
	\[\xymatrix{
		\wt{\mathbb{L}}\langle 1\rangle_4(S^1 \times \RP^3)\ar[r]^-{u} \ar[d]^\alpha  & \wt{\Omega}^{\spin}_{4,\deg 1}(S^1 \times \RP^3)/H^1(\pi;\Z/2)  \ar[d]\\
		\wt{L}_4(\Z\pi) & \wt\Omega^{\spin}_{4,\deg 1}(B\pi)/H^1(B\pi;\Z/2) \ar[l]_-{\mathcal{O}}}
	\]
	The diagram commutes since $\pi_2(S^1 \times \RP^3)=0$ and thus the surgery obstruction of a degree one map $M\to S^1 \times \RP^3$ that is an isomorphism on fundamental groups is precisely $(\pi_2(M),\lambda_M)$.
	
	Since $\wt\Omega^{\spin}_{4,\deg 1}(B\pi)/H^1(B\pi;\Z/2)$ has two elements by \cref{lem:action}, the map $\mathcal{O}$ to $L_4(\Z\pi)$ is either trivial or not: we know that $S^1 \times \RP^3$ represents a bordism class mapping to $0 \in \wt{L}_4(\Z\pi)$ under $\mathcal{O}$, again since $\pi_2(S^1 \times \RP^3)=0$.
	By \cref{lem:assembly}, the assembly map $\alpha$ is nontrivial. Thus the map $\mathcal{O} \colon \wt\Omega^{\spin}_{4,\deg 1}(B\pi)/H^1(B\pi;\Z/2)\to \wt L_4(\Z\pi)$ is injective.
	It follows that spin manifolds $M_1$ and $M_2$ with the same signature, fundamental group $\pi=\Z\times\Z/2$, and $0\neq c_*[M_i]\in H_4(\pi;\Z)$ are stably diffeomorphic if and only if $[\lambda_{M_1}]=[\lambda_{M_2}]\in \wt{L}_4(\Z\pi)$.
Finally, note that the signature is detected by the image of the  intersection forms of $M_1$ and $M_2$ in the unreduced $L$-group $L_4(\Z\pi)$. It follows that spin $4$-manifolds $M_1$ and $M_2$ with fundamental group $\pi$ and $c_*[M_i] \neq 0$ are stably diffeomorphic if and only if $[\lambda_{M_1}]=[\lambda_{M_2}]\in L_4(\Z\pi)$.
This completes the proof of the lemma.
\end{proof}

We know that if $\pi_2(M)$ is stably free, i.e.\ $\pri(M) \neq 0$, then $\beta(M) \in \{\pm 1 \}$ or $\beta(M) \in \{\pm 3\}$.

\begin{cor}
If $\pi_2(M)$ is stably free,  the stable isometry class of the equivariant intersection form is induced from a form over $\Z$ if and only if $\beta(M) \in \{\pm 1\}$.
\end{cor}

\begin{proof}
  If $\beta(M) \in \{\pm 1\}$ then there is a simply connected, spin 4-manifold $X$ such that $M$ is stably diffeomorphic to $(S^1 \times \RP^3) \# X$. Since $\pi_2(S^1 \times \RP^3)=0$, it follows that $(\pi_2(M),\lambda_M)$ is stably isometric to a form induced from $\Z$.

  On the other hand, suppose that $\lambda_M$ is induced from a form over $\Z$, $Q$ say.  Since $\lambda_M$ is nonsingular, $Q$ is unimodular, and since $M$ is spin it follows that $Q$ is even.  Let $X$ be a closed, smooth, simply connected, spin 4-manifold with intersection form stably isometric to $Q$. This exists because the signature of $Q$ equals the signature of $M$, which is 0 modulo 16, and therefore $Q$ can be stably realised as the intersection form of a smooth 4-manifold, without violating Rochlin's theorem.  We have that $[\lambda_M] = [\lambda_{(S^1 \times \RP^3) \# X}] \in L_4(\Z\pi)$, and so by \cref{lemma:stably-diffeo-iff-lambdas} we know that $M$ is stably diffeomorphic to $(S^1 \times \RP^3) \# X$. It follows that $\beta(M) \in \{\pm 1\}$ as desired.

\end{proof}

\section{The case \texorpdfstring{$\beta(M)\in\{0,\pm2,4\}$}{of even Brown invariant}}
This covers cases~\eqref{item:ZxZ2-3} to \eqref{item:ZxZ2-5} of \cref{thm:ZxZ2}.
The elements in $\wt \Omega_4^{\spin}(B\pi)\cong \Z/8$ with trivial image in $H_4(\pi;\Z)$ but nontrivial secondary invariant correspond to $\beta(M)=\pm 2$.

	Let $(C_*,d_*)$ be the 2-periodic free $\Z\pi$-module resolution of $\Z$:
	\[\xymatrix @R-0.7cm @C+1cm {\cdots \ar[r]^{1-T} & \Z\pi \ar[r]^{1+T} & \Z \pi \ar[r]^{1-T} & \Z\pi \\
		& \oplus & \oplus &  \\
		\cdots \ar[r]_{1+T} \ar[uur]_{1-t}^{d_3} & \Z\pi \ar[r]_{1-T} \ar[uur]_{1-t}^{d_2} & \Z\pi \ar[uur]_{1-t}^{d_1} &
	}\]

\begin{lemma}\label{lem:ker-coker-identified-I-I*}
We have that $\ker d_2 \cong I \pi =:I$ and $\coker d^2 \cong I^*$.
\end{lemma}

\begin{proof}
	We have $\ker d_2\cong \im d_3\cong \coker d_4 \cong \coker d_2\cong \im d_1\cong I\pi=:I$.
	By Shapiro's lemma we have $H^2(\pi;\Z\pi)=H^2(\Z;\Z[\Z])=0$ and $H^3(\pi;\Z\pi)=H^3(\Z;\Z[\Z])=0$. Then it follows from \cref{lem:dual-identification} that $\coker d^2\cong I^*$.
\end{proof}

Recall that $c_*[M]=0$ implies that the stable extension
\[0\to \ker d_2\to C^M_2\oplus \pi_2(M)\xrightarrow{(p_1,p_2)} \coker d^2\to 0\]
splits, and thus by \cref{lem:ker-coker-identified-I-I*} we see that $\pi_2(M)$ is stably isomorphic to $I\oplus I^*$. Let $\cH_I$ be the hyperbolic form on $I\oplus I^*$.

To identify the intersection form on the null-bordant element, we need the following lemma, which is well-known to the experts.
	\begin{lemma}\label{lem:general-int-form}
		Let $N$ be a 4-dimensional thickening of a $2$-complex $K$. Then the double $D(N) := N \cup_{\partial N} N$ of $N$ has \[\pi_2(D(N))\cong H_2(K;\Z\pi)\oplus H^2(K;\Z\pi)\]
and the following equivariant intersection form.
\begin{enumerate}[(i)]
  \item Restricted to $H^2(K;\Z\pi)$, $\lambda_{D(N)}$ is trivial.
  \item Restricted to $H_2(K;\Z\pi)$, $\lambda_{D(N)}$ is isometric to the equivariant intersection form $\lambda_N$ on $\pi_2(N)$.
  \item For $x\in H^2(K;\Z\pi),y\in H_2(K;\Z\pi)$ we have $\lambda_{D(N)}(x,y)=\langle x,y\rangle$.
\end{enumerate}
 	\end{lemma}
	\begin{proof}
		Let $j\colon N\to D(N)$ be the inclusion and let $i\colon D(N)\to N\times I$ be the inclusion given by identifying $D(N)$ with the boundary of $N\times I$. Then $i\circ j$ is a homotopy equivalence and both $i$ and $j$ are isomorphisms on fundamental groups. We have $H_3(N\times I;\Z\pi)=0$ and by Poincar\'e-Lefschetz duality $H_2(N\times I,D(N);\Z\pi)\cong H^3(N\times I;\Z\pi)=0$. Thus we have a short exact sequence from the long exact sequence of the pair.
		\[0\to H_3(N\times I,D(N);\Z\pi)\xrightarrow{\partial} H_2(D(N);\Z\pi)\xrightarrow{i_*} H_2(N\times I;\Z\pi)\to 0.\]
		Again by Poincar\'e-Lefschetz duality, the boundary map $\partial$ is identified with the map $i^*\colon H^2(N\times I;\Z\pi)\to H^2(D(N);\Z\pi)$. Identifying $H_2(N;\Z\pi) \cong H_2(N \times I;\Z\pi)$ via $i_*\circ j_*$, the map $j_* \colon H_2(N;\Z\pi) \to H_2(D(N);\Z\pi)$  provides a splitting of $i_*$. Thus we have an isomorphism
		\[H_2(N;\Z\pi)\oplus H^2(N\times I;\Z\pi)\xrightarrow{(j_*,PD\circ i^*)}H_2(D(N);\Z\pi)\cong \pi_2(N).\]
		The first part of the lemma then follows from the homotopy equivalence $N\simeq K$.
		
		Under the isomorphism $H^2(D(N);\Z\pi)\xrightarrow{PD;\cong}H_2(D(N);\Z\pi)\cong \pi_2(D(N))$, the equivariant intersection form pulls back to the form $\langle -\cup-,[D(N)]\rangle$. On the image of $i^*\colon H^2(N\times I;\Z\pi)\to H^2(D(N);\Z\pi)$ we have
		\[\langle i^*x\cup i^*y,[D(N)]\rangle=\langle x\cup y,i_*[D(N)]\rangle=\langle x\cup y,0\rangle=0.\]
		On the image of $i_*\colon H_2(N;\Z\pi)\to H_2(D(N);\Z\pi)$ we use the geometric computation of the intersection form to see that $\lambda_{D(N)}$ pulls back to $\lambda_{N}$ as claimed.
		
		Given $x\in H^2(K;\Z\pi),y\in H_2(K;\Z\pi)$ we have \[\lambda_{D(N)}(x,y)=\langle i^*(x)\cup PD^{-1}(j_*y),[M]\rangle=\langle i_*x,j_*y\rangle=\langle j^*i^*x,y\rangle.\]
		Since we used $j^*i^*$ for the identification of the cohomology of $N\times I$ with that of $N$, this completes the proof of the lemma.
	\end{proof}

\begin{lemma}\label{lem:null-bordant-ZxZ2}
	The equivariant intersection form on the double of the 2-skeleton of $S^1\times \RP^3$, representing the trivial element of $\Omega_4^{\spin}(B(\Z \times \Z/2))$, is given by the even form $(I\oplus I^*,\cH_I)$.
\end{lemma}

\begin{proof}
		Taking the standard CW-structure on $S^1\times \RP^3$, the 2-skeleton has $H_2((S^1\times \RP^3)^{(2)};\Z\pi)\cong I$ and $H^2((S^1\times \RP^3)^{(2)};\Z\pi)\cong I^*$ by the above calculations. We thicken the 2-skeleton to obtain a manifold $N$ with boundary. Consider immersed spheres $\alpha$ and $\beta$ in $N$ that intersect transversely. We can view $\alpha$ and $\beta$ as elements of $\pi_2(S^1\times \RP^3)$,  and then since  $\pi_2(S^1\times \RP^3)=0$, the intersection product must be trivial. It follows that the intersection form on $\pi_2(N)$ is trivial. The lemma now follows from \cref{lem:general-int-form}.
\end{proof}

Let $(s_1,s_2)\colon \coker d^2\to C_2\oplus \pi_2(M)$ be some splitting. Since the \cref{secondary-property} holds for all abelian groups by \cref{sec:abelian-groups}, it holds for $\Z\times\Z/2$. Therefore $\msec(M)=0$ if and only if the equivariant intersection form $\lambda_M$ restricted along $s_2$ is even. We have already seen that on the stably diffeomorphism class given by $\msec=\ter=0$, corresponding to the null-bordant element considered in \cref{lem:null-bordant-ZxZ2}, the equivariant intersection form is even. Also taking connected sum with a simply-connected spin manifold in order to change the signature does not affect whether the intersection form is even.  It follows that $\lambda_M$ is even if and only if $\msec(M)=0$.  This covers the cases~\eqref{item:ZxZ2-3} and~\eqref{item:ZxZ2-4} of \cref{thm:ZxZ2}.
That is, we have shown that the parity of the intersection form distinguishes between $\beta(M) \in \{\pm 2\}$ and $\beta(M) \in \{0,4\}$.

\section{The case \texorpdfstring{$\beta(M)\in\{0,4\}$}{of Brown invariant 0 mod 4}}

It remains to relate the $\tau$ invariant in \eqref{item:ZxZ2-5} to the $\ter$ invariant from the Atiyah-Hirzebruch spectral sequence.	
It will be useful to identify $I^*$ with $I^w$, where $I^w$ is the submodule of $\Z\pi$ generated by $1-t$ and $1+T$.

\begin{lemma}
	The following map is an isomorphism
	  \begin{align*}
	    \theta_{I^w} \colon I^w & \xrightarrow{\cong} \Hom_{\Z\pi}(I,\Z\pi) \\
	    1-t &\mapsto \Bigg\{\begin{array}{rcl}
	      1-t &\mapsto & 1-t  \\
	      1-T &\mapsto& 1-T
	    \end{array}\Bigg\}\\
	    1+T & \mapsto
	    \Bigg\{\begin{array}{rcl}
	      1-t &\mapsto & 1+T  \\
	      1-T &\mapsto& 0
	    \end{array}\Bigg\}.
	  \end{align*}
\end{lemma}
\begin{proof}
	First note that the composition
		\[I^w\xrightarrow{\theta_{I^w}}\Hom_{\Z\pi}(I,\Z\pi)\xrightarrow{\ev_{1-t}}\Z\pi, \]
where $\ev_{1-t}(f) := f(1-t)$,  is the standard inclusion and thus $\theta_{I^w}$ is injective.
		
		Let $f\in \Hom_{\Z\pi}(I,\Z\pi)$. Then $f$ is determined by $f(1-t)$ and $f(1-T)$. We have $(1-t)f(1-T)=f((1-t)(1-T))=(1-T)f(1-t)$. As multiplication with $(1-t)$ is injective as a $\Z\pi$-module homomorphism $\Z\pi \to \Z\pi$, $f(1-T)$ is determined by $f(1-t)$. Hence $\Hom_{\Z\pi}(I,\Z\pi)\xrightarrow{\ev_{1-t}}\Z\pi$ is injective. Thus to show that $\theta_{I^w}$ is surjective, it suffices to show that the image of $\ev_{1-t}$ is contained in $I^w$.  For then we have a factorisation of $\ev_{1-t} \circ \theta_{I^w}$ as
\[I^w \xrightarrow{\theta_{I^w}} \Hom_{\Z\pi}(I,\Z\pi) \xrightarrow{\ev_{1-t}} I^w \subseteq \Z\pi,\]
and the composition $I^w \to I^w$ is the identity map.  It follows that $\ev_{1-t}$ is onto $I^w$, so is an isomorphism onto $I^w$, and therefore $\theta_{I^w}$ is an isomorphism.
		
	To show that image of $\ev_{1-t}$ is contained in $I^w$, let $\Z^w$ denote $\Z$ considered as a $\Z\pi$-module, with the $\Z\pi$-action given as follows: $t$ acts trivially and $T$ acts by multiplication by $-1$. Let $\varepsilon^w\colon\Z\pi\to\Z^w$ be the map sending $t$ to $1$ and $T$ to $-1$.
		 In $\Z^w$ we have
		\[2\varepsilon^w(f(1-t))=\varepsilon^w((1-T)f(1-t))=\varepsilon^w((1-t)f(1-T))=0\]
		and thus $f(1-t)\in I^w$ as needed.
		Note that we have shown that the evaluation map $\ev_{1-t}\colon\Hom_{\Z\pi}(I,\Z\pi) \to I^w$ is an inverse of $\theta_{I^w}$.
%
%
\end{proof}

If the secondary invariant of $M$ is trivial, then by \cref{remark-section-as-required} there are $k,n\in\mathbb{N}$ and an isomorphism $\psi\colon I\oplus I^w\oplus \Z\pi^n\to \pi_2(M\#^k(S^2\times S^2))$, such that $\lambda_{M\#^k(S^2\times S^2)}$ restricted to $I^w$ along $\psi$ is trivial.

\begin{lemma}\label{lemma:tau-comp-ZxZ2}
	We have that $\mu(\psi(0,1+T,0))=0$ and $\psi(0,1+T,0)$ is $\RPT$-characteristic. Furthermore, $\tau(\psi(0,1+T,0))=0 \in \Z/2$ if and only if $\ter(M)=0$. In particular, $\tau(\psi(0,1+T,0))$ is independent of the choice of $\psi$.
\end{lemma}

This will complete the proof of \eqref{item:ZxZ2-5} of \cref{thm:ZxZ2}, using that $\phi$ from \cref{thm:ZxZ2} agrees with $\theta_{I^w}(1+T)$.

We need some lemmas on Tate cohomology before starting the proof of \cref{lemma:tau-comp-ZxZ2}.
Let $\Z/2$ act on a $\Z\pi$-module $P$ by multiplication by $T$. The Tate cohomology groups are $2$-periodic and
\begin{align*}
	\wh H^0(\Z/2;P) &= \ker(1-T)/\im(1+T)\\
	\wh H^1(\Z/2;P) &= \ker(1+T)/\im(1-T).
\end{align*}
The next lemma follows from the definition.

\begin{lemma}\label{lem:tate-cohomology-I-Iw}
	The Tate cohomology has the following properties:
	\begin{enumerate}[(i)]
		\item $\wh H^{2i}(\Z/2;I^w) \cong \wh H^{2i+1}(\Z/2;I) \cong \Z/2$ generated by $1+T$ and $1-T$ respectively;
		\item $\wh H^{2i}(\Z/2;\Z\pi) = \wh H^{2i+1}(\Z/2;\Z\pi) = \wh H^{2i}(\Z/2;I) = \wh H^{2i+1}(\Z/2;I^w) =0$; and
		\item $\wh H^{2i}(\Z/2;P \oplus P') \cong \wh H^{2i}(\Z/2;P) \oplus \wh H^{2i}(\Z/2;P')$.
	\end{enumerate}
\end{lemma}

\begin{lemma}
	\label{lem:lift}
	Let $X$ be a CW model for $B\pi$ with 2-skeleton $K$. Let $i \colon K \to X$ be the inclusion map. Then any element $\alpha\in H^2(K;\Z\pi)$ generating $\wh H^0(\Z/2;H^2(K;\Z\pi))\cong \Z/2$ maps to $i^*(Tt)\in H^2(K;\Z/2)$ under the reduction map.
\end{lemma}

\begin{proof}
	We first prove the statement in the case where $X=S^1 \times \RP^\infty$ with the standard CW-structure and a single choice of generator $\alpha$ in $H^2(K;\Z\pi)$, where $K$ is the $2$-skeleton of $X$. The cellular cochain complex of $K$ is given by:
	\[\xymatrix @R-0.7cm @C+1cm {\Z\pi \ar[r]^{1-T} \ar[ddr]_{1-t^{-1}} & \Z\pi  \ar[r]^{1+T} \ar[ddr]^{1-t^{-1}} & \Z \pi \\
		& \oplus & \oplus   \\
		&  \Z\pi \ar[r]_{1-T}  & \Z\pi.
	}\]
	The generator $A$ of bottom right $\Z\pi$ summand, in $C^2(K;\Z\pi)$, maps to $i^*(Tt)$ under reduction to $\Z/2$. We need to show that $A$ can be taken as our choice of $\alpha$, that is $A$ generates $\wh H^0(\Z/2;H^2(K;\Z\pi)) \cong \wh H^0(\Z/2;I^w) \cong \Z/2$.
	We can extend $C^*(K;\Z\pi)$ to a resolution for $\Z^w$ as follows:
	\[\xymatrix @R-0.7cm @C+1cm {\Z\pi \ar[r]^{1-T} \ar[ddr]_{1-t^{-1}} & \Z\pi  \ar[r]^{1+T} \ar[ddr]^{1-t^{-1}} & \Z \pi \ar[ddr]^{1-t^{-1}} & &  \\
		& \oplus & \oplus & &   \\
		&  \Z\pi \ar[r]_{1-T}  & \Z\pi \ar[r]_{1+T} & \Z\pi \ar[r]^{\varepsilon^w} & \Z^w.
	}\]
	Write $\delta = (1-t^{-1},1+T) \colon \Z\pi^{\oplus 2} \to \Z\pi$. The identification
	\[H^2(K;\Z\pi) \cong C^2(K;\Z\pi)/\im d^2 \cong C^2/ \ker \delta \cong \im \delta \cong \ker \varepsilon^w \cong I^w\]
	sends $A$ to $1+T$, which generates $\wh H^{0}(\Z/2;I^w)$ by \cref{lem:tate-cohomology-I-Iw}. This completes the proof for the fixed CW model and our choice  $\alpha$.
	
	For the general case note that any two generators of $\wh H^0(\Z/2;H^2(K;\Z\pi))$ differ by $(1+T)x$ for some $x\in H^2(K;\Z\pi)$. Under the reduction map $H^2(K;\Z\pi) \to H^2(K;\Z/2)$ any element of the form $(1+T)x$ is mapped to zero. Hence the statement is true for all generators if and only if it is true for one given generator. Since any two 2-skeleta $K$, $K'$ become homotopy equivalent after taking a wedge product with sufficiently many copies of $S^2$, the general case follows from the case $X=S^1 \times \RP^\infty$.
\end{proof}

\begin{proof}[Proof of \cref{lemma:tau-comp-ZxZ2}]
	Let $N:=M\#^k(S^2\times S^2)$ and let
	\[\psi\colon I\oplus I^w\oplus {\Z\pi}^{\oplus n}\to \pi_2(N)\]
	be an isomorphism such that $\lambda_N$ restricted to $I^w$ along $\psi$ is trivial as assumed for the statement of \cref{lemma:tau-comp-ZxZ2}.
	Consider the extension
	\begin{equation}\label{eqn:extension-ZxZ2-tau-pf}
		0\to \ker d_2^N\xrightarrow{(i_1,i_2)} C_2^N\oplus\pi_2(N)\xrightarrow{\rho_1+\rho_2} \coker d_3^N\to 0.
	\end{equation}
	Let $F, F'$ be free modules such that there is an isomorphism
	\[\theta\colon \coker d_3^N\oplus F\to I^w\oplus F'.\]
	Let $\iota\colon \coker d_3^N\to \coker d_3^N\oplus F$ and $j\colon I^w\to I\oplus I^w\oplus \Z\pi^{\oplus 2n}$ be the inclusions and let $p\colon I^w\oplus F'\to I^w$ be the projection. Consider the composition
	\begin{multline*}
		s_2:=\psi\circ j\circ p\circ \theta\circ \iota\colon \coker d_3^N \xrightarrow{\iota} \coker d_3^N \oplus F \xrightarrow{\theta} I^w\oplus F' \xrightarrow{p} I^w \\  \xrightarrow{j} I\oplus I^w\oplus \Z\pi^{\oplus 2n} \xrightarrow{\psi} \pi_2(N).
	\end{multline*}
	We want to define a map $s_1\colon \coker d_3^N\to C_2^N$ such that $(s_1,s_2)\colon \coker d_3^N\to C_3^N\oplus\pi_2(N)$ is a section of the extension~\eqref{eqn:extension-ZxZ2-tau-pf}.
	
	It will be helpful to consider the extensions in the next diagram, which we may identify using the isomorphism $\theta \colon \coker d_3^N\oplus F\to I^w\oplus F'.$
	\[
	\xymatrix @C+0.6cm {
		0\ar[r]  & \ker d_2^N\ar[r] \ar[d]^= & C_2^N\oplus\pi_2(N)\oplus F\ar[r]^-{(\rho_1+\rho_2)\oplus \id_F} \ar[d]^=  & \coker d_3^N\oplus F\ar[r] \ar[d]^{\theta}_{\cong}  & 0  \\
		0\ar[r]  & \ker d_2^N\ar[r]  & C_2^N\oplus\pi_2(N)\oplus F\ar[r]^-{\rho'_1+\rho'_2+\rho'_3} & I^w\oplus F'\ar[r]  & 0
	} \]
	%
	%
	Here $\rho_k'$, for $k \in \{1,2,3\}$ denote the resulting maps defined by the diagram,
	with $\rho_1' \colon C_2^N \to I^w \oplus F'$, $\rho_2' \colon \pi_2(N) \to I^w \oplus F'$, and $\rho_3' \colon F \to I^w \oplus F'$.
	Now define a map $s_2'$ by
	\[s_2':=\psi\circ j\circ p\colon I^w\oplus F'\xrightarrow{p} I^w \xrightarrow{j} I\oplus I^w\oplus \Z\pi^{\oplus 2n} \xrightarrow{\psi}  \pi_2(N).\]
	The maps $s_2$ and $s_2'$ fit into the following commutative diagram:
	\[\xymatrix{\coker d_3^N \ar[r]^-{\iota} \ar[d]_-{s_2} & \coker d_3^N \oplus F \ar[d]^{\theta}_{\cong} \\
		\pi_2(N) & I^w \oplus F'  \ar[l]_-{s_2'} \ar[d]^-{p} \\
		I \oplus I^w \oplus \Z\pi^{\oplus 2n} \ar[u]^-{\psi}_{\cong} & I^w \ar[l]_-{j}}\]
	
	As above let $\Z/2$ act on $I\oplus I^w\oplus {\Z\pi}^{\oplus 2n}$ by multiplication with $T$. Then by \cref{lem:tate-cohomology-I-Iw} we have:
	\begin{multline*}
		\wh H^0(\Z/2;I\oplus I^w\oplus {\Z\pi}^{\oplus 2n}) = \ker(1-T)/\im(1+T) \\ \cong \wh H^0(\Z/2;I)\oplus \wh H^0(\Z/2;I^w)\oplus \wh H^0(\Z/2;{\Z\pi}^{\oplus 2n})\cong 0\oplus \Z/2\oplus 0.
	\end{multline*}
	Similarly,
	\begin{multline*}
		\wh H^0(\Z/2;\coker d_3^N)\cong \wh H^0(\Z/2;\coker d_3^N\oplus F) \cong \wh H^0(\Z/2;I^w\oplus F') \\ \cong \wh H^0(\Z/2;I^w) \oplus \wh H^0(\Z/2; F') \cong \Z/2 \oplus 0.
	\end{multline*}
	Since $\ker d_2^N$ is stably isomorphic to $I$, we have $\wh H^0(\Z/2;\ker d_2^N)=0$. Thus the map $C_2^N\oplus\pi_2(N)\to \coker d_3^N$ induces an isomorphism
	\[ \Z/2 \cong \wh H^0(\Z/2;\pi_2(N)) \to \wh H^0(\Z/2;\coker d_3^N) \cong \Z/2.\]
	These groups are generated by
	$\psi(0,1+T,0)\in \pi_2(N)$
	and $\theta^{-1}(1+T,0)\in I^w\oplus F'$ respectively. It follows that $\rho'_2\circ s_2'(1+T,0) \in I^w \oplus F'$ differs from $(1+T,0)$ by an element of the form $(1+T)(x_1,x_2)$ for some $(x_1,x_2)  \in I^w \oplus F'$.  That is:
	\begin{equation}\label{eqn:rho2s2eqn}
		\rho'_2\circ s_2'(1+T,0)=(1+T,0)+(1+T)(x_1,x_2) \in I^w \oplus F'.
	\end{equation}
	Multiply both sides by $(1-t)$ and apply linearity of $\rho'_2\circ s_2'$ to deduce that
	\[(1+T)\rho'_2\circ s_2'(1-t,0)=(1+T)(1-t,0)+(1+T)(1-t)(x_1,x_2)  \in I^w \oplus F'.\]
	Therefore $\rho'_2\circ s_2'(1-t,0) - (1-t,0) - (1-t)(x_1,x_2) \in \ker (1+T)$. Combining this with $\wh H^1(\Z/2;I^w\oplus F')=\ker (1+T)/\im(1-T)=0$ from \cref{lem:tate-cohomology-I-Iw}, it follows that
	\begin{equation}\label{eqn:rho2s2eqn-two}
		\rho'_2\circ s_2'(1-t,0)=(1-t,0)+(1-t)(x_1,x_2)+(1-T)(x_1',x_2')  \in I^w \oplus F'
	\end{equation}
	for some $(x'_1,x_2')\in I^w\oplus F'$.
	
	Let $q \colon C_2^N \to \coker d_3^N$ be the quotient map. The diagram
	\[\xymatrix{C_2^N \oplus F \ar[rr]^-{\rho_1' + \rho_3'} \ar[rd]_-{q \oplus \id} & & I^w \oplus F'  \\  & \coker d_3^N \oplus F \ar[ur]^-{\theta}_-{\cong} & }\]
	commutes. Therefore $\rho_1' + \rho_3' \colon C_2^N\oplus F \xrightarrow{} I^w\oplus F'$ is surjective, so we can pick preimages $(y_1,y_2),(y_1',y_2')\in C_2^N\oplus F$ of $(x_1,x_2),(x_1',x_2')$, respectively, so that in $I^w \oplus F'$ we have:
	\begin{equation}\label{eqn:x1x2y1y2-lifts}
		\rho_1'(y_1) + \rho_3'(y_2) = (x_1,x_2) \text{ and } \rho_1'(y_1') + \rho_3'(y_2') = (x_1',x_2').
	\end{equation}
	Also pick a lift $r\colon F'\to C_2^N\oplus F$ of the inclusion $i_2 \colon F'\to I^w\oplus F'$, satisfying $(\rho_1' + \rho_3') \circ r = i_2$.
	
	Define
	\begin{equation}\label{eqn:defn-of-s1-prime}
		\begin{aligned}
			s'_1 \colon I^w\oplus F' & \to C_2^N \oplus F \\
			(1+T,0) &\mapsto -(1+T)(y_1,y_2) \\
			(1-t,0) &\mapsto -(1-t)(y_1,y_2)-(1-T)(y_1',y_2')\\
			(0,x) & \mapsto r(x).
		\end{aligned}
	\end{equation}
	Then it follows from a straightforward computation using \eqref{eqn:rho2s2eqn}, \eqref{eqn:rho2s2eqn-two}, \eqref{eqn:x1x2y1y2-lifts}, and \eqref{eqn:defn-of-s1-prime} that  $(s_1',s_2')\colon I^w\oplus F'\to C_2^N\oplus \pi_2(N)\oplus F$ is a section of $\rho'_1+\rho'_2+\rho'_3$. Let \[s_1\colon \coker d_3^N\to C_2^N\] be $s_1'\circ \theta\circ \iota$ composed with the projection $C_2^N\oplus F\to C_2^N$. Then we claim that \[(s_1,s_2)\colon \coker d_3\to C_2^N\oplus \pi_2(N)\] is a section of $\rho_1+\rho_2$ and the restriction of $\lambda_N$ to the image of $s_2 \colon \coker d_3^N \to \pi_2(N)$ is trivial.
	
	To see that $(s_1,s_2)$ is a section of $\rho_1+\rho_2$ we consider the following diagram of maps.
	\[\xymatrix @R+1cm @C+1cm {0 \ar[r] & \coker d_3^N \ar@/^/[r]^{\iota} \ar@/_/[d]_-{(s_1,s_2)} & \coker d_3^N \oplus F \ar[r]^{\theta}_{\cong} \ar@/^/[l]^{p'} & I^w \oplus F' \ar@/_/[dl]_{(s_1',s_2')} \\
		0 \ar[r] & C_2^N \oplus \pi_2(N) \ar[r]_{\iota'} \ar@/_/[u]_{\rho_1 + \rho_2} & C_2^N \oplus \pi_2(N) \oplus F \ar[u]^{(\rho_1+\rho_2)\oplus \Id_F} \ar@/_/[ur]_{\rho_1' + \rho_2' + \rho_3'}
	}\]
	Here $\iota'$ is the inclusion map and $p'$ is the projection. Ignoring $(s_1,s_2)$, $(s_1',s_2')$, and $p'$, the diagram commutes. It also commutes if we ignore $\rho_1 + \rho_2$, $(\rho_1+\rho_2)\oplus \Id_F$, $\rho_1' + \rho_2' + \rho_3'$, and $p'$.    In addition $p' \circ \iota = \Id$ and $(\rho_1' + \rho_2' + \rho_3') \circ (s_1',s_2') = \Id$.  Therefore we have:
	\[\Id = (\rho_1' + \rho_2' + \rho_3') \circ (s_1',s_2') = \theta \circ ((\rho_1+\rho_2)\oplus \Id_F) \circ (s_1',s_2')\]
	where the second equation uses the right hand triangle. Conjugating with $\theta$ yields:
	\begin{equation}\label{eqn:functional-equation-for-splitting}
		\Id =   ((\rho_1+\rho_2)\oplus \Id_F) \circ (s_1',s_2')\circ \theta \colon \coker d_3^N \oplus F \to  \coker d_3^N \oplus F.
	\end{equation}
	Thus
	\begin{align*}
		(\rho_1 + \rho_2) \circ (s_1,s_2) &=  p' \circ \iota \circ (\rho_1 + \rho_2) \circ (s_1,s_2)  \\ &= p' \circ ((\rho_1+\rho_2)\oplus \Id_F) \circ \iota' \circ (s_1,s_2) \\
		&= p' \circ ((\rho_1+\rho_2)\oplus \Id_F) \circ  (s_1',s_2') \circ \theta \circ \iota \\ &= p' \circ \iota = \Id.
	\end{align*}
	We used $p' \circ \iota = \Id$, the left hand square, commutativity of the trapezium involving $(s_1,s_2)$ and $(s_1',s_2')$, \eqref{eqn:functional-equation-for-splitting}, and finally $p' \circ \iota = \Id$ again. This completes the proof that $(s_1,s_2)$ is a section of $\rho_1+\rho_2$.  The image of $s_2 \colon \coker d_3^N \to \pi_2(N)$ is contained in the image  $\psi(I^w)$ on which $\lambda_N$ vanishes by assumption. It follows that $\lambda_N$ vanishes on $s_2(\coker d_3^N)$ as claimed.
	
	Now that we have found an appropriate section, we complete the proof of \cref{lemma:tau-comp-ZxZ2}.
	Let $\alpha\in \coker d^N_3$ be an element generating $\wh H^0(\Z/2;\coker d_3^N)$. Then $s_2(\alpha)\in \pi_2(M)$ has trivial self-intersection number and is $\RPT$-characteristic by \cref{lem:ter}. Hence we can compute $\tau(s_2(\alpha))\in \Z/2$.
	The tertiary property holds for $\pi = \Z\times\Z/2$, and the group $E^3_{2,2} \cong H_2(\pi;\Z/2)/d_{2} = \Z/2$ in which the tertiary invariant lives is generated by $(Tt)^*$.  Thus using \cref{lem:lift} and the construction of $\tau_{M,s}$, we see that $\ter(M)=0$ if and only if $\tau_{M,s}(Tt)=\tau(s_2(\alpha))=0$.
	
	Note that $\theta\circ \iota(\alpha)=(1+T,0)+(1+T)(x,y)$ for some $(x,y)\in I^w\oplus F'$, and thus $s_2(\alpha)=\psi(0,1+T,0)+(1+T)\psi(0,x,0)$. Since $\lambda_M$ vanishes on the image of $\psi|_{I^w}$, $\tau(s_2(\alpha))=\tau(\psi(0,1+T,0))$ by \cref{lem:augtau}. Since the $\tau$-invariant on $s_2(\alpha)$ is well-defined in $\Z/2$, it follows that also $\psi(0,1+T,0)$ has to be $\RPT$-characteristic.  By the previous paragraph we deduce that $\tau(\psi(0,1+T,0)) = 0 \in \Z/2$ if and only if $\ter(M)=0$.
\end{proof}



\backmatter
\bibliographystyle{amsalpha}
\bibliography{classification}

\end{document}